\documentstyle[10pt]{article}
\input amssymb.sty
\input amssym.def
\input amssym
\input epsf

\newtheorem{theorem}{Theorem}[section]

\newtheorem{lemma}[theorem]{Lemma}
\newtheorem{proposition}[theorem]{Proposition}
\newtheorem{corollary}[theorem]{Corollary}
\newtheorem{conjecture}[theorem]{Conjecture}
\newtheorem{definition}[theorem]{Definition}

\newtheorem{construction}[theorem]{Construction}

\begin{document}
\newcommand{\Z}{{\Bbb Z}}
\newcommand{\R}{{\Bbb R}}
\newcommand{\Q}{{\Bbb Q}}
\newcommand{\C}{{\Bbb C}}
\newcommand{\lms}{\longmapsto}
\newcommand{\lra}{\longrightarrow}
\newcommand{\hra}{\hookrightarrow}
\newcommand{\ra}{\rightarrow}

\begin{titlepage}
\title{Periods and mixed motives }
\author{A. B. Goncharov}
\end{titlepage}
\date{}

\maketitle
\addtocounter{page}{+1}

\tableofcontents

\section  {Introduction} 

{\bf Abstract}. We define  motivic  
multiple polylogarithms and prove the 
double shuffle relations for them. We  use
this to study the motivic fundamental group $\pi_1^{\cal M}({\Bbb G}_m - \mu_N)$, 
where $\mu_N$ is the group of all $N$-th roots of unity,  
and relate its structure 
to the  geometry and topology of modular varieties
$$
Y_1(m; N):= 
\Gamma_1(m;N)\backslash GL_m(\R)/\R^*_+\cdot O_m \qquad \mbox{for $m=1,2,3,4,...$}.
$$
  Using this we get several new results about the action of the Galois group on $\pi_1^{(l)}({\Bbb G}_m - \mu_N)$, and 
the values of multiple polylogarithms at $N$-th roots of unity.  

To prove these results some new tools are developed, including the following: 

1)  We construct, under certain assumptions, framed mixed motives corresponding to 
periods. Here is the main  example. Suppose that $X$ is an $n$-dimensional 
 complex variety,  
$A$ and $B$ are two divisors on $X$,  $\omega_A \in \Omega^n_{\rm log}(X-A)$ and 
$\Delta_B$ is a relative cycle providing a class in $H_n(X, B)$. 
Consider the integral  
\begin{equation} \label{1.31.02.5}
\int_{\Delta_B}\omega_A 
\end{equation}  
If it is convergent we construct, under some technical conditions on $A$ and $B$, 
a framed mixed motive $m(X; [\omega_A]; [\Delta_B])$ 
whose period is given by integral (\ref{1.31.02.5}). 

2) We prove a specialization theorem: if $(\omega_{A(\varepsilon)}, 
\Delta_{B(\varepsilon)})$ 
is a perturbation of the data defining (\ref{1.31.02.5}) then
 in the Hodge/l-adic realization, under certain assumptions on 
$(A(\varepsilon), B(\varepsilon))$,  the specialization of 
 $m(X; [\omega_{A(\varepsilon)}]; [\Delta_{B(\varepsilon)}])$ at $\varepsilon = 0$
is equivalent (as a framed object) 
to $m(X; [\omega_A]; [\Delta_B])$. 

These results allow us to construct the  
framed mixed Tate motives corresponding to 
multiple polylogarithms and study the corresponding Hopf algebra. 

In general a mixed motive corresponding to a period 
is not uniquely defined. However  the {\it equivalence class 
of  framed mixed motive} is supposed to be uniquely defined by a period. 
The specialization theorem is used to prove 
that different constructions lead to the same equivalence class of framed objects. 
This plays  a decisive role in the proof of the motivic double shuffle relations. 
 
An immediate consequence of these results is a motivic proof of the 
week version of Zagier's conjecture on $\zeta_F(n)$-- see section 10.13.
 The other 
proves were given in [BD2] and [DJ].  

This paper is the second half of [G7].

 {\bf 1. The double shuffle relations}. Recall ([G0-1]) that the 
multiple polylogarithm functions are defined for $|x_i| <1$ by 
the power series
\begin{equation} \label{1.31.02.1}
{\rm Li}_{n_1, ..., n_m}(x_1, ..., x_m):= \sum_{k_1 < ... < k_m }
\frac{x_1^{k_1} ... x_m^{k_m}}{k_1^{n_1} ... k_m^{n_m}}
\end{equation}
They can be written as iterated integrals as follows. Set (assuming $a_i \not = 0$)
\begin{equation} \label{5*zx}
{\rm I}_{n_1, ..., n_m}
(a_1, ..., a_m):=  \int_{0}^{1} \underbrace {\frac{dt}{t-a_1}  \circ 
... \circ 
\frac{dt}{t} }
_{n_{1} \quad  \mbox {\rm times}} \circ \quad
... \quad \circ \underbrace 
{\frac{dt}{t-a_m} \circ ... 
\circ \frac{dt}{t}}_{ n_{m} \quad  \mbox {\rm times}}
\end{equation}
Then  (see s. 12 of [G1] or theorem 2.2 in [G7])
\begin{equation} \label{1.16.02.1z}
{\rm Li}_{n_1, ..., n_m}(x_1, ..., x_m)= (-1)^m{\rm I}_{n_1, ..., n_m}
(a_1, ..., a_m)
\end{equation}
where
\begin{equation} \label{1.16.02.1zsd}
a_1:= (x_1 ... x_m)^{-1}, ..., a_m:= x^{-1}_m
\end{equation}

The multiple polylogarithms satisfy two sets of relations, 
called the double shuffle relations. The first, called below the ${\rm I}$-relations, 
are obtained by applying the shuffle product formula 
$$
\int_{\gamma}\omega_1 \circ ... \omega_p \cdot \int_{\gamma}\omega_{p+1} \circ ... \omega_{p+q} 
= 
\sum_{\sigma \in \Sigma_{p.q}}
\int_{\gamma}\omega_{\sigma(1)} \circ ... \omega_{\sigma(p+q)}
$$
(where $\Sigma_{p.q}$ 
is the set of  shuffles of 
 $\{1, ..., p\}$ and $\{p+1, ..., p+q\}$) to the iterated integrals (\ref{5*zx}). For example 
\begin{equation} \label{2.2.02.1}
{\rm I}_1(a_1) \cdot {\rm I}_1(a_2) = {\rm I}_{1,1}(a_1, a_2) + {\rm I}_{1,1}(a_2, a_1)
\end{equation}
The second set of relations, 
called below the ${\rm Li}$-relations, are obtained by multiplying the power series 
(\ref{1.31.02.1}). The simplest of them is this:
 \begin{equation} \label{2.1.02.4}
{\rm Li}_{1}(x) \cdot {\rm Li}_{1}(y) = {\rm Li}_{1,1}(x,y) + {\rm Li}_{1,1}(y,x) + 
{\rm Li}_{2}(xy)
\end{equation}

The integral (\ref{5*zx}) is divergent if $a_m =1, n_m=1$. We regularize it 
by replacing $\int_0^{1}$ to $\int_0^{1-\varepsilon}$ and taking the constant term of the asymptotic expansion in $\log (\varepsilon)$. The regularized integrals obviously obey 
the ${\rm I}$-shuffle relations. However to keep the natural form of the 
${\rm Li}$-shuffle relations one needs to use a different regularization,   
and then compare the regularizations. 
  For the multiple $\zeta$-values, i.e. 
when $x_i =1$, the double shuffle relations 
were considered by Zagier [Z1].  

If 
$|x_i| > 1$ the double shuffle relations are not even  well defined as relations between 
numbers because of multivaluedness of multiple polylogarithms. 
What is more important, these relations   are only an analytic reflection of deeper 
 properties of the corresponding 
 geometric  objects. Here is how we handle them.

{\bf 2. The story on the Hodge level}. Recall the commutative graded 
Hopf algebra ${\cal H}_{\bullet}$ 
of the equivalence classes of framed Hodge-Tate structures (see for instance 
chapter 3 of [G7]).  
In chapter 5 of [G7] we defined a  framed Hodge-Tate structure
\begin{equation} \label{2.1.02.1}
{\rm I}^{\cal H}_{n_1, ..., n_m}
(a_1, ..., a_m) \in {\cal H}_w;  \quad w:= n_1 + ... +n_m, \quad a_i \in \C^*
\end{equation}
whose period  is given by the iterated integral (\ref{1.16.02.1z}), 
regularized if needed. 
By lemma 6.6 in [G7]  
elements (\ref{2.1.02.1}) satisfy the ${\rm I}^{\cal H}$-shuffle relations. 

For {\it generic} parameters $x_i \in \C^*$ 
we define
$$
\widetilde {\rm Li}^{\cal H}_{n_1, ..., n_m}
(x_1, ..., x_m) := {\rm I}^{\cal H}_{n_1, ..., n_m}
(a_1, ..., a_m)  
$$
assuming that 
$x_i$'s and $a_i$'s are related by  (\ref{1.16.02.1zsd}). 
 If 
we use this definition  for  all  $x_i \in \C^*$,
 the ${\rm Li}^{\cal H}$-shuffle relations will not hold:   
for instance (\ref{2.2.02.1}) 
is not  consistent with (\ref{2.1.02.4}) for  $x_1 = x_2 =1$. 
Therefore we  use a different regularization, suggested in theorem 7.1 
in  [G3], see also s. 2.10 in   [G7]:

\begin{definition} \label{2.8.02.1x}
For arbitrary   $x_i \in \C^*$  set 
\begin{equation} \label{2.1.02.5}
{\rm Li}^{\cal H}_{n_1, ..., n_m}
(x_1, ..., x_m):= 
{\rm Sp}_{\partial/\partial \varepsilon}\Bigl(\widetilde 
{\rm Li}^{\cal H}_{n_1, ..., n_m}(x_1(1-\varepsilon ), ..., x_m(1-\varepsilon ))\Bigr)
\end{equation}
We set ${\rm Li}^{\cal H}_{n_1, ..., n_m}(x_1, ..., x_m) =0$ if $x_1 ... x_m =0$.  
\end{definition} 
Here the expression under the sign ${\rm Sp}_{\partial/\partial \varepsilon}$ is understood as a unipotent variation of Hodge-Tate structures on a small disc with coordinate $\varepsilon$,
 punctured at $\varepsilon =0$. We apply to it the Verdier specialization functor 
 ${\rm Sp}_{\varepsilon = 0}$, getting a unipotent variation over the punctured tangent space at 
$\varepsilon =0$, 
and then take its fiber at the tangent vector $\partial/\partial \varepsilon$. 
The functor ${\rm Sp}_{\partial/\partial \varepsilon}$ 
is the  composition of these two functors.

The specialization theorem proved in chapter 3 implies that, assuming (\ref{1.16.02.1zsd}),
$$
{\rm Li}^{\cal H}_{n_1, ..., n_m}
(x_1, ..., x_m) = {\rm I}^{\cal H}_{n_1, ..., n_m}
(a_1, ..., a_m) \quad \mbox{provided $n_m>1$ or $x_m \not =1$}
$$  
i.e. when the corresponding integral is convergent. 

\begin{theorem} \label{2.1.02.7}
 The elements ${\rm Li}^{\cal H}_{n_1, ..., n_m}(x_1, ..., x_m)$ satisfy the 
${\rm Li}^{\cal H}$-shuffle relations for  any  $x_i \in \C$. 
\end{theorem}

For the precise form of the ${\rm Li}$--shuffle relations 
 see theorem \ref{4.16.01.13saq}. 
It is proved in  chapters 7-9.  We first prove it  for generic parameters $x_i$. There are two 
 different methods how to do this. We explain  these methods 
in chapter 8. To demonstrate how they 
 work we 
give in chapter 8 two proofs of the  ${\rm Li}^{\cal H}$-shuffle relations 
for the double logarithm.  
In this case all 
essential problems  are already present, while  
the combinatorics is still very simple. In chapter 9 we employ  one of the methods to 
get theorem \ref{2.1.02.7} for generic parameters $x_i$. 
In the appendix we outline a  proof of the same result via the other method. 
Then using the specialization theorem we deduce from this 
that theorem \ref{2.1.02.7} is valid for all parameters $x_i$.

To get the double shuffle relations it remains to compare explicitly 
the ${\rm Li}^{\cal H}$- and ${\rm I}^{\cal H}$-elements 
when  $n_m=1$ and $x_m =1$, i.e. when the corresponding integrals are divergent.
For this we need   the following.

{\it A Hodge version of the asymptotic expansions 
of divergent integrals}.  The equivalence classes of framed 
 unipotent variations of Hodge-Tate structures 
on  ${\C}^*$ form a commutative  graded Hopf algebra, denoted  
$
{\cal A}^{\cal H}_{\bullet}({\Bbb A}^1 - \{0\})
$ (see ch. 3 [G7]). Let $\varepsilon$ be a natural coordinate on 
${\Bbb A}^1 -\{0\}$ and $\log^{\cal H}\varepsilon$  the canonical variation 
of framed Hodge-Tate structures on $\C^*$ with the  period function  $\log\varepsilon$. 
We can view $\log^{\cal H}\varepsilon$ as an element of 
${\cal A}^{\cal H}_{1}({\Bbb A}^1 - \{0\})$. According to 
 lemma \ref{1.13.02.2} there is an isomorphism of graded commutative algebras  
\begin{equation} \label{3456wq}
{\cal A}^{\cal H}_{\bullet}({\Bbb A}^1 - \{0\}) = {\cal H}_{\bullet} \otimes_{\Q} 
\Q[\log^{\cal H}\varepsilon]
\end{equation}

If $H_{\varepsilon}$ is a unipotent 
variation of framed Hodge-Tate structures on a punctured disc, denote by 
$[{\rm Sp}_{\varepsilon = 0}](H_{\varepsilon})$ the equivalence class 
of its specialization ${\rm Sp}_{\varepsilon = 0}(H_{\varepsilon})$. 
Then, employing the isomorphism (\ref{3456wq}), we get  an element 
$$
[{\rm Sp}_{\varepsilon = 0}](H_{\varepsilon}) = \sum_{i \geq 0} H_{(i)} \cdot 
(\log^{\cal H}\varepsilon)^i \in {\cal H}_{\bullet} \otimes_{\Q} 
\Q[\log^{\cal H}\varepsilon]
$$
 It is  the Hodge version of the asymptotic expansion 
of the period function. Its fiber over 
the tangent vector $\partial/\partial \varepsilon$ coincides with the constant 
term $H_{(0)}$.

{\it The two regularizations and the comparison theorem}. a) Denote by 
\begin{equation} \label{1.13.02.3u}
 {\rm I}^{{\cal H}}_{n_1, ..., n_{m}}(a_1: ...: a_{m-1}: 
a_{m}: (1-\varepsilon))
\end{equation} 
the framed Hodge-Tate structure corresponding to the iterated 
integral (\ref{5*zx}) where $\int_0^1$ changed to $\int_0^{1-\varepsilon}$. 
It has been defined in chapter 5 of [G7]. Consider the 
unipotent variation of framed Hodge-Tate structures over a small 
punctured  disc with  coordinate 
$\varepsilon$, whose fiber at $\varepsilon$ is (\ref{1.13.02.3u}). 
Let  
\begin{equation} \label{1.13.02.3uz}
\widehat {\rm I}^{{\cal H}}_{n_1, ..., n_{m}}(a_1, ..., a_{m})
\in {\cal A}^{\cal H}_{\bullet}({\Bbb A}^1 - \{0\})  
\end{equation}
be the equivalence class of the  framed 
unipotent variation of  Hodge-Tate structures on 
the punctured tangent space obtained by specialization of the family 
(\ref{1.13.02.3u}). 
Its fiber over $\partial/\partial \varepsilon$  is  
the framed Hodge-Tate structure (\ref{2.1.02.1}), 
by the very definition of (\ref{2.1.02.1}). 

b) Set 
\begin{equation} \label{2.1.02.5x}
\widehat {\rm Li}^{\cal H}_{n_1, ..., n_m}
(x_1, ..., x_m):= 
\end{equation}
$$
[{\rm Sp}_{\varepsilon=0}]\Bigl({\rm Li}^{\cal H}_{n_1, ..., n_m}(x_1(1-\varepsilon ), 
..., x_m(1-\varepsilon ))\Bigr) \in {\cal A}^{\cal H}_{\bullet}({\Bbb A}^1 - \{0\})  
$$
The constant terms of the expansions (\ref{1.13.02.3uz})  and (\ref{2.1.02.5x}) 
are given by the elements 
(\ref{2.1.02.1}) and (\ref{2.1.02.5}) which we want to compare. 
An explicit formula for the coefficients 
of  expansion (\ref{1.13.02.3uz}) via multiple polylogarithms 
is given by lemma 6.7 or proposition 2.14 in [G7]. So to compute  
 ${\rm Li}^{{\cal H}}$ via ${\rm I}^{{\cal H}}$  
we  need to  compare 
 $\widehat {\rm Li}^{{\cal H}}$  with  
$\widehat {\rm I}^{{\cal H}}$. 

Set $\zeta^{\cal H}(n):= {\rm Li}^{\cal H}_n(1)$ and 
similarly with hats. For example  
$$
\widehat\zeta^{\cal H}(1) = \widehat {\rm Li}^{\cal H}(1) = {\rm Li}^{\cal H}(1-\varepsilon) =
-\log^{{\cal H}}(\varepsilon)
$$

\begin{definition} \label{1.13.02.10} ${\Bbb L}$ is the ${\cal H}_{\bullet} $-linear endomorphism
$$
{\Bbb L}:  {\cal H}_{\bullet} \otimes_{\Q} 
\Q[\log^{\cal H}\varepsilon] \lra {\cal H}_{\bullet} \otimes_{\Q} 
\Q[\log^{\cal H}\varepsilon] 
$$
 determined by  the following equality of formal power series in $u$: 
$$
\sum_{n \geq 0} {\Bbb L}\Bigl( (\log^{\cal H}\varepsilon)^n\Bigr) \cdot \frac{u^n}{n!}
:= {\rm exp}\left( - \sum_{n=1}^{\infty} 
(-1)^n \frac{\widehat \zeta^{\cal H}(n)}{n}u^n\right)
$$
\end{definition}

\begin{theorem} \label{1.13.02.5s} For any $x_i \in \C^*$ one has, assuming  (\ref{1.16.02.1zsd}), 
$$
\widehat {\rm Li}^{{\cal H}}_{n_1, ..., n_{m}}(x_1, ..., x_{m}) = {\Bbb L} \circ 
\widehat {\rm I}^{{\cal H}}_{n_1, ..., n_{m}}(a_1, ..., a_{m})
$$
\end{theorem} 

Our proof of theorem \ref{1.13.02.5s} follows the approach developed for the 
 multiple polylogarithms functions in  theorem 7.1 in [G3], see also s. 2.10 in [G7]. 

A similar comparison formula for the multiple $\zeta$-values was obtained 
by Zagier (it is documented in [IK]) and later independently 
by  Boute de Monville   
(and documented in [R]). 
However 
their approach  seems to have no Hodge/motivic interpretation.

Summarizing, combining lemma 6.6 from [G7]  giving  
the ${\rm I}^{\cal H}$-shuffle relations, 
theorem \ref{2.1.02.7}, 
and theorem \ref{1.13.02.5s} we get the double shuffle relations 
for the equivalence classes of 
framed Hodge-Tate structures corresponding to multiple polylogarithms. 

We suggest that when $x_i\in \mu_N$  these relations, 
combined with the distribution relations, should provide most of (but in general not all) 
 the relations 
between the multiple polylogarithm Hodge-Tate structures with $x_i \in \mu_N$. 
In a contrast with this, if  $x_i \in \C$  
they  provide only a little part of the relations.

{\bf 3. The motivic and \'etale versions of the story}. 
 Let $F$ be a number field. Then there exists an abelian category ${\cal M}_T(F)$ 
of of mixed Tate motives over $F$ (see chapter 5 of [G9] or [L1]). 
The objects of this category carry canonical weight filtration $W_{\bullet}$. 
 The category ${\cal M}_T(F)$  is equipped with   
canonical fiber functor 
$$
\omega: X \lra \oplus_n {\rm Hom}(\Q(-n), {\rm Gr}^W_{2n}X)
$$
This fiber functor  provides an equivalence 
of the category ${\cal M}_T(F)$ with the category of finite dimensional graded 
comodules over the 
fundamental Hopf algebra ${\cal A}_{\bullet}(F)$. 
Here ${\cal A}_{\bullet}(F)$ is a commutative Hopf algebra over $\Q$ 
graded by the integers $n \geq 0$. The elements of ${\cal A}_n(F)$ 
can be understood as  the equivalence classes of the framed objects 
in the category ${\cal M}_T(F)$, see ch. 3 of [G7].

Suppose that  $a_i \in F^*$ and $n_i$ are 
arbitrary positive integers. We define geometrically 
 framed mixed Tate motives over $F$
\begin{equation} \label{2.3.02.1}
{\rm I}^{\cal M}_{n_1, ..., n_m}(a_1, ..., a_m) \in {\cal A}_w(F); \quad w:= n_1 + ... +n_m
\end{equation}
and 
\begin{equation} \label{2.3.02.wq1}
{\rm Li}^{\cal M}_{n_1, ..., n_m}(x_1, ..., x_m) \in {\cal A}_w(F); 
\end{equation}
For any complex embedding $\sigma: F \hookrightarrow \C$ the Hodge realization 
of (\ref{2.3.02.1}) 
is given by ${\rm I}^{\cal H}_{n_1, ..., n_m}(\sigma(a_1), ..., \sigma(a_m))$, 
and similarly   for the elements (\ref{2.3.02.wq1}). 
We prove that the motivic multiple polylogarithms 
satisfy the double shuffle  relations. The proof of the ${\rm I}^{\cal H}$-relations 
given in [G7]  is motivic. For generic $x_i$ 
the proof of the ${\rm Li}$-relations given in chapter 9 is also motivic. 
However in some special cases the result  is deduced from the Hodge 
realization using the injectivity of regulators. 

Theorem \ref{1.13.02.5s} provides an explicit formula expressing 
${\rm Li}^{\cal H}$- via the ${\rm I}^{\cal H}$-elements. It 
implies, using the in jectivity of regulators, 
  the same kind of formula  relating the ${\rm Li}^{\cal M}$- and 
${\rm I}^{\cal M}$-elements.  

The ${\rm Li}$-shuffle relations are obvious for the power series (\ref{1.31.02.1}), but 
deliver  surprisingly nontrivial information on the motivic level. 
 For example  if $x_i$ are N-th roots of unity they  provide the most sofisticated information  
about the action of the motivic Galois group 
on the motivic fundamental group $\pi_1^{\cal M}({\Bbb G}_m - \mu_N)$.
In the l-adic realization they provide 
 constraints on the image of the absolute Galois group 
acting by automorphisms of 
$\pi_1^{(l)}({\Bbb G}_m - \mu_N)$. 

Suppose now that $F$ is {\it an arbitrary field such that $\mu_{l^{\infty}} \not \in F^*$}. 
Then we define framed l-adic mixed Tate ${\rm Gal}(\overline F/F)$-modules 
\begin{equation} \label{2.3.02.wq1gt}
{\rm Li}^{\rm et}_{n_1, ..., n_m}(x_1, ..., x_m) 
\end{equation}
related to multiple polylogarithms and 
prove the l-adic version of the double shuffle relations for them. 
When $x_i$ are elements of a number field $F$ the framed 
${\rm Gal}(\overline F/F)$-modules (\ref{2.3.02.wq1gt}) are 
given by the l-adic realization of their motivic counterparts 
(\ref{2.3.02.wq1})

{\bf 4. The higher cyclotomy and cyclotomic Hopf algebras}.  
Let   $\zeta_N := {\rm exp}(2\pi i/N)$. We define 
$
{\cal Z}_{w}^{\cal M}(\mu_N)
$ 
as the $\Q$-vector subspace of ${\cal A}_w(\Q(\zeta_N))$ generated 
by the elements 
\begin{equation} \label{2.2.02.5}
{\rm Li}^{{\cal M}}_{n_1, ..., n_{m}}(\zeta^{\alpha_1}_N, ..., \zeta^{\alpha_m}_N) \in {\cal A}_w(\Q(\zeta_N)), \quad w= n_1 + ... n_m
\end{equation}   and set 
$
{\cal Z}_{\bullet}^{\cal M}(\mu_N):= \oplus_{w \geq 0}{\cal Z}_{w}^{\cal M}(\mu_N)
$.

\begin{theorem} \label{4.16.01.s} ${\cal Z}_{\bullet}^{\cal M}(\mu_N)$ is a graded Hopf algebra.
\end{theorem} 

Recall 
$$
S_N := {\rm Spec}(\Z[\zeta_N][\frac{1}{N}])
$$
By the very definition 
 ${\cal Z}_{1}^{\cal M}(\mu_N)$ is generated by cyclotomic $N$-units 
$1 - \zeta_N^{\alpha}$, so 
\begin{equation} \label{1.20.02.5}
{\cal Z}_{1}^{\cal M}(\mu_N) \stackrel{\sim}{=}{\cal O}^*(S_N)\otimes \Q
\end{equation} 
We call ${\cal Z}_{\bullet}^{\cal M}(\mu_N)$ the {\it level 
$N$ motivic cyclotomic Hopf algebra}, and suggest its study 
 as the  goal for higher analog of the theory of cyclotomic units.

Recall the commutative Hopf algebra ${\cal A}_{\bullet}(S_N)$   classifying
the abelian category ${\cal M}_T(S_N)$ of 
mixed Tate motives over the scheme $S_N$ (see  chapter 3 of [G7]  or [DG]).

\begin{theorem} \label{4.16.01.qa} 
${\cal Z}_{\bullet}^{\cal M}(\mu_N)$ is a 
Hopf subalgebra of  ${\cal A}_{\bullet}(S_N)$.
\end{theorem}

Theorem \ref{4.16.01.qa} provides a strong bound from above on the size of the cyclotomic Hopf algebra 
${\cal Z}_{\bullet}^{\cal M}(\mu_N)$. Indeed, the Hopf algebra  ${\cal A}_{\bullet}(S_N)$ 
has the following structure. 
Consider the free graded Lie algebra generated by the vector spaces 
$K_{2n-1}(S_N)\otimes \Q$ sitting in the degree $-n$, when $n = 1, 2, 3, ... $. Then 
${\cal A}_{\bullet}(S_N)$ is 
isomorphic to the graded dual to the universal enveloping algebra 
of this graded Lie algebra. In particular the algebra structure of 
${\cal A}_{\bullet}(S_N)$ is described as follows. 

\begin{proposition} \label{5.26.02.1}
The algebra ${\cal A}_{\bullet}(S_N)$ is isomorphic to the tensor algebra 
generated by  the  graded $\Q$-vector space
\begin{equation} \label{2.14.02.3}
\oplus_{n=1}^{\infty}K_{2n-1}(S_N)\otimes \Q
\end{equation}
 and equipped with  commutative multiplication provided by the shuffle product.
\end{proposition}

By Borel's theorem  
$$
{\rm dim} K_{2w-1}(\Z) \otimes \Q= \left\{ \begin{array}{ll}
0 &  \quad \mbox{$w$: even} \\
 1 &    \quad \mbox{$w>1$: odd} \end{array} \right.
 $$
$$
{\rm dim} K_{2w-1}(S_2) \otimes \Q= \left\{ \begin{array}{ll}
0 &  \quad \mbox{$w$: even} \\
 1 &    \quad \mbox{$w$: odd} \end{array} \right.
 $$
and for $N>2$:
$$
{\rm dim} K_{2w-1}(S_N) \otimes \Q =  \left\{ \begin{array}{ll}
\frac{\varphi(N)}{2} &  \quad \mbox{$ w>1$} \\
\frac{\varphi(N)}{2}+ p(N)-1 &  \quad \mbox{$ w=1$} \end{array} \right.
 $$
where $p(N)$ is the number of prime factors of $N$. 

Recall ([G7]) the $\Q$-vector space 
$\widetilde {\cal Z}_w(\mu_N)$ spanned over $\Q$ 
by $(2 \pi i)^{-w}$ times the values of weight $w$ 
multiple 
polylogarithms at $N$-th roots of unity. (We take into account all  brunches of the corresponding 
multivalued analytic functions).  
Then $\widetilde {\cal Z}(\mu_N) := \cup \widetilde {\cal Z}_w(\mu_N)$ 
is a commutative algebra filtered by the weight.

\begin{theorem} \label{4.16.01.qoin} a) 
There exists canonical surjective algebra homomorphism 
\begin{equation} \label{5.26.02.2}
 {\cal Z}_{\bullet}^{\cal M}(\mu_N) \lra {\rm Gr}_{\bullet}^W\widetilde {\cal Z}(\mu_N)
\end{equation}

b) The algebra  ${\rm Gr}_{\bullet}^W\widetilde {\cal Z}(\mu_N)$ is a 
subquotient of the algebra ${\cal A}_{\bullet}(S_N)$.
\end{theorem}

The part b) follows immediately from the part a) and theorem \ref{4.16.01.qa}.

A considerable part of the proof of this theorem is contained in [G7]. 
In particular we constructed there a homomorphism (\ref{5.26.02.2}) 
with ${\cal Z}_{\bullet}^{\cal M}(\mu_N)$ replaced by ${\cal Z}_{\bullet}^{\cal H}(\mu_N)$. 
Combined with the given above explicit description of the algebra 
${\cal A}_{\bullet}(S_N)$, theorem  \ref{4.16.01.qoin} 
implies strong estimates from above for ${\rm dim}{\rm Gr}_{w}^W\widetilde{\cal Z}(S_N)$ 
 discussed in [G5]. In particular in the case $N=1$ 
this completes the proof of theorem 1.2 in [G5]. 
Another proof in the $N=1$ case  was suggested by T. Terasoma [T].

{\bf 5. The cyclotomic and dihedral Lie algebras}. 
The space of indecomposables
$$
{\cal C}_{\bullet}^{\cal M}(\mu_N):= 
\frac{{\cal Z}_{>0}^{\cal M}(\mu_N)}{{\cal Z}_{>0}^{\cal M}(\mu_N)\cdot {\cal Z}_{>0}^{\cal M}(\mu_N)}
$$
has a natural structure of a graded Lie coalgebra over $\Q$. 
The cyclotomic Lie algebra ${\rm C}_{\bullet}^{\cal M}(\mu_N)$ is defined as its  graded dual.  

There is a depth filtration $F^{D}_{\bullet}$ on the cyclotomic  Hopf algebra such that the 
depth $\leq d$ part
is generated by the elements (\ref{2.2.02.5})  
with  $m \leq d$. It induces the depth filtration on the cyclotomic Lie coalgebra and algebra. 
The depth filtration is compatible with the weight grading. 
Taking the  associate graded for the depth filtration we get  bigraded Lie coalgebra and Lie algebra 
over $\Q$:
$$
{\cal C}^{\cal M}_{\bullet, \bullet}(\mu_N):= {\rm Gr}^D\left({\cal C}_{\bullet}^{\cal M}(\mu_N)\right), \quad 
{\rm C}^{\cal M}_{\bullet, \bullet}(\mu_N):= {\rm Gr}^D\left({\rm C}_{\bullet}^{\cal M}(\mu_N)\right)
$$

Recall the dihedral Lie algebra ${\rm D}_{\bullet, \bullet}(G)$ of a commutative  group $G$. 
It was defined  [G3] or  chapter 4 of [G4] as the graded dual to the dihedral Lie 
coalgebra ${\cal D}_{\bullet, \bullet}(G)$. This Lie coalgebra 
  is generated by the symbols 
$$
{\rm I}_{n_1, ..., n_m}(g_1, ..., g_m), \quad g_i \in G, \quad n_i >0
$$ 
subject to certain relations. The relations  reflect the double shuffle relations 
projected 
onto the space of decomposables and  considered modulo the lower depth terms, 
and the distribution relations. 
Let us define a map 
$
  {\cal D}_{\bullet, \bullet}(\mu_N) \lra {\rm C}^{\cal M}_{\bullet, \bullet}(\mu_N)
$
by setting  ($a_i \in \mu_N$) 
\begin{equation} \label{1.21.02.1f}
{\rm I}_{n_1, ..., n_m}(a_1, ..., a_m) \lms 
\end{equation}
$${\rm I}_{n_1, ..., n_m}^{\cal M}(a_1, ..., a_m) 
\quad \mbox{modulo  depth $<m$ terms and products}
$$
where the right hand side is defined via (\ref{2.3.02.1}). 

\begin{theorem} \label{4.16.01.pi} Formula (\ref{1.21.02.1f}) provides a well defined  surjective 
homomorphism of bigraded Lie coalgebras
$$
\nu_{\bullet, \bullet}^{\cal M}(\mu_N): {\cal D}_{\bullet, \bullet}(\mu_N) \lra 
{\cal C}^{\cal M}_{\bullet, \bullet}(\mu_N) 
$$
\end{theorem} 

This theorem reflects the following two basic facts: 

a) we have  the motivic double shuffle and distribution 
relations, 

b) 
the cobracket in the dihedral Lie coalgebra is compatible 
with the coproduct of the elements (\ref{2.1.02.1}) explicitly computed in chapter 6 of [G7]. 
(In fact our definition of the cobracket in [G3] was suggested by that calculation).

\begin{theorem} \label{2.14.02.2} There are canonical isomorphisms
$$
{\cal D}_{w,1}^{\cal M}(\mu_N)    = K_{2w-1}(\Q(\zeta_N))\otimes \Q = {\cal C}_{w,1}^{\cal M}(\mu_N)
$$
\end{theorem}

For the first isomorphism see in s. 7.2 of [G4]. The second 
follows from the results of [BD2], or  easily deduced from 
the results of this paper, see section 10.12.

Combining theorem \ref{2.14.02.2} with theorem \ref{4.16.01.qa} and description of the 
fundamental Hopf algebra ${\cal A}_{\bullet}(S_N)$, 
we conclude that the Lie algebra ${\rm C}_{\bullet}^{\cal M}(\mu_N)$ is generated 
by its depth one component given by the graded vector space (\ref{2.14.02.3}). 

\begin{conjecture} \label{dwhn1dds}
The map $\nu_{\bullet, \bullet}^{\cal M}(\mu_N)$ is an isomorphism 
for either $N=1$, or $N=p$ is a prime and $w=m$.
\end{conjecture}

\begin{theorem} \label{2.14.02.1} Conjecture \ref{dwhn1dds} is true for $m \leq 3$. 
Thus the double shuffle relations provide all the relations between 
the elements ${\rm Li}^{\cal M}_{n_1, ..., n_m}(\zeta_p^{\alpha_1}, ...,  \zeta_p^{\alpha_m})$ 
for $m \leq 3$ if $N=1$, or if $N=p$ is a prime and $n_i =1$. 
\end{theorem}
This follows from the results of [G4] and chapter 10 below. 
The results of [G10] give a strong evidence for conjecture \ref{dwhn1dds} for $m=4$.

{\bf 6. The  cyclotomic Lie algebras and   geometry of modular varieties}. Theorem \ref{4.16.01.pi} is the crucial step on the way to a mysterious correspondence between 
 the structure of 
the cyclotomic Lie algebras and geometry and topology of modular varieties. 
We defer the detailed discussion of the geometric aspects of this  correspondence to [G10]. 
Several constructions have been described in [G3], [G4] 
and especially [G5]. Here is a very brief outlook. 

The rank $m$ modular complex $
{M}_{(m)}^{\bullet}$ is a complex of $GL_m(\Z)$-modules
$$
M_{(m)}^{1} \stackrel{\partial}{\lra} M_{(m)}^{2}
\stackrel{\partial}{\lra}
... \stackrel{\partial}{\lra} M_{(m)}^{m}
$$
which we  defined in [G3], see also s. 2.5 of [G4]. 
 If $m=2$ it is isomorphic to the chain complex of the classical modular triangulation 
of the hyperbolic plane.

Denote by
$\Lambda_{(m, w)}^*{\cal G}_{\bullet \bullet}$ the  bidegree $(m,w)$ part of the
standard cochain complex of a bigraded  Lie (co)algebra ${\cal G}_{\bullet \bullet}$.
Set $$
\widehat {\cal C}^{\cal M}_{\bullet, \bullet}(\mu_p):= 
{\cal C}^{\cal M}_{\bullet, \bullet}(\mu_p) \oplus \Q_{1,1}
$$
 where $\Q_{1,1}$ 
is the one dimensional Lie coalgebra of the  weight-depth  $(1,1)$. 

\begin{theorem} \label{dwhn1d}
For  $m>1$ there exists canonical surjective map of  complexes
$$
    \mu^*_{m;w}(N):  M^*_{(m)} \otimes_{\Gamma_1(m;N)} S^{w-m}{\rm V}_m  \lra
    \Lambda_{(m, w)}^*\widehat {\cal C}^{\cal M}_{\bullet, \bullet}(\mu_N)\,.
$$
\end{theorem}

This theorem is an immediate consequence of theorem \ref{4.16.01.pi} and the following 
result proved in [G3], see also [G10]:
for  $m>1$ there is canonical surjective map of  complexes
$$
    \eta^*_{m;w}(N):  M^*_{(m)} \otimes_{\Gamma_1(m;N)} S^{w-m}{\rm V}_m  \lra
    \Lambda_{(m, w)}^*\widehat {\mathcal D}_{\bullet \bullet}(\mu_N)
$$
It is an isomorphism when  $N=1$, or when $N=p$ is a prime and $w=m$.

\begin{conjecture} \label{dwhn1dd}
Let $N=1$, or $N=p$ is a prime and $w=m$. Then the map $\mu^*_{m;w}(N)$ is an isomorphism.
\end{conjecture}

Let $\Gamma_1(m; p) \subset GL_m(\Z)$ be the  subgroup stabilizing the vector 
$(0, ..., 0, 1)$ modulo $N$. Denote by $\varepsilon_m$ the  one dimensional 
$GL_m(\Z)$-module given by the determinant.  
The results above allow to replace ${\cal G}^{(l)}_{\bullet, \bullet}(\mu_N)$ 
by ${\rm C}^{\cal M}_{\bullet, \bullet}(\mu_N)$ everywhere  in  [G5] where 
 the structure of the Galois Lie algebra 
${\cal G}^{(l)}_{\bullet, \bullet}(\mu_N)$ was related  to 
 modular varieties. Here is a sample.  

\begin{theorem} \label{1.30.02.7} One has for a prime $p \geq 5$
$$
H^i_{(2,2)}(\widehat {\rm C}^{\cal M}_{\bullet, \bullet}(\mu_p))  =  
H^{i-1}(\Gamma_1(2;p), \varepsilon_2); \quad i = 1,2
$$
$$
H^i_{(3,3)}(\widehat {\rm C}^{\cal M}_{\bullet, \bullet}(\mu_p)) =  
 H^{i}(\Gamma_1(3;p), \Q); \quad i = 1,2,3
$$
\end{theorem}

{\bf 7. The motivic torsor of path 
${\cal P}^{\cal M}({\Bbb G}_m - \mu_N; v_0, v_1)$ and the cyclotomic Hopf algebra}. 
Let $t$ be the canonical coordinate on $P^1 - \{0,1,\infty\}$ and 
$v_0, v_1$ the corresponding tangential base points at $0$ and $1$. 
Denote by  $v_{\infty}$ the tangential base point at infinity 
provided by the local parameter $t^{-1}$. 
The motivic torsor of path  ${\cal P}^{\cal M}({\Bbb G}_m - \mu_N; v_0, v_1)$ 
is a pro-object in the abelian category 
of mixed Tate motives over $S_N$. It is defined in [DG]. Its more explicit geometric 
construction 
is given in section 6.3 below. Namely, consider 
the dual to the associate graded for the weight filtration 
\begin{equation} \label{2.6.02.2}
\Bigl({\rm Gr}^W{\cal P}^{\cal M}({\Bbb G}_m - \mu_N; v_0, v_1)\Bigr)^{\vee}
\end{equation}
as an Ind-object in the category of pure Tate motives. 
Recall that ${\cal A}_{\bullet}(S_N)$ is a 
Hopf algebra in the same category. 
In section 6.3 we define explicitly 
a coaction of 
the 
 Hopf  algebra ${\cal A}_{\bullet}(S_N)$ on (\ref{2.6.02.2}). Therefore 
we get an Ind-object in the category of mixed Tate motives over $S_N$. 
An immediate corollary of this construction is  the following important result.
 
\begin{theorem} \label{2.6.02.1} The Ind-object (\ref{2.6.02.2}) 
is a comodule over the cyclotomic Hopf subalgebra 
${\cal Z}^{\cal M}_{\bullet}(\mu_N) \subset {\cal A}_{\bullet}(S_N)$. 
Moreover ${\cal Z}^{\cal M}_{\bullet}(\mu_N)$
acts cofreely on (\ref{2.6.02.2}). 
\end{theorem}
 
The proof  uses the properties of the Hodge version of 
the motivic torsor of path established 
in theorem 1.6 in [G7].  

The grading of ${\cal Z}^{\cal M}_{\bullet}(\mu_N)$ provides an action of ${\Bbb G}_m$ 
on the prounipotent 
algebraic group scheme ${\rm Spec}{\cal Z}^{\cal M}_{\bullet}(\mu_N)$. 
Theorem \ref{2.6.02.1} just means that the image of the  motivic Galois group acting
on  ${\cal P}^{\cal M}({\Bbb G}_m - \mu_N; v_0, v_1)$ 
is given by the semidirect product of ${\Bbb G}_m$ and 
 ${\rm Spec}{\cal Z}^{\cal M}_{\bullet}(\mu_N)$.  

\begin{theorem} \label{2.6.02.wqw1} a) The motivic Galois group acts
on  ${\cal P}^{\cal M}({\Bbb G}_m - \mu_N; v_0, v_1)$ and on 
$\pi_1^{\cal M}({\Bbb G}_m - \mu_N; v_{\varepsilon})$, where 
$\varepsilon \in \{0,1, \infty\}$, 
via  the same quotient. 

b) This quotient is the semidirect product of 
${\Bbb G}_m$ and  ${\rm Spec}{\cal Z}^{\cal M}_{\bullet}(\mu_N)$. 
\end{theorem}

Observe that the map $z \lms z^{-1}$ provides an isomorphism 
$\pi_1^{\cal M}({\Bbb G}_m - \mu_N; v_0) \stackrel{\sim}{=}
\pi_1^{\cal M}({\Bbb G}_m - \mu_N; v_{\infty})$.
Therefore the part b) of this theorem follows from the part a) 
and theorem \ref{2.6.02.1}. 

{\bf 8. Applications to the Galois action on $\pi_1^{(l)}({\Bbb G}_m - \mu_N; 
v_{\infty})$}. 
Theorem \ref{2.6.02.wqw1} in the l-adic realization provides an important new 
information about the action of the absolute Galois group on the pro-l completion 
$\pi_1^{(l)}({\Bbb G}_m - \mu_N; v_{\infty})$
of the fundamental group. This combined  with the results outlined in s. 1.6 allows 
to relate the structure of the Galois module 
$\pi_1^{(l)}({\Bbb G}_m - \mu_N; v_{\infty})$ 
with the geometry of modular varieties. Here is a more detailed account. 

Let ${\Bbb L}^{(l)}_N$ be the l-adic pro-Lie algebra corresponding 
to $\pi_1^{(l)}({\Bbb G}_m - \mu_N, v_{\infty})$ via Maltsev's theory. The 
Galois group ${\rm Gal}(\overline \Q/\Q(\zeta_{l^{\infty}N}))$ acts by automorphisms 
of $\pi_1^{(l)}$ and hence of  ${\Bbb L}^{(l)}_N$, providing a  homomorphism
$$
{\rm Gal}(\overline \Q/
\Q(\zeta_{l^{\infty}N})) \lra {\rm Aut}{\Bbb L}^{(l)}_N
$$
Let ${\cal G}_N^{(l)}$ be the linearization of the image of this map:
$$
{\cal G}_N^{(l)}:= \lim_{\longleftarrow}{\rm Lie}\left({\rm Im}\left({\rm Gal}(\overline \Q/
\Q(\zeta_{l^{\infty}N})) \lra {\rm Aut}{\Bbb L}^{(l)}_N[m]\right)\right)
$$
where $L[m]$ stays for the quotient of a Lie algebra $L$ by the m-th lower central series 
ideal. Therefore 
$$
{\cal G}_N^{(l)} \subset {\rm Der} {\Bbb L}^{(l)}_N
$$
The natural weight filtration on ${\Bbb L}^{(l)}_N$ 
coincides with the lower central series filtration. It induces the weight 
filtration $W_{\bullet}$ on ${\rm Der} {\Bbb L}^{(l)}_N$, and hence on ${\cal G}_N^{(l)}$. 
The weight filtration admits a splitting compatible with the depth 
 filtration (see [G4]). There is also the 
 depth filtration on ${\Bbb L}^{(l)}_N$ and hence on ${\cal G}^{(l)}_N$.

The l-adic realization of $\pi_1^{\cal M}({\Bbb G}_m - \mu_N, v_{\infty})$ 
is given by the Galois module  ${\Bbb L}^{(l)}_N$. 

\begin{theorem} \label{5.2.02.2} a) There is canonical isomorphism of Lie algebras
$$
{\rm Gr}_{\bullet}^W{\cal G}^{(l)}_N = {\rm C}^{\cal M}_{\bullet}(\mu_N)\otimes \Q_l
$$

b) This isomorphism is compatible with the depth  filtration. 
\end{theorem}

It follows from theorem \ref{2.6.02.1}, the Tannakian formalism and general properties of the l-adic realization functor on the category ${\cal M}_T(F)$.

Let ${\cal G}_{\bullet, \bullet}^{(l)}(\mu_N)$ be the 
associate graded for the depth and weight filtrations. 
Theorem \ref{5.2.02.2} implies  an isomorphism 
$$
{\cal G}_{\bullet, \bullet}^{(l)}(\mu_N) = {\rm C}^{\cal M}_{\bullet, \bullet}(\mu_N) \otimes \Q_l
$$
Combining this with theorem \ref{4.16.01.pi} we get the following result 
 stated in conjecture 1.1 in [G4]: 

\begin{theorem} \label{5.2.02.3} There is canonical inclusion
$$
{\cal G}_{\bullet, \bullet}^{(l)}(\mu_N) \hookrightarrow 
{\rm D}_{\bullet, \bullet}(\mu_N)\otimes \Q_l
$$
\end{theorem}

Now we turn to more technical aspects of the paper.

{\bf 9.  Periods and framed mixed motives}. The integral (\ref{1.31.02.5}) is a typical example of 
 a period. Here are the simplest examples:
$$
\log (2) = \int_{1}^{2}\frac{dt}{t}; \quad  \zeta(2) = \int_{0 \leq t_1 \leq t_2 \leq 1}
\frac{dt_1}{1-t_1 }\wedge \frac{dt_2}{t_2}
$$
A period is supposed to come from the Hodge realization of 
a uniquely defined 
equivalence class of framed mixed motive. 

Let 
\begin{equation} \label{5.25.02.1}
\int_{\Delta_B} \omega_A
\end{equation}
be a convergent integral of a differential $n$-form 
$\omega_A$ on a complex variety $X$ with logarithmic singularities at a divisor $A$ 
over a relative cycle 
${\Delta_B}$ representing 
a class in $H_n(X, B)$, where $B$ is a divisor on $X$. 
Let us suppose that 
\begin{equation} \label{5.25.02.10}
\mbox{the divisors $A$ and $B$ have no common irreducible components}
\end{equation}
and the divisor $A \cup B$ is a 
normal crossing divisor. Then setting $B_A:= B - (B \cap A)$
  we define the corresponding 
mixed motive by  
\begin{equation} \label{5.2.02.1}
 H^n(X-A, B_A)
\end{equation} 
It has a natural framing coming from classes $[\omega_{A}]; [\Delta_{B}]$. 
Denote the corresponding framed mixed motive 
as well as its equivalence class by $m(X; [\omega_{A}]; [\Delta_{B}])$. 

In general construction of 
the corresponding framed object is neither unique nor
 obvious even the Hodge realization - 
see the next section for discussion of $\zeta(2)$. Here is a general method. 

Let $\pi: \widehat X \lra X$ be a blow up of $X$ transforming the 
divisor $A\cup B$ in $X$  to a normal crossing divisor $\widehat D$ in 
$\widehat X$. It provides an isomorphism 
$$
\pi^0: \widehat X - \widehat D \stackrel{\sim}{\lra} X - A \cup B
$$ 

Let $\partial \Delta_B$ be the boundary of  $\Delta_B$,  
so $\partial \Delta_B \subset B(\C)$. Then the open part 
$\Delta^0_B:= \Delta_B - \partial \Delta_B$ sits inside of the set of complex points of 
$X - A \cup B$. Using the isomorphism $\pi^0$ we 
identify $\Delta^0_B$ with an open cell in the set of complex points of 
$\widehat X - \widehat D$. Denote by $\overline \Delta^0_B$ its 
 topological closure there. The Zariski closure 
of the boundary $\partial \overline \Delta^0_B$ is contained in $\widehat D$. 
Let $\widehat B$ be the smallest union of the irreducible components of $\widehat D$ 
containing it. 

The form $\omega_A$ is a regular differential 
form on $X-A$ with logarithmic singularities at $A$. 
The isomorphism $\pi^0$ 
 identifies $\omega_A$  it with a regular differential form $\widehat \omega_A$ 
on $\widehat X - \widehat D$. It automatically has logarithmic singularities 
along  the divisor $\widehat D$. Let $\widehat A$ be the actual divisor of singularities 
of  $\widehat \omega_A$. Suppose that integral (\ref{5.25.02.1}) is convergent. 
Then  one can show that the crucial condition (\ref{5.25.02.10}) is satisfied. 
A detailed elaboration of this construction in the case of the multiple zeta values 
see in [GM]. It works equally well for all multiple polylogarithms. 

However in this paper, aiming to the specialization theorem, we choose a different 
approach. It uses perverse sheaves  on $X$ instead of  blow ups of $X$ to
 construct Hodge/\'etale realizations of the mixed motive related to a period. 
Here is a more detailed account. 

In chapter 2, for so-called admissible pair of divisors  $(A, B)$, we produce 
a framed mixed Hodge structure whose period is given by 
the integral (\ref{1.31.02.5}). More precisely, we construct a perverse 
sheaf of geometric origin ${\cal F}_{A,B}^*$ on $X$ such that 
$$
H^0(X, {\cal F}_{A,B}^*)
$$
is equipped with a natural framing corresponding to the given cohomology classes 
$[\omega_A] \in H^n(X-A)$ and $[\Delta_B] \in H_n(X; B)$. 
We define ${\cal F}_{A,B}^*$ as $H^0$ for the perverse $t$-structure of
 an object ${\cal F}_{A,*, B}^{\bullet} \in D^b_{\rm Sh}(X)$.  

Let $\pi: X \lra {\rm Spec}(F)$ be the canonical projection. 
In chapter 3  we show that, assuming   
 a condition on  $A \cup B$ 
valid in all examples through the paper, the object 
${\rm R}\pi_* {\cal F}_{A,*, B}^{\bullet}$ 
is of geometric origin. This means that 
they are realizations of certain objects of the triangulated category ${\cal D}{\cal M}_F$ 
of mixed motives over $F$ constructed in [V] and [L]. 

If $F$ is a number field and 
the divisor $ A \cup B$ provides a Tate stratification of $X$ we define 
$
H^0(X; {\cal F}^{*}_{A,B})
$
as an object of the abelian category ${\cal M}_T(F)$ 
of mixed Tate motives over $F$ constructed in chapter 5 of [G9].

Often there are   many 
constructions of framed objects with the same period, and it is important to know 
that they lead to the same equivalence class of framed objects. 
Here is how we can face such a situation.  Suppose that we have a deformation 
$(\omega_{A(\varepsilon)}, 
\Delta_{B(\varepsilon)})$ such that  
$A(\varepsilon) \cup B(\varepsilon)$ 
for $\varepsilon \not = 0$ is a normal crossing divisor satisfying (\ref{5.25.02.10}).
Suppose that as $\varepsilon \to 0$ the limit 
of the data $(A(\varepsilon), B(\varepsilon), [\omega_{A(\varepsilon)}], 
[\Delta_{B(\varepsilon)}])$ exists and coincides with the one 
 $(A, B, [\omega_{A}], 
[\Delta_{B}])$. So in particular 
\begin{equation} \label{errra}
\lim_{\varepsilon \to 0}\int_{\Delta_{B(\varepsilon)}}\omega_{A(\varepsilon)} = 
\int_{\Delta_{B}}\omega_{A}
\end{equation} 
Then the specialization 
${\rm Sp}_{\partial/\partial \varepsilon}m(X; [\omega_{A(\varepsilon)}]; [\Delta_{B(\varepsilon)}])$ 
is a framed object whose period is  given by  (\ref{errra}).   
Taking different perturbations we get different mixed objects. 
Nevertheless in chapter 4 we prove that,  
under certain assumptions,  the equivalence classes of the 
obtained framed objects are the same.  (In fact there is no need to assume that 
$A(\varepsilon) \cup B(\varepsilon)$ is a normal crossing divisor for $\varepsilon \not = 0$.)

We construct the motives corresponding to periods by 
resolving singularities and then using  construction (\ref{5.2.02.1}).

Our specialization theorem provides an efficient way to 
compute the coproduct for the framed mixed object 
corresponding to a degenerate admissible configuration $(A,B)$: 
take a generic deformation $(A(\varepsilon), B(\varepsilon))$, 
compute explicitly the coproduct for  $\varepsilon \not = 0$, and 
take the specialization  when $\varepsilon \to 0$. 

Here are some applications. 
In chapter 5 we apply the results above to provide motivic treatment of Aomoto polylogarithms. 
In chapter 6 we use the results of chapters 2-4  to study the torsor of path between 
the tangential base points 
on a curve. In particular we define the framed mixed Tate motives corresponding to 
the classical polylogarithms and compute the coproduct for framed these objects, 
see section  10.12. This implies  
the weak version of Zagier's conjecture on $\zeta_F(n)$.

 {\bf 10. An example: the motivic double logarithm}. To explain 
the nature of the problem we discuss below 
the simplest nontrivial example: construction of the framed mixed Tate motive 
corresponding to double logarithm 
\begin{equation} \label{iop}
{\rm I}_{1,1}(a_1, a_2):= \int_{0 \leq t_1 \leq  t_2 \leq 1}
\frac{dt_1}{t_1-a_1}\wedge\frac{dt_2}{t_2-a_2} 
\end{equation}

i) Let us suppose first that 
\begin{equation} \label{poi}
a_1, a_2 \not \in \{0,1\}
\end{equation}
Consider the following 
two divisors in the plane $X$ with coordinates $(t_1, t_2)$:
$$
A:= \{t_1=a_1\} \cup \{t_2= a_2\}; \qquad B:= \{t_1=0\} \cup \{t_1=t_2\} \cup \{t_2= 1\}
$$
\begin{center}
\hspace{4.0cm}
\epsffile{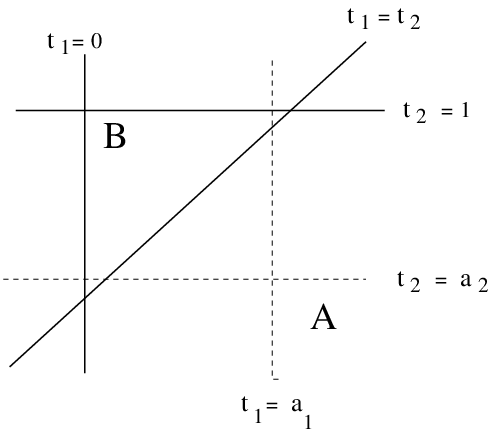}
\end{center}
Let $B_A:= B - B \cap A$. Set 
$$
m(0; a_1, a_2; 1) := H^2(X -A, B_A)
$$
Then, if $F$ is a number field and 
$a_1, a_2 \in F$, this is a mixed Tate motive over $F$ with a  natural framing. 
One component of the framing is provided by 
$$
[\frac{dt_1}{t_1-a_1} \wedge \frac{dt_2}{t_2-a_2}]: \Q(-2) \lra  
{\rm Gr}^W_4H^2(X -A) = {\rm Gr}^W_4H^2(X - A, B_A)
$$
The other is given by canonical isomorphisms
$$
\Q(0)= {\rm Gr}^W_0H_2(X, B_A) = 
{\rm Gr}^W_0H_2(X - A, B_A) 
$$
In the Betti realization it looks as follows. 
Choose an embedding $F \hookrightarrow \C$, so we may assume 
$a_1, a_2 \in \C$. 
Let $\gamma: [0,1] \to {\C} - \{a_1, a_2\}$ be a path between $0$ and $1$. 
Then 
$
\gamma(\Delta_2)  \subset \C^2 - A(\C)
$ 
provides a relative homology class $[\gamma(\Delta_2)]$ generating 
${\rm Gr}^W_0H_2({\C}^2 - A(\C), B_A(\C))$.

ii) Now suppose that 
$$
a_1 \not = 0; \quad a_2 \not = 1
$$
This is precisely the convergence condition for the 
integral (\ref{iop}). 

Observe that  $H^2(X - A, B_A)$  is always a mixed Tate motive. However it has ${\rm Gr}^W_0=0$ 
if $a_1 = 0$ 
or $a_2 = 1$, e.g. for $\zeta(2)$. Indeed, in these cases one of the vertices of the triangle $\gamma(\Delta_2)$ has been removed, so there is no  relative cycle 
needed as a frame component. 

\begin{center}
\hspace{4.0cm}
\epsffile{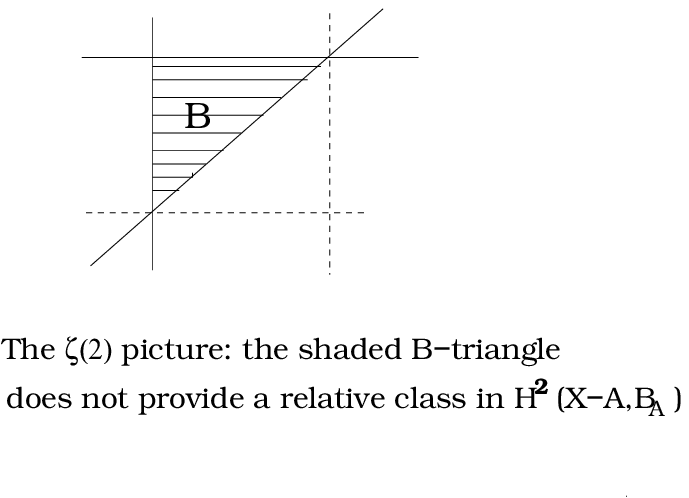}
\end{center}

In any realization one can overcome this problem by setting 
\begin{equation} \label{poi1q}
m(0; a_1, a_2; 1):= H^2(X, {\cal F}_{A,B})
\end{equation}
where ${\cal F}_{A,B}$ is an object of 
the corresponding derived category of sheaves on $X$ 
defined as follows. Take the constant sheaf $\Q(0)$ on 
$X - (A \cup B)$ and extend it first by full direct image 
into $A - (A\cap B)$, and after that by $Rj_!$ into $B$, 
where $j: X - B \hookrightarrow X$. 
Then we get a framed mixed object 
in any realization which is isomorphic to the realization $m(0; a_1, a_2; 1)$ if (\ref{poi}) is valid. 
However since so far there is no theory of motivic sheaves 
we can not interpret the right hand side of 
(\ref{poi1q}) directly  as an object of the triangulated category of motives 
over $F$. 

Therefore we choose a different strategy. Let $\widehat X$  be 
the blow up of the plane $ X$ at 
the vertices $(0,0)$ and $(1,1)$ of the triangle $B$. The preimage 
$\widehat B$ of the triangle $B$ 
is a divisor having shape of the pentagon.  

\begin{center}
\hspace{4.0cm}
\epsffile{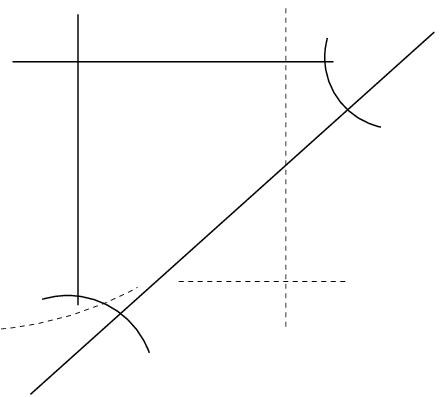}
\end{center}
Denote by $\widehat A$ the strict preimage 
of  $A$ on the blow up. 
Then $\widehat A$ does not contain any vertex of the algebraic 
pentagon $\widehat B$. 
It follows that the mixed Tate motive
$$
m(0; a_1, a_2; 1):= H^2(\widehat X - \widehat A, \widehat B_A ); \qquad \widehat B_A:= \widehat B 
- (\widehat B \cap \widehat A)
$$
does have a natural framing. It is canonically isomorphic 
to (\ref{poi1q}). 

In particular when $a_1=1, a_2=0$ we constructed a motivic avatar for $\zeta(2)$.

iii) When $a_1 = 0$ or $a_2 = 1$ we define  the Hodge or l-adic realization of 
$m(0; a_1, a_2; 1)$ by using the specialization of $m(0; x, y; 1)$ when $x\to 0$ or $y\to 1$. 
The computation of the specialization given in lemma 6.7 in [G7] 
(see also propositions 2.14 and 2.15 loc. cit.) dictates the formulae
$$
m(0; 0, 1; 1) = m(0; 0, 0; 1) = m(0; 1, 1; 1) = 0, 
$$
$$
m(0; 0, a_2; 1) = - m(0; a_2, 0; 1); \quad  m(0;  a_1, 1; 1) = 
- m(0; 1, a_1; 1)
$$
We use them as a {\it definition} of the framed mixed Tate motives $m(0; a_1, a_2; 1)$ 
in the case $a_1 = 0$ or $a_2 = 1$.

{\bf Acknowledgement}. It is my pleasure to thank J. Bernstein for useful discussions 
and help. I appreciate his advise to split the huge original manuscript 
into the pieces provided [G7], this paper, and [G10]. 

The final draft of this paper was prepared during my stay at the University Paris VII, 
MPI(Bonn), 
IHES (Bures sur Yvette) and MSRI (Berkeley). I am very grateful to these 
institutions for hospitality and support. 
This work was supported by the NSF grants DMS-9800998 
and  DMS-0099390.

\begin{center}
\hspace{4.0cm}
\epsffile{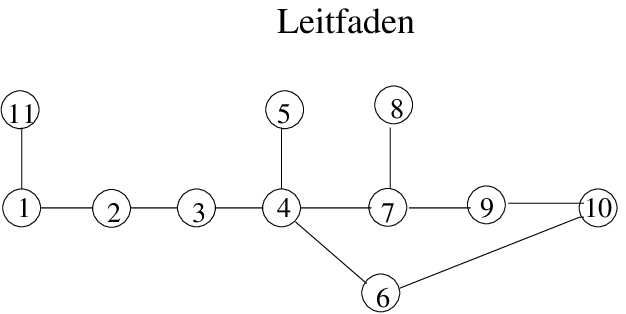}
\end{center}

\section{Framed perverse sheaves related to  periods}


{\bf 1.  The language of perverse sheaves [BBD]}.  
Let $X$ be a variety over a filed $F$. Below $D^b_{\rm Sh}(X)$ stays either for Saito's 
bounded derived category $D_{\cal H}^b(X)$ of mixed Hodge sheaves on $X$ (when $F = \C$), 
or for the bounded derived category $D_{et}^b(X)$ 
of constructible \'etale 
l-adic  sheaves on $X$. 
For a regular irreducible variety $Z$  denote by 
$\delta_{Z}(k)$  the  perverse sheaf $\Q_Z[{\rm dim}Z](k)$ where $\Q_Z = p^*\Q$ (or $p^*\Q_l$ in 
$D_{et}^b(X)$).  
Here $*(k)$ means the Tate twist. For a map $f: X \lra Y$ we denote by $f_!,
f_*, f^!, f^*$ the functors between the derived categories
$D^b_{\rm Sh}(X)$ and $ D^b_{\rm Sh}(Y)$. 
 The notation $R^nf_!{\cal F}$  is used for the $n$-th
cohomology group of  $f_!{\cal F}$. Denote by ${\cal P}_{\rm sh}(X)$ the category of perverse sheaves in $D^b_{\rm Sh}(X)$. 

Let $i: Y \hra X$,  $U:= X-Y$ and  $j:U  \hra X$. 
Then there are  standard exact triangles
\begin{equation} \label{1.11.07.11}
{\cal G}^{\bullet} \lra j_*j^*{\cal G}^{\bullet} \lra i_!i^!{\cal G}^{\bullet}[1]
\end{equation}
\begin{equation} \label{1.11.07.12}
i_*i^*{\cal G}^{\bullet}[-1]  \lra j_!j^!{\cal G}^{\bullet}  \lra {\cal G}^{\bullet}
\end{equation}
Now let $Y$ be a codimension $m$ regular subvariety of a regular variety $X$. 
Then  
\begin{equation} \label{11.15.01.1}
i^!\Q_X = \Q_Y[-2m](-m), \qquad i^*\Q_X = \Q_Y
\end{equation}
So
\begin{equation} \label{11.15.01.1a}
i^!\delta_X = \delta_Y[-m](-m), \qquad i^*\delta_X = \delta_Y[m]
\end{equation}
Therefore 
there are
exact triangles
\begin{equation} \label{1.11.01.1}
\delta_X \lra j_* \delta_U \lra \delta_Y[-m+1](-m)
\end{equation} 
\begin{equation} \label{1.11.01.q1}
\delta_Y[m-1] \lra j_! \delta_U \lra \delta_X 
\end{equation}
In particular  if $Y$ is a regular 
divisor in  $X$ there are short 
exact sequences of perverse sheaves
\begin{equation} \label{1.21.01.8}
0 \lra \delta_X \lra j_*(\delta_{U}) \lra \delta_Y(-1) \lra 0
\end{equation}
\begin{equation} \label{1.21.01.8se}
0 \lra \delta_Y \lra j_!(\delta_{U}) \lra \delta_X \lra 0
\end{equation}

The  
perverse mixed Hodge sheaves  are equipped with  canonical weight
filtration $W_{\bullet}$.  We use the standard convention for the weights: 
${\rm w}(X[n]) = {\rm w}(X) +n$.  
By purity for a regular projective variety $P$ and a pure
perverse sheaf ${\cal F}_a$ of weight $a$ the weight of  $H^i(P, {\cal F}_a)$ is $i+a$. 
The weight of $\delta_Y(-k)$ is ${\rm dim}(Y) + 2k$.

Let ${\cal O}_1$ and ${\cal O}_2$ be two sets of isomorphism classes of 
objects of a triangulated category. 
We define a new set ${\cal O}_1 \ast {\cal O}_2$ of isomorphism classes  of  objects  so that 
$Y \in {\cal O}_1 \ast {\cal O}_2$  if and only if there exists 
an exact triangle $X \lra Y \lra Z$ with $X \in {\cal O}_1$ and 
$Z \in {\cal O}_2$. 
We say that an object is glued from  the set of objects $\{A_i\}$ 
if it is isomorphic to an object from $\{A_i\} \ast ... \ast \{A_i\}$. 

\begin{lemma} \label{11.14.01.1} Suppose that a perverse \'etale $l$-adic 
sheaf ${\cal F}$ on $X$ is glued from the objects $\delta_{C_i}(m)$ where 
$X$ and $C_i$'s are regular and defined over a 
 field $F$. Then  ${\cal F}$ 
has a natural weight filtration.
\end{lemma}

{\bf Proof}. If $X$ is defined over a finite field $F_q$ 
there is a weight filtration on perverse sheaves, see [BBD]. 
We use the reduction to a finite field  to deduce the lemma from this. 
One can find a model of $X$ and $C_i$'s over a ring ${\cal O}$ which is 
finitely generated over $\Z$. 
Choose a generic $F_q$-point $i_{\varepsilon}: {\rm Spec}F_q \hra 
 {\rm Spec}{\cal O}$. 
Denote by $\overline X, \overline C_i$ the reduction 
of the schemes $X, C_i$ modulo 
 the maximal ideal of ${\cal O}$ corresponding to $i_{\varepsilon}$. 
For 
generic $F_q$-points $\varepsilon$ the restriction functor $i_{\varepsilon}^*$ 
commutes with 
the standard functors. In particular this implies that 
$$
Ext^1_X(\delta_Y(n), \delta_Z(m)) = 
Ext^1_{\overline X}(\delta_{\overline Y}(n), \delta_{\overline Z}(m))
$$
The group on the right vanishes when ${\rm w}(\delta_Y(n)) \leq {\rm w}(\delta_Z(m))$.
It remains to use lemma 5.3.6 in [BBD]. 
The lemma follows.

{\bf 2. Perverse sheaves corresponding to an admissible pair of divisors}. 
{\it The stratification of $X$ defined by a divisor $D$}.  
 Let $D$ be a divisor whose 
 irreducible components $D_i$, $i \in \Lambda$, are divisors. 
For each point $x \in X$ let $I_x$  be the set parametrizing 
all irreducible components of $D$ containing $x$. 
The subvariety formed by all $x \in X$ with given  $I_x = I$ 
is called the $I$-{\it stratum} defined by $D$. 

{\it Regular stratifications}. 
We say that the stratification defined by the divisor 
$D$ is {\it regular} 
if  all $I$-strata  are irreducible, and  the closure of any stratum 
 is a regular subvariety. The closure of the 
$I$-stratum is $D_I := \cap_{i \in I}D_i$.

{\it Admissible pairs of divisors and their faces}. 
Let $X$ be a regular projective variety over a field $F$ of dimension $n$,  
$A$ and $B$ are two divisors. 
Let
$$
A = A_1 \cup ... \cup A_s; \quad B = B_1 \cup ... \cup B_t
$$
be the decomposition of $A$ and $B$ into irreducible components. For subsets
$$
I = \{i_1, ..., i_k\} \subset \{1, ..., s\}; \quad J= \{j_1, ..., j_l\} \subset \{1, ..., t\}
$$
we set $A_I = A_{i_1} \cap ... \cap A_{i_k}$ and $B_J = B_{j_1} \cap ... 
\cap B_{j_l}$ and call them the 
$A$- and $B$-{\it faces}. 

\begin{definition} \label{pp}
A pair $(A, B)$ is {\it admissible} if the divisor $A \cup B$ provides a regular 
stratification of $X$, and no $A$- and $B$-faces  coincide. 
\end{definition}

We allow the $A$-faces to be  subvarieties of  $B$-faces of positive codimension, and vice versa. 
We assume that the pair $(A, B)$ is admissible. 

{\it Pure $A$ and $B$- strata}. Set
$$
A^0_I := A_I - \{ \mbox{union of all codimension $>0$ strata of $A_I$}
\}
$$
and similarly $B_J^0$. We say that $A^0_I$ and $B_J^0$ are {\it pure} $A$- and $B$- 
{\it strata}. 

Observe that  in general $A$- or $B$- faces are not strata. They 
are closures of the corresponding pure $A$- and $B$- strata.

Let $U:= X - (A\cup B)$. We define an object ${\cal F}^{\bullet}_{A,B} \in 
D_{\rm Sh}^b(X)$  as 
a certain extension of 
$\delta_U$ to $X$. This is done by an inductive procedure
consisting of $n= {\rm dim} X$ steps. Namely, consider open subsets
 $$
U(k):= X - \mbox{all codimension $>k$ strata}, 
$$
 so 
$$
U = U(0) \hra U(1) \hra  ... \hra U(n-1) \hra U(n) = X
$$
 Then 
$U(k) - U(k-1)$ is the disjoint union of the codimension $k$ strata of three  types: 
pure $A$-strata, pure $B$-strata, 
and mixed strata. 

We define the  {\it extension to codimension $k$ strata} functor 
$$
E^*_k: D^b_{\rm Sh}(U(k-1)) \lra D^b_{\rm Sh}(U(k))
$$
 as the extension by $j_*$ to all pure codimension $k$ $A$-strata and mixed  strata, 
followed by  the extension by $j_!$ to all codimension $k$ pure $B$-strata. 
Precisely, consider an open subset 
$$
U^A_k := U(k-1) \cup \mbox{all codimension $k$ pure $A$-strata and mixed strata} =
$$
$$
U(k) - \mbox{all codimension $k$ pure $B$-strata}
$$
Then we have 
$$
U(k-1) \quad \stackrel{j_{k}^A}{\hra} \quad U^A_k \quad \stackrel{j_{k}^B}{\hra} \quad U(k);
\qquad E^*_k:= j_{k!}^{B}  j_{k*}^A 
$$
Set
$$
{\cal F}^{\bullet}_{A,*,B}:= \quad E^*_n  ...   E^*_2
   E^*_1 \delta_U \quad = \quad 
j_{n!}^B j_{n*}^{A} ...  j_{1!}^B j_{1*}^{A} \delta_U
$$
and consider the perverse sheaf 
$$
{\cal F}^*_{A,B}:= H^0_{\tau}{\cal F}^{\bullet}_{A, *, B}
$$
where $H^{*}_{\tau}$ are the cohomology groups with respect to the $t$-structure 
on $D_{\rm Sh}^b(X)$ providing the abelian category of  perverse sheaves on $X$.

One can define a similar functor $E_k^!$ 
by using the $j_{!}$-extension to the mixed 
strata. In particular set
$$
{\cal F}^!_{A,B}:= H^0_{\tau}\Bigl(E^!_n  ...   E^!_2
   E^!_1 \delta_U \Bigr)
$$
More generally, on can define  similar objects  by using either 
$E_k^!$ or $E_k^*$ for each given $k$. 

{\bf Remark}.  Since the strata of different types are disjoint, 
the extensions by $j_!$ and by $j_*$  commute. Therefore 
$$
\ast{\cal F}^!_{A,B} ={\cal F}^*_{B, A} 
$$
where $\ast$ is the duality functor.

{\it The frame morphisms}. 
Let $j_A:X -A \hra X$, $j_B: X - B \hra X$. 
Set
$$
{\cal F}_A:= j_{A*}\delta_{X -A}; \qquad {\cal F}_B:= j_{B!}\delta_{X - B}
$$

\begin{lemma} \label{11.7.01.1}
There are canonical morphisms
\begin{equation} \label{11.5.01.2}
{\cal F}_B \stackrel{\beta}{\lra}{\cal F}^!_{A,B} 
\stackrel{c}{\lra}{\cal F}^*_{A,B}\stackrel{\alpha}{\lra} {\cal F}_{A} 
\end{equation}
\end{lemma}

{\bf Proof}. Let us define the map $\alpha$. We define by induction the maps
\begin{equation} \label{11.7.01.14}
\alpha_k: {\rm Res}_{U(k)}{\cal F}^{\bullet}_{A,*,B} \lra {\rm Res}_{U(k)}{\cal F}_{A}
\end{equation}
and set $\alpha:= H^0_{\tau}\alpha_n$. The map $\alpha_0$ is the obvious isomorphism. 
The  map $\alpha_k$ is constructed by combining the following  
 two  general observations, which  
will also play an 
important role below (especially in chapter 3).

i) Suppose that ${\cal G}^{\bullet} \in D^b_{\rm sh}(X)$ and there is a map
\begin{equation} \label{12.24.01.1}
{\rm Res}_{U(k-1)}{\cal G}^{\bullet} \lra {\rm Res}_{U(k-1)}{\cal F}_A
\end{equation}
Then we have a map
\begin{equation} \label{12.24.01.11}
j_{k*}^A{\rm Res}_{U(k-1)}{\cal G}^{\bullet} \lra {\rm Res}_{U^A(k-1)}{\cal F}_A
\end{equation}
given by applying $j_{k*}^A$ to (\ref{12.24.01.1}):
$$
j_{k*}^A{\rm Res}_{U(k-1)}{\cal G}^{\bullet} \stackrel{}{\lra} 
j_{k*}^A{\rm Res}_{U(k-1)}{\cal F}_A = 
{\rm Res}_{U^A(k-1)}{\cal F}_A
$$

ii) Suppose that ${\cal G}^{\bullet} \in D^b_{\rm sh}(X)$ and there is a map
\begin{equation} \label{12.24.01.112}
{\rm Res}_{U^A(k-1)}{\cal G}^{\bullet} \lra {\rm Res}_{U^A(k-1)}{\cal F}_A
\end{equation}
Then we have a map
\begin{equation} \label{12.24.01.12}
j_{k!}^B{\rm Res}_{U^A(k-1)}{\cal G}^{\bullet} \lra {\rm Res}_{U(k)}{\cal F}_A
\end{equation}
given as a composition 
$$
j_{k!}^B{\rm Res}_{U^A(k-1)}{\cal G}^{\bullet} \stackrel{(\ref{12.24.01.112})}{\lra} 
j_{k!}^B{\rm Res}_{U^A(k-1)}{\cal F}_A = j_{k!}^Bj_{k}^{B!}{\rm Res}_{U(k)}{\cal F}_A  \lra {\rm Res}_{U(k)}{\cal F}_A
$$
The last map  is the adjunction morphism.

Combining the constructions i) and ii) 
for ${\cal G}^{\bullet} = {\cal F}^{\bullet}_{A,*,B}$ we get 
the map $\alpha_k$.

The map $\beta$ is defined 
using the duality $\ast$ interchanging the role of $A$ and $B$ divisors.  
The map $c$ is defined  using the canonical morphism $j_!\lra j_*$ 
for each extension to a mixed 
stratum. The lemma is proved.

{\bf 3. The structure of  $H^0(X, {\cal F}_{A,B}^*)$}. 
We are going to study the following perverse sheaves:
\begin{equation} \label{11.7.01.3}
{\cal F}^*_{A,B}, {\cal F}^!_{A,B}, {\cal F}_{A}, {\cal F}_{B}
\end{equation}

\begin{proposition} \label{10.24.01.2} Suppose that $(A,B)$ is an admissible 
 pair of divisors on $X$. Then 
each of the objects (\ref{11.7.01.3}) is glued from $\delta_C(-m)$,  
where   $C$ runs 
through the 
strata of the stratification defined 
by $A \cup B$, and $0 \leq m \leq  {\rm codim}_X C$. 
\end{proposition}

{\bf Proof}. 
If $F: {\cal C}_1 \lra {\cal C}_2$ is a triangulated functor between 
triangulated categories, an object $X$ is glued from the objects $\{A_i\}$, and 
every object $F(A_i)$ is glued from the set of 
objects $\{B_j\}$, then $F(X)$ is glued from $\{B_j\}$. So using 
(\ref{1.11.01.1})-(\ref{1.11.01.q1}) we see that 
each of the objects (\ref{11.7.01.3}) is glued from $\delta_C(-m)[p]$ 
where $0\leq m \leq {\rm codim}_XC$. 
Further, one easily proves by induction that 
if an object of an abelian category ${\cal A}$ 
is glued from the objects $\{A_i[p_i]\}$ of the derived category $D^b({\cal A})$ 
 where $A_i \in {\cal A}$ then it is glued 
from those of these objects which belong to ${\cal A}$, i.e. have $p_i=0$. 
Since objects  (\ref{11.7.01.3}) are from the abelian category of 
the perverse sheaves 
the proposition follows.

Let us set
$$
{\rm dw}(\delta_C(-m)) = {\rm dim}_XC + m = {\rm w}(\delta_C(m)) - m
$$

\begin{proposition} \label{11.7.01.131} Suppose that $Y$ and $Z$ are two  
closed irreducible regular 
subvarieties of a regular variety $X$. Then in  ${\cal P}_{\rm Sh}(X)$ one has 
\begin{equation} \label{11.26.01.5}
Hom_X(\delta_Y(n), \delta_Z(m)) = 0 \quad \mbox{unless $Y=Z$, $n=m$}
\end{equation}
$$
Ext_X^1(\delta_Y(n), \delta_Z(m)) = 0 \quad \mbox{if ${\rm dw}(\delta_Y(n)) < 
{\rm dw}(\delta_Z(m))$}
$$
$$
Ext_X^1(\delta_Y(-n), \delta_Z(-m)) = 0 \quad \mbox{if $m>n$}
$$
\end{proposition}

{\bf Proof}. 
The objects $\delta_Y(n)$ are irreducible, thus we have the first statement. 
The case $Y=Z$ is clear from the weight considerations. So we may assume $Y \not = Z$. 
We claim that 
$$
Ext_X^1(\delta_Y(n), \delta_Z(m)) = \left\{ \begin{array}{ll}
0 &  \quad \mbox{if $|d_{Y,Z}|>1$ } \\
0& \quad \mbox{if  none of the varieties $Y,Z$ contains the other}\end{array} \right.
$$
The second claim is obvious, so we may assume without loss of generality 
that $Z \subset Y$. Let $d_W = {\rm dim}W$ and $d_{Y,Z}:= d_Y - d_Z$. 
Then one has 
$$
Ext_X^i(\delta_Y(n), \delta_Z(m)) = Ext_X^i(\delta_Z(-m), \delta_Y(-n+d_{Y,Z})) = 
H^{i-d_{Y,Z}}(Z, \Q(m-n)) 
$$
Indeed,  to compute the first $Ext$ observe that 
$$
Hom_{X}(\Q_Y[d_Y](n), \Q_Z[d_Z+i](m)) = 
$$
$$
Hom_{Y}(\Q_Y, \Q_Z[d_Z-d_Y+i](m-n)) = 
H^{i-d_{Y,Z}}(Z, \Q(m-n))
$$
The second equality follows by the  duality since 
\begin{equation} \label{11.12.01.s2}
\ast \delta_C(-m) = \delta_C(m + {\rm dim}C)
\end{equation}
So to prove the proposition it remains to consider the case when 
$|d_{Y,Z}| = 1$. The weight considerations shows that 
\begin{equation} \label{11.12.01.1}
Ext_X^1(\delta_Y(n), \delta_Z(m)) = 0 \quad \mbox{if $d_{Y,Z} = 1$ and  $m-n < 0$}
\end{equation}
and 
\begin{equation} \label{11.12.01.2}
Ext_X^1(\delta_Y(n), \delta_Z(m)) = 0 \quad \mbox{if $d_{Y,Z} = - 1$ and  $m-n \leq  0$}
\end{equation}
Indeed, under these assumptions 
 $$
{\rm w}(\delta_Y(n)) - {\rm w}(\delta_Z(m))  = d_{Y,Z} + 2(m-n) < 0
$$ 
   So
the $Ext^1$-groups vanish. Observe also that 
(\ref{11.12.01.2}) can be deduced from (\ref{11.12.01.1}) by duality thanks to 
(\ref{11.12.01.s2}). 
This gives the third statement 
of the proposition. 
To check the second  notice that by the assumption in the second statement
$$
{\rm dw}(\delta_Y(n)) - {\rm dw}(\delta_Z(m)) = d_{Y, Z}+ (m-n) <0
$$
This means that 
$$
d_{Y, Z}= -1, m-n \leq 0 \quad \mbox{or} \quad d_{Y, Z}= 1, m-n < -1
$$ 
In the first case $Ext^1$ vanishes by  (\ref{11.12.01.2}), and 
in the second by  (\ref{11.12.01.1}). 
The second statement and hence the proposition are proved.

\begin{corollary} \label{11.9.01.211} Let ${\cal M}$ be an object of ${\cal P}_{\rm Sh}(X)$ 
on a regular variety $X$. 
Suppose that  ${\cal M}$ is glued from the objects 
$\delta_C(n)$, where $C$'s are regular closed irreducible subvarieties of $X$. Then 
${\cal M}$  
carries two canonical increasing  filtrations indexed by integers: 

a)  The  dimension-weight filtration $DW_{\bullet}$ such that 
${\rm Gr}^{DW}_k{\cal M}$ is a direct sum  of the objects $\delta_C(-m)$ with 
${\rm dim}C + m = k$. 

b) The filtration $T_{\bullet}$ by the Tate twist   
such that 
${\rm Gr}^{T}_k{\cal M}$ is a direct sum  of the objects $\delta_C(-k)$. 

c) The duality $\ast$ interchanges these two filtrations, so that
$$
\ast{\rm Gr}^{T}_k{\cal M} = {\rm Gr}^{DW}_{-k}\ast{\cal M}
$$
\end{corollary}
 
{\bf Proof}. The parts a) and b) follow from proposition 
\ref{11.7.01.131} and lemma 5.3.6 in [BBD]. 
Namely, for the   dimension-weight filtration we use the  second 
statement of the proposition, and for the filtration by the Tate twist the third. 
The part c) follows from a) and b) using formula (\ref{11.12.01.s2}). 
The corollary is proved. 

Recall that $n={\rm dim}X$. 

\begin{proposition} \label{11.7.01.12} The morphisms  (\ref{11.5.01.2}) induce isomorphisms 
\begin{equation} \label{11.19.01.1}
{\rm Gr}^{DW}_n({\cal F}^{!}_{A,B}) \stackrel{\sim}{\lra} 
{\rm Gr}^{DW}_n({\cal F}^{*}_{A,B}) \stackrel{\sim}{\lra} {\rm Gr}^{DW}_n({\cal F}_{A})
\end{equation}
\begin{equation} \label{11.19.01.2}
{\rm Gr}^{T}_0({\cal F}_{B}) \stackrel{\sim}{\lra} 
{\rm Gr}^{T}_0({\cal F}^{!}_{A,B}) \stackrel{\sim}{\lra} {\rm Gr}^{T}_0({\cal F}^*_{A,B})
\end{equation}
\end{proposition}

{\bf Proof}. The second line of isomorphisms follows from the first one by the duality using 
corollary \ref{11.9.01.211}. 

Let us prove the last isomorphism in (\ref{11.19.01.1}). 
Let us show by induction that 
the maps $\alpha_k$, see (\ref{11.7.01.14}), induce  isomorphisms on 
${\rm Gr}^{DW}_nH^0_{\tau}(-)$. If $k=0$ this is clear. 

For a set  $\{\delta_{C_i}(-m_i)[p_i]\}$ 
let  
\begin{equation} \label{11.19.01.4}
{\rm dw}_+(\{\delta_{C_i}(-m_i)[p_i]\}) = 
{\rm max}_i({\rm dim}C_i + m_i)
\end{equation}
 Suppose that ${\cal M}$ is glued from the  objects $\{\delta_{C_i}(-m_i)[p_i]\}$. 
There are many sets of such objects 
generating ${\cal M}$ 
(for instance we can  add an object to a list). Let 
${\rm dw}_+({\cal M})$ be the smallest 
value invariant (\ref{11.19.01.4}) can have for such a generating  set. 

\begin{lemma} \label{11.15.01.10} Let $Y$ be a regular subvariety of a 
regular variety $X$, $U:= X-Y$, and  $j:U \hra X$. 
Suppose that ${\cal M}$ is glued from $\delta_C(-m)[p]$ where all subvarieties $C$
 and $C \cap Y$ are regular. 
Then 
$$
{\rm dw}_+(j_*{\cal M}) \leq  {\rm dw}_+({\cal M}), \quad  
{\rm dw}_+(j_!{\cal M}) \leq  {\rm dw}_+({\cal M})
$$
 \end{lemma}

{\bf Proof}. Follows immediately from (\ref{1.11.07.11}), (\ref{1.11.07.12}) 
and (\ref{11.15.01.1a}). 
The lemma is proved. 

By lemma \ref{11.15.01.10} and using (\ref{1.11.01.1})-(\ref{1.11.01.q1}) we see by induction on $k$ 
that ${\cal G}^{\bullet}= {\cal G}^{\bullet}_*$ satisfies the following condition: 

for $U = U(k)$ and  $U = U^A(k)$, $k \geq 0$,  one has 
\begin{equation} \label{12.26.01.1}
{\rm dw}_+\Bigl(H^i_{\tau}({\rm Res}_{U}
{\cal G}^{\bullet})\Bigr) \leq \left\{ \begin{array}{ll}
n &   \\
n-1 & \quad  \mbox{for $i < 0$}\end{array} \right. 
\end{equation}

\begin{lemma} \label{12.25.01.1} Let ${\cal G}^{\bullet} \in D^b_{\rm sh}(X)$. 
Assume that ${\cal G}^{\bullet}$ is glued from $\delta_C(-m)[p]$ where 
$C$ are strata of a regular stratification of $X$ and 
$0 \leq m \leq {\rm codim}_X C$. Let us also assume (\ref{12.26.01.1}). Then 

i) If map (\ref{12.24.01.1}) induces an isomorphism on ${\rm Gr}^{WD}_nH^0_{\tau}(-)$ 
then map (\ref{12.24.01.11}) has the same property.

ii) If map (\ref{12.24.01.112}) induces an isomorphism on ${\rm Gr}^{WD}_nH^0_{\tau}(-)$ 
then map (\ref{12.24.01.12}) has the same property.
\end{lemma}

{\bf Proof}. i) 
We claim that 
\begin{equation} \label{11.18.01.1n}
{\rm Gr}^{DW}_n H^0_{\tau}\Bigl(j^A_{k*}{\rm Res}_{U(k-1)}
{\cal G}^{\bullet}\Bigr) \stackrel{\sim}{=} 
{\rm Gr}^{DW}_n H^0_{\tau}\Bigl(j^A_{k*}H^0_{\tau}({\rm Res}_{U(k-1)}
{\cal G}^{\bullet})\Bigr) 
\end{equation}
Indeed, it follows from (\ref{12.26.01.1}) and  lemma \ref{11.15.01.10}  
that  the left object in 
(\ref{11.18.01.1n}) is isomorphic to
$$
{\rm Gr}^{DW}_n H^0_{\tau}\Bigl(j^A_{k*}{\tau}_{\geq 0}({\rm Res}_{U(k-1)}
{\cal G}^{\bullet})\Bigr) 
$$
On the other hand by (\ref{12.26.01.1}) and (\ref{1.11.01.1})-(\ref{1.11.01.q1}) 
$$
{\rm Gr}^{DW}_n H^0_{\tau}\Bigl(j^A_{k*}{\tau}_{\geq 1}({\rm Res}_{U(k-1)}
{\cal G}^{\bullet})\Bigr) =0
$$
Isomorphism (\ref{11.18.01.1n}) follows from these two facts. 

Similarly, by (\ref{12.26.01.1}) and lemma \ref{11.15.01.10} 
the right object in (\ref{11.18.01.1n}) is isomorphic to 
$$
{\rm Gr}^{DW}_n H^0_{\tau}\Bigl(j^A_{k*} {\rm Gr}^{DW}_n H^0_{\tau}  ({\rm Res}_{U(k-1)}
{\cal G}^{\bullet})\Bigr) 
$$
By the induction assumption it is isomorphic to 
$$
{\rm Gr}^{DW}_n H^0_{\tau}\Bigl(j^A_{k*} {\rm Gr}^{DW}_n {\rm Res}_{U(k-1)}
{\cal F}_{A}\Bigr) 
$$
(Observe that ${\cal F}_{A}$ is a perverse sheaf). By 
lemma \ref{11.15.01.10} this object, and hence the right 
object in (\ref{11.18.01.1n}), are isomorphic to 
$$
{\rm Gr}^{DW}_n H^0_{\tau}\Bigl(j^A_{k*}{\rm Res}_{U(k-1)}
{\cal F}_{A}\Bigr) 
$$
It remains to notice that 
$
j^A_{k*}{\rm Res}_{U(k-1)}
{\cal F}_{A} = {\rm Res}_{U^A(k-1)}
{\cal F}_{A}
$. 
The statement i) is  proved. 

To prove statement ii) we follow the same pattern,  
changing $U^A(k-1)$ to $U(k-1)$ and $j^A_{k*}$ to $j^B_{k!}$. 
In the end we need to show in addition that 
$$
{\rm Gr}^{DW}_n H^0_{\tau}\Bigl(j^B_{k!}j^{B!}_{k}{\rm Res}_{U(k)}
{\cal F}_{A}\Bigr) = {\rm Gr}^{DW}_n H^0_{\tau}\Bigl({\rm Res}_{U(k)}
{\cal F}_{A}\Bigr)
$$
There is a standard exact triangle where $i_k^B$ is the closed embedding complementary 
to $j_k^B$:
$$
i_{k*}^Bi_k^{B*}{\rm Res}_{U(k)}{\cal F}_{A}[-1]  \lra j_{k!}^Bj_k^{B!}{\rm Res}_{U(k)}{\cal F}_{A}\lra {\rm Res}_{U(k)}{\cal F}_{A}
$$
So it is sufficient to show that the left object is glued from 
the objects 
$\delta_C(-m)[p]$ with ${\rm dim}C + m < n$. Since ${\cal F}_{A}$ and hence 
${\rm Res}_{U(k)}{\cal F}_{A}$ are glued from such objects with 
${\rm dim}C + m \leq n$, and the stratification defined by $A \cup B$ is regular, 
the statement follows from the second part of  (\ref{11.15.01.1a}). 
The lemma is proved. 

The proof of the first isomorphism in (\ref{11.19.01.1}) is even simpler: 
it uses only the argument similar to the one in the end of the proof of part ii). 
The proposition is proved. 

 Let  ${\rm dp}(A)$ (resp. ${\rm dp}(B)$) 
be the codimension of the 
smallest pure $A$-strata (resp. pure $B$-strata) in $X$. 

\begin{theorem} \label{11.5.01.1}  
In the category ${\cal P}_{\rm Sh}(X)$ 

a) the morphisms (\ref{11.5.01.2})  induce canonical isomorphisms
\begin{equation} \label{1.21.01.23}
{\rm Gr}^W_{2n}H^*(X, {\cal F}^!_{A,B}) \stackrel{\sim}{\lra} 
{\rm Gr}^W_{2n}H^*(X, {\cal F}^*_{A,B})\stackrel{\sim}{\lra}
{\rm Gr}^W_{2n}H^*(X, {\cal F}_{A}) 
\end{equation}
\begin{equation} \label{1.21.01.22}
{\rm Gr}^W_{0}H^*(X, {\cal F}_B) \stackrel{\sim}{\lra} 
{\rm Gr}^W_{0}H^*(X, {\cal F}^!_{A,B}) \stackrel{\sim}{\lra} 
{\rm Gr}^W_{0}H^*(X, {\cal F}^*_{A,B})
\end{equation}

b) there are canonical isomorphisms
$$
{\rm Gr}^W_{n+{\rm dp}(A)}H^0(X, {\cal F}^!_{A,B}) \stackrel{\sim}{\lra} 
$$
\begin{equation} \label{11.19.01.6}
{\rm Gr}^W_{n+{\rm dp}(A)}H^0(X, {\cal F}^*_{A,B})\stackrel{\sim}{\lra} 
{\rm Gr}^W_{n+{\rm dp}(A)}H^0(X, {\cal F}_{A}) 
\end{equation}
$$
{\rm Gr}^W_{n-{\rm dp}(B)}H^0(X, {\cal F}_B) \stackrel{\sim}{\lra} 
$$
\begin{equation} \label{11.19.01.7}
{\rm Gr}^W_{n-{\rm dp}(B)}H^0(X, {\cal F}^!_{A,B}) \stackrel{\sim}{\lra} 
{\rm Gr}^W_{n-{\rm dp}(B)}H^0(X, {\cal F}^*_{A,B})
\end{equation}
\end{theorem}

{\bf Proof}. a)
Let us calculate the groups 
${\rm Gr}^W_wH^p(X, {\cal F})$ using the spectral
sequence for the weight filtration on the perverse sheaf ${\cal F}$. 
Later on ${\cal F}$ will be one of the sheaves (\ref{11.7.01.3}). 
The $E_1$-term of the spectral sequence consists of the groups
$H^i(X; {\rm Gr}^W_a{\cal F})$ with the differential 
\begin{equation} \label{1.21.01.4}
H^i(X; {\rm Gr}^W_a {\cal F})  \lra H^{i+1}(X; {\rm Gr}^W_{a-1}{\cal F})
\end{equation}
Observe that by purity the weight of $H^a(X, {\rm Gr}^W_b{\cal F})$ is $a+b$. 
So the higher differentials are zero by the weight considerations. 
The $E_1$-term of the spectral sequence consists of complexes sitting on the lines 
where the weight 
$a+b$ is a constant.  
In particular the
weight zero and $2n$ parts of this spectral sequence are the complexes formed
by the groups 
\begin{equation} \label{1.21.01.7}
H^{-a}(X; {\rm Gr}^W_a{\cal F})\qquad \mbox{and} \quad H^{a}(X; {\rm Gr}^W_{2n-a}{\cal F})
\end{equation}
respectively, 
with the differential provided  by (\ref{1.21.01.4}).

If $C$ is a regular subvariety of $X$ then 
$$
H^a(X, \delta_C(-m)) =0 \quad \mbox{for $|a|> {\rm dim}C$}
$$
It follows that 
$$
H^{-a}(X; {\rm Gr}^W_a{\cal F}) = H^{-a}(X; {\rm Gr}^W_a(T_0{\cal F}))
$$
Indeed, by proposition \ref{10.24.01.2} the object ${\rm Gr}^W_a{\cal F}$ 
is a direct sum of $\delta_C(-m)$ where ${\rm dim} C+2m=a$. So
if $m>0$ then  $H^{-a}(X; \delta_C(-m))=0$. This and proposition \ref{11.7.01.12} gives 
 (\ref{1.21.01.22}). 
The isomorphism  (\ref{1.21.01.23}) follows by duality.  
The part a) is proved.  

b) {\it The shape of the $E_1$-term of the spectral sequence.} 
If ${\cal F}$ is one of the sheaves (\ref{11.7.01.3}) then 
$$
H^a(X, {\rm Gr}^W_b{\cal F}) =0 \quad \mbox{if $b\pm a >2n$ or $b\pm a <0$}
$$
So the $E_1$-term of the spectral sequence is located inside of the diamond 
$0 \leq b \pm a \leq 2n$. Moreover, it 
is bounded by the lines $n-{\rm dp}(B) \leq b\leq n+{\rm dp}(A)$. 
A particular  interesting case is  ${\rm dp}(A) = {\rm dp}(B) =n$. 

\begin{center}
\hspace{4.0cm}
\epsffile{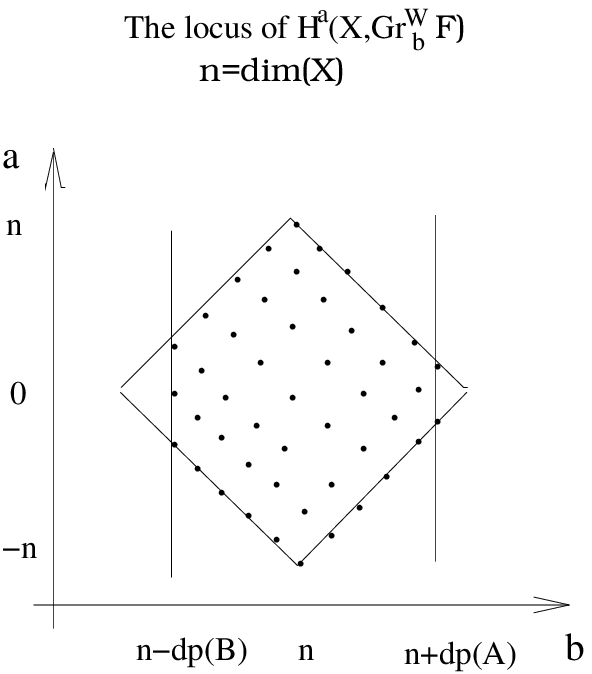}
\end{center}

It follows that  for any perverse sheaf ${\cal F}$ glued from 
$\{\delta_C(-m)\}$ with $0 \leq m \leq {\rm codim}C \leq {\rm dp}(A)$ one has 
$$
{\rm Gr}^W_{n+{\rm dp}(A)}H^0(X, {\cal F}) = 
{\rm Gr}^W_{n+{\rm dp}(A)}H^0(X, {\rm Gr}^{DW}_n{\cal F})
$$ 
The isomorphisms (\ref{11.19.01.6}) follow from this and proposition \ref{11.7.01.12}. 
The isomorphisms (\ref{11.19.01.7}) are obtained by duality.
The part b) and hence  
the theorem are proved. 

{\bf 4. The framing in the Tate case}. 
A regular projective variety $X$ is a 
{\it Tate variety} if the motive of $X$ is a direct sum 
of the Tate motives $\Q(m)$. For example a flag variety $G/P$ where $P$ is a 
parabolic subgroup of a simple Lie group is Tate. The product of Tate varieties is  a Tate variety. 
If all strata of a divisor $D$ on $X$ (including $D_{\emptyset} = X$) 
are Tate varieties we say 
that $D$ provides a Tate stratification of $X$. 
For example a union of hyperplanes  provides a Tate stratification of $P^n$. 

If $A \cup B$ provides a Tate stratification of $X$ and $(A,B)$ is an admissible pair then 
obviously applying the functor 
$H^0(X; -)$ to any of the perverse sheaves (\ref{11.7.01.3}) we get a mixed Tate 
object over the point. 

Recall (see, say,  [G7] or [G8]) 
that a framed object in a mixed Tate category ${\cal C}$ is defined as a triple 
$(H, v_{2n}, f_0)$ where $v_{2n}: \Q(-n) \lra {\rm Gr}^W_{2n}H$ and 
$f: {\rm Gr}^W_{0}H \lra \Q(0)$ are non zero morphisms. We define 
 an equivalence of framed mixed Tate objects as the finest equivalence relation such that 
a morphism compatible with frames is an equivalence. The equivalence classes 
of objects framed by $\Q(0)$ and $\Q(-n)$ form an abelian group 
${\cal A}_n({\cal C})$, and ${\cal A}_n({\cal C}):= \oplus_{n\geq 0}{\cal A}_n({\cal C})$
has a commutative Hopf algebra structure. 

Suppose ${\rm dp}(A) = {\rm dp}(B) = n$. Let us assume that
we are given non zero classes
$$
[\omega_A] \in {\rm Hom}(\Q(-n), {\rm Gr}^W_{2n}H^0(X, {\cal F}_{A})); \qquad [\Delta_B] \in 
{\rm Hom}({\rm Gr}^W_{0}H^0(X, {\cal F}_B), \Q(0))
$$
Then using canonical isomorphisms (\ref{11.19.01.6})-(\ref{11.19.01.7}) 
we transfer them to 
distinguished classes 
$$
[\omega_A] \in {\rm Hom}(\Q(-n), {\rm Gr}^W_{2n}H^0(X, {\cal F}^*_{A,B})) =  
{\rm Hom}(\Q(-n), {\rm Gr}^W_{2n}H^0(X, {\cal F}^!_{A,B}))
$$
and 
$$
[\Delta_B] \in {\rm Hom}({\rm Gr}^W_{2n}H^0(X, {\cal F}^*_{A,B}), \Q(0)) =  
{\rm Hom}({\rm Gr}^W_{2n}H^0(X, {\cal F}^!_{A,B}), \Q(0))
$$
Therefore we get equivalent framed mixed Tate objects 
\begin{equation} \label{11.14.01.10}
(H^0(X, {\cal F}^*_{A,B}), [\omega_A], [\Delta_B]) \quad \sim \quad 
(H^0(X, {\cal F}^!_{A,B}, [\omega_A], [\Delta_B])
\end{equation}

{\bf 5. The framing in the general case}. The general 
 definition of framed objects in a mixed category see in chapter 2 of 
[G8]. 
In this paper we employ an adapted version of this 
definition 
sufficient for our  goals. The reader is invited to formulate 
and prove the results of the next chapters using the definition given in [G8].  

Let ${\cal P}$ be a semisimple category. We assume that each simple object of ${\cal P}$ 
is labelled by an integer called the  weight. Let ${\cal C}_{\cal P}$ 
be a mixed category such that every its object 
carries canonical  weight filtration $W_{\bullet}$ so that 
${\rm Gr}^W_k$ is a direct sum of weight $k$ objects of 
${\cal P}$. 
Let us assume that we are given a fiber functor $F: {\cal P} \lra {\rm Vect}$ to the 
category of finite dimensional vector spaces. Then a framing on a mixed  object 
$H$ is a triple $(H, v, f)$ defined by a pair of non zero vectors 
$$
v \in F({\rm Gr}^W_pH), \quad f \in F({\rm Gr}^W_qH)^{\vee}; \qquad p \geq q
$$
where $V^{\vee}$ is the vector space dual to $V$.

Let us return to our situation. The two categories we have in mind are 
 the categories of mixed Hodge structures and  $l$-adic Galois modules equipped with a 
weight filtration.  
The functor $F$ assigns to an object of one of these categories the underlying vector space. 
We define a framing 
on $H^0(X, {\cal F}^*_{A,B})$ as a triple 
\begin{equation} \label{11.14.01.10q}
(H^0(X, {\cal F}^*_{A,B}), [\omega_A], [\Delta_B])
\end{equation}
where the classes
$$
[\omega_A] \in {\rm Gr}^W_{n + {\rm dp}({A})}H^0(X, {\cal F}_{A}); \qquad [\Delta_B] \in 
\Bigl({\rm Gr}^W_{n -{\rm dp}({B})}H^0(X, {\cal F}_B)\Bigr)^{\vee}
$$
are transformed  via canonical isomorphisms (\ref{11.19.01.6})-(\ref{11.19.01.7}) 
 to 
distinguished classes 
$$
[\omega_A] \in {\rm Gr}^W_{n + {\rm dp}({A})}H^0(X, {\cal F}^*_{A, B}); \qquad [\Delta_B] \in 
\Bigl({\rm Gr}^W_{n -{\rm dp}({B})}H^0(X, {\cal F}^*_{A,B})\Bigr)^{\vee}
$$
We define a framing on  $H^0(X, {\cal F}^!_{A,B})$ using the same two  classes. 
So the map $c$ from lemma \ref{11.7.01.1} 
is an equivalence of the framed objects corresponding to 
$H^0(X, {\cal F}^*_{A,B})$ and $H^0(X, {\cal F}^!_{A,B})$.

{\bf Remark}. Both $X-A$ and $X-B$ are affine varieties, so 
$H^i_{DR}(X-A) = H^i\Omega_{\rm log}^{\bullet}(X-A)$. Thus the class $[\omega]$ 
is represented by a  form $\omega \in \Omega_{\rm log}^{n}(X-A)$. 
So the equivalence class of the framed objects (\ref{11.14.01.10q}) always 
encodes the properties of the period $\int_{\Delta}\omega$, where 
$\Delta$ is a relative cycle. 

{\bf Remark}. If $X$ is a regular projective curve then $H^0(X, {\cal F}^*_{A,B})$ 
is a $1$-motive in the sense of Deligne [D4]. 

\section{Motivic setting}

Let $\pi: X \lra {\rm Spec}(F)$ be the canonical projection. 
In this section we show that, assuming   
 a condition on  $D$, the objects 
$\pi_* {\cal F}^{\bullet}_{A,*,B}$ and $\pi_* {\cal F}^{\bullet}_{A,!,B}$ in 
$D^b_{\rm sh}({\rm Spec}(F))$ 
are of geometric origin. 
If $F$ is a number field and 
the divisor $D = A \cup B$ provides a Tate stratification of $X$ we define 
$$
H^0(X; {\cal F}^{\bullet}_{A,*,B}) \quad \mbox{and} \quad H^0(X; {\cal F}^{\bullet}_{A,!,B})
$$
as objects of the abelian category ${\cal M}_T(F)$ 
of mixed Tate motives over $F$ 

{\bf 1. A blow up theorem}. Let $D$ be a divisor on $X$. 
We assume that all irreducible components of $D$ are divisors, and 
the stratification defined by $D$  is regular.

By a {\it sequence of centered at strata blow ups} we understand the following procedure. 
Blow up the closure of a stratum $S$ of $D$, getting a divisor $D_1$ 
on a variety $X_1$. Then blow up the  closure of a stratum $S_1$ of $D_1$, 
getting a divisor $D_2$ 
on  $X_2$, ..., blow up the closure of  a stratum $S_{m-1}$ of 
$D_{m-1}$, 
getting a divisor $D_m$ 
on a variety $X_m$.

\begin{theorem} \label{12.27.01.1} Suppose that $D$ is a divisor 
 locally isomorphic to a union of hyperplanes. Then there exists a centered at strata 
sequence of blow ups such that the 
last blow up delivers a normal crossing divisor. 
\end{theorem}

{\bf Proof}. A stratum $S$ of $D$ is called a good stratum 
 if  $D$ is a normal crossing 
divisor near $S$. 
Otherwise it is  a bad stratum. 

\begin{construction} \label{12.31.01.100} 
Let us blow up all bad vertices of $D$, then closures of preimages 
of all bad $1$-dimensional strata, after this closures of preimages 
of all bad $2$-dimensional strata, and so on. Let $\widehat X$ be the 
resulting variety and  $\widehat D $ the 
corresponding divisor 
on $\widehat X$.
\end{construction}

\begin{lemma} \label{ooppl}
$\widehat D$ a normal crossing divisor on $\widehat X$. 
\end{lemma} 

{\bf Proof}.
We use the induction on  ${\rm dim}X$. 
Let $p: \widehat X \lra X$ be the natural projection and $x \in \widehat X $. 
Let $k$ be the dimension of the  stratum containing $p(x)$. Suppose that $k>0$. 
Since our configuration is locally isomorphic to a family of hyperplanes 
we have a product situation near $p(x)$. Precisely, 
 there exists a regular divisor $D_Y$ in a variety $Y$ 
which is locally isomorphic to a family of hyperplanes, such that  
near $p(x)$ the pair $(D, X)$ is locally isomorphic to the one 
$({\Bbb A}^k \times D_Y, {\Bbb A}^k \times Y)$.  
Moreover the pair $(\widehat D, \widehat X)$ near $x$ 
is locally isomorphic to  $({\Bbb A}^k \times \widehat D_Y, 
{\Bbb A}^k \times \widehat Y)$ where $(\widehat D, \widehat Y)$ 
is obtained from $(D, Y)$ by our construction. Therefore by induction 
 $\widehat D$ a normal crossing divisor 
near such $x$. 

Suppose now that $k=0$, i.e. $p(x)$ is bad vertex of $D$. We blow up $p(x)$. 
Let $P$ be the special fiber of the blow up. It intersects transversally the strict preimages 
of other strata of $D$. 
Intersecting $P$ with them  
we get a divisor 
$D'$ in 
$P$. It is isomorphic to a family of hyperplanes in $P^{n-1}$. 
Now continue our construction, i.e. blow up other bad vertices, then bad edges, and so on,  
and see how it affects neighborhood of 
 $P$ and its preimages after the  blow ups. The preimage $\widehat P$ 
of  $P$ i`n $\widehat X$ is obtained by applying our construction 
to the divisor  
$D'$ in 
$P$. Therefore the divisor $\widehat D_P$ in $\widehat P$ 
defined as the intersection of $\widehat P$ with the closure of $\widehat D - \widehat P$ 
 is a 
normal crossing divisor by the induction assumption. 
The closure of the divisor $\widehat D - \widehat P$ 
intersects $\widehat P$ transversally, providing a bijective correspondence 
$S \lra S \cap P$ between the proper 
strata of $\widehat D_P$ and  the ones of $\widehat D$ intersecting $\widehat P$. 
 It follows that 
$\widehat D$ is a normal crossing divisor near $\widehat P$. 

The lemma and hence the theorem are proved.

{\bf 2. A comparison theorem.}  Suppose that $(A,B)$ is an admissible pair of divisors on $X$.  
Let us blow up the closure of a stratum $S$ of $D:= A \cup B$. 
We get a new variety $\widehat X$ equipped with a projection $p:\widehat X \lra X$. 
Then the preimage $S_1$ of $S$ 
and strict preimages of the irreducible components of the divisor $D$ 
are the irreducible components of the divisor $D_1 = p^{-1}D$. Let us assign the 
$A$- and $B$- labels to the irreducible components of $D_1 $, presenting it as  
$D_1 = A_1 \cup B_1$. The strict preimages of 
the irreducible components of  $D$ inherit their  labels 
from  $D$. So we need only to  label  
 $S_1$. We declare it an $A$- (resp. $B$-) component if $S$ is 
a pure $A$- (resp. $B$-) stratum. If $S$ is a mixed stratum we have a freedom to 
 declare it 
either $A$- or $B$-component. 
It is easy to see that for either choices $(\widehat A, \widehat B)$ is 
 an admissible pair of divisors.  
If we  choose always 
to declare blow ups of mixed strata as $A$-components, we denote by 
$(\widehat A^*, \widehat B^*)$ the final pair of divisors on  $\widehat X$. 
If we  always 
declare blow ups of mixed strata as $B$-components, denote the final pair of divisors by 
$(\widehat A^!, \widehat B^!)$. 

\begin{lemma} \label{12.31.01.1} Suppose that $A \cup B$ is a normal crossing divisor. 
Then there are canonical isomorphisms
$$
{\cal F}^{\bullet}_{A, !, B} \stackrel{\sim}{=}
{\cal F}^{\bullet}_{A, *, B} \stackrel{\sim}{=} {\cal F}^{*}_{A, B}
\stackrel{\sim}{=} {\cal F}^{!}_{A, B}
$$
\end{lemma}

{\bf Proof}. It is  clear from (\ref{1.21.01.8}) - (\ref{1.21.01.8se}) that these objects are perverse sheaves. It remains to show that extending to a mixed stratum by 
 $j_!$ we get the same object as when we extend by  $j_*$. 
Indeed, let $A_1$ be one of the $A$-components and 
 $j_1: X - A_1 \hra X $. Then if ${\cal F}_{A, B}$ is ${\cal F}^{!}_{A, B}$ or ${\cal F}^{*}_{A, B}$
one has 
$$
j_{1*} {\rm Res}_{X - A_1} {\cal F}_{A, B} = {\cal F}_{A, B}
$$ 
and a similar statement is valid for $j_{1!}$. 
The lemma follows by induction. In fact for a normal crossing divisor the extension
 of the constant sheaf on $X - (A \cup B)$ by $j_*$ to the $A$-components and by 
$j_!$ to the $B$-components does not depend on the order of the extensions.

Let $\widehat D = \widehat A \cup \widehat B$ be a 
splitting on $A$- and $B$-components obtained by the described above procedure, 
where the blow ups of mixed strata are labeled 
any way we want. Then 
$\widehat A \subset \widehat p^{-1} A$ and 
$\widehat B \subset \widehat p^{-1} B$. 
Therefore there are canonical 
morphisms
$$
{\cal F}^*_{\widehat  A} \lra {\cal F}^*_{\widehat p^{-1} A}, \qquad 
{\cal F}^!_{\widehat p^{-1} B}\lra {\cal F}^!_{\widehat  B}
$$
(We put $!$ and $*$ over the ${\cal F}$-sheaves to distinguish between  the extension 
functors, i.e. $j_!$ or $j_*$, used to define them. So ${\cal F}^*_{A} = {\cal F}_{A}$ 
and ${\cal F}^!_{A} = {\cal F}_B$.)

Since $\widehat p$ is proper,  we get canonical maps
$$
H^0(\widehat X, {\cal F}^*_{\widehat  A}) \lra 
H^0(\widehat X, {\cal F}^*_{\widehat  p^{-1} A}) = H^0(X, {\cal F}^*_{A})
$$
$$
H^0(\widehat X, {\cal F}^!_{\widehat p^{-1}B})
= H^0(X, {\cal F}^!_{B}) \lra H^0(X, {\cal F}^!_{\widehat  B})
$$
Applying ${\rm Gr}^W$, and then dualizing the second map,  we get morphisms
$$
\alpha_{\widehat  A}:  
{\rm Gr}^W_{n +{\rm dp}({A})}(X, {\cal F}^*_{\widehat  A})  
\lra 
{\rm Gr}^W_{n +{\rm dp}({A})}H^0(X, {\cal F}^*_{A})^{\vee} 
$$
$$
\beta_{\widehat  B}:  
\Bigl({\rm Gr}^W_{n -{\rm dp}({B})}(X, {\cal F}^!_{\widehat  B})\Bigr)^{\vee}  
\lra 
\Bigl({\rm Gr}^W_{n -{\rm dp}({B})}H^0(X, {\cal F}^!_{B})\Bigr)^{\vee} 
$$

Consider a framing on 
$H^0(X; {\cal F}^*_{A,  B})$ provided by nonzero classes 
$$
[\omega_{A}] \in 
{\rm Gr}^W_{n + {\rm dp}({A})}H^0(X, {\cal F}_{A}); \qquad 
[\Delta_{B}] \in  
\Bigl({\rm Gr}^W_{n -{\rm dp}({B})}
H^0(X, {\cal F}_{B})\Bigr)^{\vee}
$$

Denote by ${\cal F}_{A, B}$ the perverse sheaf provided by lemma \ref{12.31.01.1}. 
If  $\widehat D = \widehat A^* \cup \widehat B^* = \widehat A^! \cup \widehat B^! $ 
is a normal crossing divisor then  lemma \ref{12.31.01.1} provides perverse sheaves 
 ${\cal F}_{\widehat A^*,  \widehat B^*}$ and 
 ${\cal F}_{\widehat A^!,  \widehat B^!}$.

\begin{theorem} \label{10.3.01.1b} i) Assume that construction \ref{12.31.01.100} 
delivers a normal crossing divisor (e.g. $D$ is locally isomorphic to a family of hyperplanes).
Then there are canonical isomorphisms
\begin{equation} \label{12.30.01.1}
\widehat p_*  {\cal F}_{\widehat A^*,  \widehat B^*} = {\cal F}^{\bullet}_{A, *, B}; \qquad 
 \widehat p_* {\cal F}_{\widehat A^!,  
\widehat B^!} = {\cal F}_{A, !, B}^{\bullet}
\end{equation}
Moreover the maps $\alpha_{\widehat A^*}$, $\beta_{\widehat B^*}$, 
as well as $\alpha_{\widehat A^!}$, $\beta_{\widehat B^!}$ are isomorphisms.  
 
ii)  There is an  equivalence of framed objects:
\begin{equation} \label{12.30.01.11}
\Bigl(H^0(\widehat X; {\cal F}_{\widehat A^*,  \widehat B^*}),  
\alpha^{-1}_{\widehat A^*}[\omega_{A}], 
 \beta^{-1}_{\widehat B^*}[\Delta_{B}]\Bigr) 
\sim \Bigl(H^0(X; {\cal F}^*_{A,  B}),  [\omega_{A}], [\Delta_B]\Bigr)
\end{equation}
iii) There is a similar statement  with $*$ changed to  $!$. 
\end{theorem}

{\bf Proof}. i) Thanks to the construction for each 
bad stratum of $D$ there is a unique irreducible component of the divisor $\widehat D$ projecting
 onto the closure of this stratum. 
Let us show that if there is an isomorphism over $U(k)$
$$
{\rm Res}_{U(k)}\widehat p_*{\cal F}_{\widehat A^*,  \widehat B^*} = {\rm Res}_{U(k)}{\cal F}^{\bullet}_{A, *,  B} 
$$
then the same is true over $U(k+1)$. Set $\widehat X_k:= \widehat p^{-1}U(k)$

Let $\{\widehat D_i\}$ be the irreducible components of $\widehat D$ 
projecting to the bad strata 
in $U(k+1) - U(k)$, and $\{\widehat D_i^0\}$ their restrictions to $\widehat X_{k+1}$. 
Take one of them,  $\widehat D_1^0$. Suppose it is an $A$-component. 
Let $$
\widehat j_1: \widehat X_{k+1} - \widehat D^0_1  \hra \widehat X_{k+1}; \qquad 
j_1: U(k+1) - \widehat p(\widehat  D^0_1) \hra U(k+1)
$$ 
Since $\widehat D$ is a normal crossing divisor we can apply lemma \ref{12.31.01.1}. So 
if ${\cal F}$ is either 
${\cal F}_{\widehat A^*,  \widehat B^*}$ or ${\cal F}_{\widehat A^!,  \widehat B^!}$
there is an isomorphism
$$
\widehat j_{1*}{\rm Res}_{\widehat X_{k+1} - \widehat D^0_1}
{\cal F} = {\rm Res}_{\widehat X_{k+1}}{\cal F}
$$
One has 
$\widehat p \circ \widehat j_1 = j_1 \circ \widehat p$. Thus 
 $\widehat p_* \widehat j_{1*} = j_{1*} \widehat p_*$. Therefore 
$\widehat p_*{\rm Res}_{\widehat X_{k+1}}{\cal F} = 
j_{1*}\widehat p_*{\rm Res}_{\widehat X_{k+1} - D_1^0}{\cal F}$. The
case when $\widehat D^0_1$ is a $B$-component is similar. 
Since the  divisors $\widehat D_i^0$ are disjoint we immediately get a similar statement 
for $\cup_i\widehat D_i^0$. 

The second statement follows from isomorphisms (\ref{12.30.01.1}) and 
theorem \ref{11.5.01.1}b). 
The part i) of the theorem is proved. 

ii) The part i) implies the equivalence 
$$
\Bigl( H^0(X; {\cal F}^{\bullet}_{A, *, B}), [\omega_A], [\Delta_B]\Bigr) \sim 
\Bigl( H^0(\widehat X; {\cal F}_{\widehat A^*, \widehat  B^*}), 
\alpha^{-1}_{\widehat A^*}[\omega_{A}], \beta^{-1}_{\widehat B^*}[\Delta_{B}]\Bigr)
$$ 
Recall that the frame morphisms for  ${\cal F}^*_{A,  B}$ were constructed 
by applying the functor $H^0_{\tau}$ to the frame morphisms 
${\cal F}_{B} \lra {\cal F}^{\bullet}_{A, *, B}\lra {\cal F}_{A}$. 
So the canonical morphism in the derived category
$
H^0_{\tau}: {\cal F}^{\bullet}_{A, *, B}\lra {\cal F}^*_{A,   B} 
$ 
provides a commutative diagram where the horizontal  arrows are the frame morphisms:
$$ 
\begin{array}{ccccc}
{\cal F}_{B}& \lra &{\cal F}^{\bullet}_{A, *, B}&\lra &{\cal F}_{A}\\
\downarrow =&&\downarrow &&=\downarrow  \\
{\cal F}_{B}&\lra &{\cal F}^*_{A,   B} & \lra &{\cal F}_{A} \end{array}
$$
Applying the functor $H^0(X;-)$ to it we get an equivalence of framed objects
$$
\Bigl( H^0(X; {\cal F}^{\bullet}_{A, *, B}), [\omega_A], [\Delta_B]\Bigr) \sim 
\Bigl( H^0(X; {\cal F}^*_{A,  B}), [\omega_A], [\Delta_B]\Bigr)
$$
The part ii) is proved. 
The part iii) is very similar. 
The theorem is proved.

{\bf 3. The mixed Tate motive corresponding to $H^0(X; {\cal F}_{A,  B})$}. 
\begin{proposition} \label{1.1.02.1}
Suppose that $F$ is a number field, and 
$A \cup B$ is a normal crossing divisor on a regular variety $X$ over $F$, 
providing a Tate stratification of $X$. Then there exists a mixed Tate motive 
$m(X; A, B)$ whose Hodge and l-adic realizations are given by 
$H^0(X; {\cal F}_{A,B})$. 
\end{proposition}

{\bf Proof}. Suppose that $X$ is a regular variety over any field $F$ 
and $Y$ is a normal crossing divisor on $X$. 
Then  there is an object $m(X,Y)$ of the triangulated category of mixed motives 
over $F$ 
such that realizations of $m^{\bullet}(X,Y)$ are given by the relative cohomology groups 
$H^*(X,Y)$. (I use a cohomological version of Voevodsky's category ${\cal D}{\cal M}_F$). 
Namely, let $\{Y_i\}$, $i \in \Lambda$, 
 be the set of all irreducible components of the divisor $Y$. 
Set $Y(k):= \cup_{|I|=k}Y_I$ and consider the following complex of regular varieties where 
$X$ is in the degree zero, and the differentials are the ones 
in a simplicial resolution:
$$
m^{\bullet}(X,Y):= ... \lra Y(2) \lra Y(1) \lra X
$$
Let ${\cal D}_T(F)$ be the full subcategory of ${\cal D}{\cal M}_F$ 
consisting of the objects glued from the the pure Tate motives $\Q(n)$ and their shifts. 
If $Y$ provides a Tate stratification of a Tate variety $X$ then 
$m^{\bullet}(X,Y)$ is an object of ${\cal D}_T(F)$. 
Finally, if $F$ is a number field then there is a canonical $t$-structure $t$ 
on the category ${\cal D}_T(F)$ providing an abelian category ${\cal M}_T(F)$ of mixed 
Tate motives over $F$ (see ch. 5 of [G9]). Set $m^i(X,Y):= H_t^i(m^{\bullet}(X,Y))$.

Applying this construction to  $Y:= B_A:= B - B \cap A$
we get mixed Tate motives 
$m^i(X; A, B_A)$. 
Then  $m^{-{\rm dim}X}(X; A, B_A)$ is just what we need. The proposition is proved.

{\bf 4. Remarks on the frame classes}. Recall a splitting on $A$- and $B$-components
 $\widehat D = \widehat A \cup \widehat B$ 
 described in s. 3.2, 
where the blow ups of mixed strata are labeled 
any way we want. For applications of theorem 
\ref{10.3.01.1b} it is useful to have a more direct description of 
the classes $\alpha^{-1}_{\widehat A}[\omega_{A}]$ and 
$\beta^{-1}_{\widehat B}[\Delta_{B}]$. 

\begin{lemma} \label{2.16.02.11} 
Suppose that $A$ is a normal crossing divisor and $\omega_A \in \Omega^n_{\rm log}(X-A)$. Then 
$\widehat p^{*}\omega_A \in \Omega^n_{\rm log}(\widehat X-\widehat A)$. 
Therefore  $\alpha^{-1}_{\widehat A}[\omega_{A}]$ is represented by 
 $[\widehat p^*\omega_{A}]$. 
\end{lemma}

{\bf Proof}. A priori the form $\widehat p^{*}\omega_A$ might have logarithmic singularities 
at $\widehat p^{-1}A$. The following general lemma tells us that 
 $\widehat p^*\omega_A$ can have singularities only 
at the blow ups of pure $A$-strata. 

\begin{lemma} \label{1.23.01.6} i) Let  $Y$ be a normal crossing divisor
  in a regular variety   $X$, and $\omega \in \Omega^n_{\rm log}(X-Y)$. 
Let $p:\widehat X \lra X$ be the  blow up of an
irreducible subvariety $Z$. 
Suppose that the generic point of $Z$ is different from the  generic points of  strata 
of $Y$. Then  $p^*\omega$ does not have a singularity 
at the special divisor of $\widehat X$. 

ii) Let $F \in K_n^M(F(X))$. Assume that $F$ has non zero 
residues  only at the components of a normal crossing divisor $Y$. 
  Then under the same conditions as in i) the residue of  $p^*F$ 
at the special divisor of the blow up is zero.
\end{lemma}

{\bf Proof}. i) Considering  the generic point of
$Z$ we reduce the problem to the case when $Z$ is a point. So let $x_1, ..., x_n$ be 
local coordinates such that $\omega = d \log (x_1) \wedge ... \wedge d
\log (x_k) \wedge \Omega$ where $\Omega$ is a regular form of positive
degree near the point
$Z = (0, ...,0)$. Let $u_1:= x_1$, $u_i: = x_i/x_1$ be the local 
coordinates on the blow up. Then $x_1 =0$ is a local equation of the special divisor. 
One has 
$$
p^* \omega = d \log (x_1) \wedge d
\log (u_2) \wedge ... \wedge d \log (u_k) \wedge p^* \Omega
$$
Since $d \log (x_1) \wedge p^*f\cdot  d(x_1 u_{k+1}) = p^*f \cdot d x_1
\wedge du_{k+1}$ where $f$ is regular near $Z$,  the statement follows. 

ii) Similar to the proof of i). Lemma \ref{1.23.01.6}, and 
hence  lemma \ref{2.16.02.11} are proved.

Now look at the similar question for the relative 
cycle $\Delta_B$ in the Betti realization. A simple geometric argument shows that 
there is an open part $\Delta^0_B \subset X-B$  of the cycle $\Delta_B$ such that 
the class $\beta^{-1}_{\widehat B}[\omega_{B}]$ can be represented by the closure of $\Delta^0_B$ 
in 
$\widehat X$.  
This is a homological version of lemma \ref{2.16.02.11}. 

{\bf 5. An unramifiedness criteria}.  Let $R$ be 
a discrete valuation ring with the fraction field $K$,  
${\cal M}$ the maximal ideal in $R$, 
 $k$ the residue field of $R$, 
and $\overline k$ its algebraic closure. For a scheme $X$ over a field $F$ set 
 $\overline X:= X \otimes_F \overline F$.

Suppose that $D$ is a proper normal crossing divisor 
in a proper regular scheme $X$ over $R$. 
Let $A$ and $B$ be unions 
of the irreducible components of $D$. We will assume that no irreducible 
components is shared by  both of them. Then there is the 
${\rm Gal}(\overline K/K)$-module 
\begin{equation} \label{2.24.02.2}
H^n_{\rm et}(\overline X- \overline A, \overline B_A; \Q_l); 
\quad B_A:= B - (A \cap B); 
\end{equation}
 We want to have a criteria for this module to be unramified.

\begin{definition} \label{3.15.02.1} Let $D$ be a normal crossing divisor 
in a  regular scheme $X$ over $R$. Assume that the pair 
$(D,X)$ is  proper over $R$. 
We say that reduction modulo ${\cal M}$ does not change the combinatorics of $(D,X)$ if 
$X$ and every stratum of $D$ are smooth over $R$, and 
the reduction map from the strata of $D$ to ones at the special fiber 
is a bijection.
\end{definition}

To check that a scheme is smooth over $R$  the following 
result is useful,  see 
proposition 3.24 in chapter I of  [M]. 
Let $Y$ be an integral scheme over 
$R$. Denote by $Y^0$ and $Y^{\eta}$ its special and generic fibers. 
Then $Y$ is smooth over $R$ if and only if the generic fiber is non empty and 
$Y^0(\overline k)$ 
and $Y^{\eta}(\overline K)$ have no singular points, 

\begin{proposition} \label{2.24.02.1} Suppose that the reduction 
modulo ${\mathcal M}$ does not change the combinatorics of $(D,X)$. Then, assuming 
$l$ is prime to the characteristic of $k$,  there is an isomorphism
\begin{equation} \label{2.24.02.3}
H^n_{\rm et}(\overline X-\overline A, \overline B_A; \Q_l)  \stackrel{\sim}{\lra}
 H^n_{\rm et}({\overline {X^0} -\overline {A^0}}, \overline {B_A^0}; \Q_l)
\end{equation}
In particular  the 
${\rm Gal}(\overline K/K)$-module (\ref{2.24.02.2}) is unramified. 
\end{proposition}

{\bf Proof}. 
By the smooth and proper base change theorem ([M]) 
for any stratum $Y$ of $X$, including $X$ itself, 
there is  an  isomorphism
$$
H^n_{\rm et}(\overline Y; \Q_l)  
\stackrel{\sim}{\lra}
H^n_{\rm et}(\overline {Y^0}; \Q_l)
$$
and hence $H^n_{\rm et}(\overline Y; \Q_l)$ is unramified. We calculate 
$H^n(X-A,B_A)$ by replacing the pair $(X-A, B_A)$ by its 
standard ([D4]) simplicial resolution $S_{\bullet}(X-A, B_A)$. The $i$-simplices 
of this simplicial scheme are given by  disjoint union of 
the $i$-fold intrersections of the irreducible components of  $B_A$. 
The result follows from the following two observations:

i) There is a well defined reduction modulo $\cal M$ 
of the simplicial scheme $S_{\bullet}(X-A, B_A)$. 

ii) The reduction modulo ${\cal M}$ induces an isomorphism on the \'etale cohomology 
for each of the schemes 
$S_{i}(X-A, B_A)$.

Indeed,  assuming i), let $S^0_{\bullet}(X-A, B_A)$ be the simplicial scheme 
obtained by reduction of the simplicial scheme 
$S_{\bullet}(X-A, B_A)$. In particular each component $S^0_{i}(X-A, B_A)$ is obtained by 
reduction of the corresponding component $S_{i}(X-A, B_A)$. 
Then there is  a map of the standard spectral sequence provided by the simplicial resolution   
which computes the left hand side 
of (\ref{2.24.02.3}) to the one computing the right hand side.
Thanks to ii) the map induces an isomorphism of the $E_1$ terms of these 
spectral sequences, so the proposition follows.  

The check the statement  i) observe that we can always reduce modulo ${\cal M}$  
each of the components of the simplicial scheme $S_{\bullet}(X-A, B_A)$. 
The maps between them are defined since the reduction modulo ${\cal M}$ does not change the combinatorics of 
$D$, and in particular  no stratum of $\overline B_A$ is supposed to be mapped to a stratum which 
has been removed. 

The statement  ii) follows from the following lemma

\begin{lemma} \label{2.24.02.5} Let $X$ be a proper scheme over $R$ and $B$ is a proper 
normal crossing divisor 
in $X$. Suppose that reduction modulo ${\cal M}$ does not change the combinatorics of $B$. 
Then there is an isomorphism
$$
H^n_{\rm et}(\overline X- \overline B; \Q_l) 
\stackrel{\sim}{\lra}H^n_{\rm et}({\overline {X^0}}- \overline {B^0}; \Q_l)
$$
\end{lemma}

{\bf Proof}. By duality it is sufficient to prove a similar isomorphism for 
$H^n(X,B)$. Then we replace $(X,B)$ by its simplicial resolution $S_{\bullet}(X,B)$. 
Since $X$ is proper this
 simplicial resolution always  has a reduction modulo ${\cal M}$ denoted 
$S^0_{\bullet}(X,B)$. 
Observe that in general $H^n(\overline {S^0_{\bullet}}(X,B); \Q_l)$ 
does not necessarily isomorphic to $H^n(\overline {X^0}, \overline {B^0}; \Q_l)$ (e.g., take 
 $X = {\Bbb A}^1$, $B = \{0\} \cup \{p\}$, and ${\cal M} = (p)$). 
However we do have the isomorphism if 
  $S^0_{\bullet}(X,B) = S_{\bullet}(X^0, B^0)$. This 
 is the case in our situation since the reduction does not change the combinatorics of $B$.
The lemma and hence the proposition are proved. 

{\bf Example}. Let $p$ be a prime number. 
Then  $H^1(P^1 - \{0, \infty\}, \{1, p\})$ is unramified outside of $p$. 
The simplicial scheme corrersponding to this cohomology group is 
$
\{1\} \cup \{p\} \lra P^1 - \{0, \infty\}
$. 
Observe that we can not define its reduction modulo $p$ since  $\{p\}$  should be mapped to $0$, which is not in $P^1 - \{0, \infty\}$!

\begin{corollary} \label{2.24.02.6} Suppose that we are in the situation described 
in proposition \ref{1.1.02.1}. 
Assume further that $(X, A, B)$ is defined and  proper over a ring of $S$-integers 
${\cal O}_S$ in the  number field $F$, and ${\cal P}$ is a prime ideal of $F$ outside of $S$. 
Then if  the reduction modulo ${\cal P}$ does not change the combinatorics 
of the divisor $A \cup B$, then the mixed Tate motive $m(X; A,B)$ is unramified at ${\cal P}$.
\end{corollary} 

{\bf Proof}. A mixed Tate motive over $F$ is unramified at 
${\cal P}$ if and only if its $l$-adic realization, where $l$ is prime to ${\cal P}$,  is 
unramified at ${\cal P}$. So the corollary follows from proposition \ref{2.24.02.1}. 



\section{A specialization theorem}

{\bf 1. The specialization functor [Ve]}. Let $Z$ be a regular 
subvariety of a regular 
variety $X$. 
Let  $N^0_ZX$ be the normal bundle  to  $Z$ in $X$, with the zero section
  removed. Set $U:= X-Z$. 
Recall the specialization functor  
$$
 {\rm Sp}: D^b_{\rm Sh}(U) \lra D^b_{\rm Sh}(N^0_ZX)
$$
In fact Verdier defined it for any pair of subvarieties 
$Z \subset X$, with $N^0_ZX$
replaced by the normal cone, and $D^b_{\rm Sh}(U)$ by $D^b_{\rm Sh}(X)$. Extending elements of $D^b_{\rm Sh}(U)$ by $j_!$ to $X$ 
and applying the Verdier's definition we come to the functor above.

The functor $ {\rm Sp}$  sends perverse
  sheaves to perverse sheaves and  commutes with the duality
  functor on $D^b_{\rm Sh}$. 

Since $Z$ is a regular subvariety of a regular variety $X$ the
specialization functor $ {\rm Sp}$ provides a
Tate functor between the corresponding mixed Tate categories on $U$
and $N^0_ZX$ 
 of unipotent variations of Hodge-Tate structures or lisse Tate l-adic sheaves.    
Indeed, $ {\rm Sp}$ is  an exact tensor functor, and 
if $Z$ is regular subvariety it obviously transforms the Tate objects,
which are the constant sheaves on $U$ tensored by $\Q(n) $ (resp $\Q_l(n)$)
and shifted to the left by ${\rm dim}U$, to the
Tate objects.

{\bf 2. A specialization theorem}. 
 Recall that  $D_i$,  $i \in \Lambda$, is the set of irreducible 
components of the divisor $D$, and for any subset $I \subset \Lambda$ we set 
$D_I:= \cap_{i \in I} D_i$. In particular $D_{\emptyset} = X$. 
The subvarieties $\{D_I\}$ are the closures of the 
strata of the stratification defined by $D$.

{\it Excellent variations}. We say that  $\{D_i(t)\}$ is a smooth variation 
 of regular 
divisors on $X$ 
parametrized by a regular base $T$ if  there are 
 regular divisors ${\cal D}_i \subset X\times T$ 
smooth with respect to the projection $p: X\times T \lra T$, so that 
 $D_i(t) = {\cal D}_i(t) \cap p^{-1}t$.  

\begin{definition} \label{11.2.01.5}
A smooth variation of the divisors $\{D_i(t)\}$  is an {\rm excellent} variation if 
for any   $I\subset \Lambda$

1)  ${\rm dim}D_I(t)$ is the same for all $t\in T$.

2)  $D_I(t)$ are regular irreducible projective subvarieties. 

3) The family  $\{D_I(t)\}$ extends to a smooth 
family of regular subvarieties parametrized by a projective 
variety $B_I$ containing $T$. 
\end{definition}

Replacing $T$ by a 
non empty Zariski open subvariety $T^0 \subset T$ 
we may always assume that 1).

{\bf Example}. $X = P^n$; $D$ is a union of hyperplanes. 
Any base of  deformation of hyperplanes $D_i$ satisfying condition 1) 
 satisfy also 
conditions 2) and 3). Indeed,  $D_I(t)$ is a plane in $P^n$, 
so one can take $B_I$ to be 
the corresponding Grassmannian.

Now suppose we have a variation of the divisors $D_i(\varepsilon)$ 
in $X$ parametrized by a 
one dimensional regular base $\Sigma$ containing a point 
$0$.  Let $\Sigma^0:= \Sigma -\{0\}$. 

Let $\Lambda = I \cup J$. 
Set 
$$
A(\varepsilon):= \cup_{i\in I}D_i(\varepsilon), \qquad 
B(\varepsilon) := \cup_{j \in J}D_j(\varepsilon)
$$
 and $(A,B):= (A(0), B(0))$. 
Assume that

i) The pair $(A,B)$ is admissible.

ii) $\{D_i(\varepsilon)\}$ 
is an excellent variation over $\Sigma^0$. 

Then 
 replacing the base by a non empty open subset of $\Sigma$ we may assume that 
 $(A(\varepsilon), 
B(\varepsilon))$ is an admissible pair for all $\varepsilon \in \Sigma$. 

Denote by ${\cal A}$ and ${\cal B}$ the divisors on $X \times \Sigma^0$ 
whose fibers over $\varepsilon \in \Sigma^0$ are $A(\varepsilon)$ and $B(\varepsilon)$. 
So there are the corresponding objects denoted 
 $\widetilde {\cal F}^{\bullet}_{\cal A,*,\cal B}$ and 
$\widetilde {\cal F}^*_{\cal A,\cal B}$  on 
$X \times \Sigma^0$. The tilde  
emphasizes that we are working over $\Sigma^0$.

We are going to apply the specialization functor  to the divisor 
\begin{equation} \label{4.6.01.1}
 i_0: X \times 0 \hookrightarrow
  X \times \Sigma
\end{equation} 
 Let $v 
\in T_{0}(\Sigma) - 0$. Denote  
by $ {\rm Sp}_{s}({\cal G}) $ the restriction of   $ {\rm Sp}({\cal G})$ to
  the section   $s := p^{-1}(v)$ of the normal bundle, which is 
  identified with $X \times 0$ in (\ref{4.6.01.1}).

We claim that in both Hodge and $l$-adic  settings the specialization 
${\rm Sp}_s\widetilde {\cal F}^*_{\cal A, \cal B}$ satisfies the same conditions as 
in proposition \ref{10.24.01.2}: 

 \begin{theorem} \label{11.1.01.3} Assume the conditions i), ii) above. 
Then in both Hodge and $l$-adic settings 
  ${\rm Sp}_s\widetilde {\cal F}^*_{\cal A,\cal B}$ is glued from the objects $\delta_C(-m)$ 
where $0 \leq m \leq {\rm codim}C$, 
where  $C$ runs 
through the limiting strata of the 
stratification ${\cal D}_J(\varepsilon)$ as $\varepsilon \to 0$. 
\end{theorem}

{\bf Remark 1}. The limiting strata consist of  the stratum $D_J$ and perhaps some other 
subvarieties. For example consider degeneration of a pair (line, plane) in $P^3$ 
degenerating to a line sitting in the limiting plane. 

{\bf Remark 2}. A similar result is valid when ${\cal F}^*_{\cal A,\cal B}$ is replaced by 
$H^i({\cal F}^{\bullet}_{\cal A,*,\cal B})$

{\bf Proof}. {\it The Hodge setting}.  
We have a supply of Hodge sheaves $\delta_{{\cal D}_J}(-m)$. 
 For all but finitely many  points $\varepsilon \in \Sigma^0(\C)$ the restriction to 
$\varepsilon$ commutes with all the standard functors. 
Let $i_{\varepsilon}: \{\varepsilon\} \hra \Sigma^0(\C)$.  
It follows that for all but finitely many  points $\varepsilon$ the 
object $i^*_{\varepsilon}\widetilde {\cal F}^*_{\cal A, \cal B}[-1]$ is isomorphic 
to ${\cal F}^*_{A(\varepsilon), B(\varepsilon)}$, 
 and hence has all the properties 
listed in proposition \ref{10.24.01.2}. Observe  that  
 $i^*_{\varepsilon}\delta_{{\cal D}_J}(-m) = \delta_{D_J(\varepsilon)}(-m)[-1]$.

Denote by ${\cal P}_{\cal H}(X)$ the category of perverse Hodge sheaves on $X$.

 Choose a sufficiently small punctured analytic disc 
$U^* \subset \Sigma^0(\C)$ with the puncture at  $0$. For  a finite covering
$\pi: \widehat U^* \lra U^*$ consider the map  $$
\widehat \pi:= 
{\rm Id}\times \pi: X \times \widehat U^* 
\lra X \times U^*
$$ 
Let $\widehat {\cal D}$ be the variation of the divisors parametrized by 
$\widehat U^*$, obtained by the base change $\pi: \widehat U^* \lra U^*$. 

\begin{lemma} \label{11.1.01.4} There exist $J$, $m \leq {\rm codim}D_J$, and 
a finite  covering 
$\pi: \widehat U^* \lra U^*$  such that there is a non zero 
element in
$$
{\rm Hom}_{{\cal P}_{\cal H}(X \times \widehat U^*)}(\widehat\pi^*\widetilde {\cal F}^*_{\cal A,\cal B}, 
\delta_{\widehat {\cal D}_J}(-m))
$$
\end{lemma}

{\bf Proof}. Let 
$p: X \times U^* \lra U^*$ be the canonical projection.  One has 
\begin{equation} \label{11.25.01.42}
{\rm RHom}_{D^b_{\cal H}(X \times U^*)}(\widetilde {\cal F}^*_{\cal A,\cal B}, 
\delta_{{\cal D}_J}(-m)) = 
\end{equation}
$$
{\rm RHom}_{D^b_{\cal H}(X \times U^*)}(\delta_{X \times U^*}, 
\ast \widetilde {\cal F}^*_{\cal A,\cal B} \otimes^!
\delta_{{\cal D}_J}(-m)) = 
$$
\begin{equation} \label{11.25.01.41}
{\rm RHom}_{D^b_{\cal H}(U^*)}(\delta_{U^*}, p_*(
\ast \widetilde {\cal F}^*_{\cal A,\cal B} \otimes^!
\delta_{{\cal D}_J}(-m))) 
\end{equation} 
The Hodge sheaves
\begin{equation} \label{11.25.01.31}
R^ip_*(
\ast \widetilde {\cal F}^*_{\cal A,\cal B} \otimes^!
\delta_{{\cal D}_J}(-m))
\end{equation} 
are restrictions to $U^*$ of the similar ones on $\Sigma^0$, which are  
smooth at the generic point of $\Sigma^0$ (this follows from the 
corresponding well known fact about $D$-modules). Thus 
removing a finite set of points from  $\Sigma^0$ we may assume that 
(\ref{11.25.01.31}) is a variation of mixed Hodge structures on $\Sigma^0$, 
and hence on  $U^*$.  

Recall that  ${\rm Hom}_{{\cal P}_{\cal H}} = {\rm R^0Hom}_{D^b_{\cal H}}$.   
The monodromy  around $0$ provides an operator $N$  acting on the vector space 
\begin{equation} \label{11.26.01.1}
{\rm R^0Hom}_{D^b_{\cal H}(\C)}(\delta_{\varepsilon}[-1], i_{\varepsilon}^* p_*(
\ast \widetilde {\cal F}^*_{\cal A,\cal B} \otimes^!
\delta_{{\cal D}_J}(-m))),
\end{equation}
The $H^0$ of the complex 
(\ref{11.25.01.41}) is isomorphic to the invariants $N$ acting on (\ref{11.26.01.1}). 
Let us show that this vector space is non zero for some $J$ and $m \leq {\rm codim}D_J$. 
Indeed, since restriction to a generic point 
 $\varepsilon \in U^*$ commutes with all the standard functors, 
we can interchange $i_{\varepsilon}^*$ with $p_*$ and other functors involved in 
(\ref{11.26.01.1}). Then 
proceeding as above we rewrite (\ref{11.26.01.1}) as a vector space 
\begin{equation} \label{11.25.01.10} 
{\rm R^0Hom}_{D^b_{\cal H}(X \times \{\varepsilon\})}(
{\cal F}^*_{A(\varepsilon),B(\varepsilon)}, 
 \delta_{{ D}_J(\varepsilon)}(-m))  
\end{equation}
By proposition \ref{10.24.01.2} there exists $J$ and 
$m\leq {\rm codim}D_J(\varepsilon)$ such that (\ref{11.25.01.10})
 is non zero. 

The regularity 
theorem for D-modules (see lecture 4 in [Be]) implies that 
variation (\ref{11.25.01.31}) on $U^*$ has a quasiunipotent monodromy around $0$. 
Indeed, the variation was obtained as a composition of several 
functors applied to $\delta_U$. Since $\delta_U$ obviously has 
quasiunipotent spectrum in the terminology of [Be], and the standard functors 
preserve this property, 
the statement follows. 
Therefore the monodromy operator $N$ acting on (\ref{11.26.01.1}) is quasiunipotent. 
If $N$ is unipotent then we are done:  
the subspace of $N$-invariants on a  non zero vector space (\ref{11.26.01.1}) is non zero. 
Replacing $U^*$ by 
its finite cover $\pi: \widehat U^* \lra U^*$ we can make the monodromy $N$ unipotent. 
The lemma is proved. 

Having the lemma we get the theorem as follows. Let us suppose first that the monodromy $N$ 
is unipotent, so we do not have to pass to a  covering  $\widehat U^*$. Then we find by induction 
a filtration $F$ of the object  $\widetilde {\cal F}^*_{\cal A, \cal B}$ such that its 
 associated graded 
are isomorphic to $\delta_{{\cal D}_J}(-m)$ for some $J$ and $m \leq 
{\rm codim}{\cal D}_J$. Since $\{D_i(\varepsilon)\}$ 
is an excellent variation over $U^*$, each of the families 
$D_J(\varepsilon)$ extends to a smooth family on $U:= U^* \cup 0$. In particular 
there is a regular reduced subscheme $ D_J(0)$, the limit of the family $D_J(\varepsilon)$
as $\varepsilon \lra 0$. 
Therefore 
$$
{\rm Sp}_s\delta_{{\cal D}_J}(-m) = \delta_{D_J(0)}(-m)
$$
Since the specialization is an exact functor we get 
the Hodge version of theorem when $N$ is unipotent. 

Let us deduce the general case from the one when $N$ is unipotent. 
Consider the  diagrams
$$
\begin{array}{cccccccc}
X \times \{0\}& \stackrel{\widehat i}{\hookrightarrow} & X \times 
\widehat U &&& N^0_{X \times \{0\}}X \times \widehat U 
& = & X \times \C^*\\
\downarrow =&&\downarrow \widehat \pi &&&\downarrow  {\overline \pi}&& \downarrow  \\
X \times \{0\}&\stackrel{i}{\hookrightarrow} &X \times U &&& 
N^0_{X \times \{0\}}X \times  U
& = & X \times \C^*
\end{array}
$$
The fibers of the punctured normal bundle to $X \times \{0 \}$ are identified with 
$\C^*$, and  $\overline \pi$ is given by $z \lra z^{{\rm deg}\pi}$ 
on the fibers. 
Then, denoting by ${\rm Sp}_{i}$ (resp. ${\rm Sp}_{\widehat i}$) the specialization functor for the 
divisor given by $i$ (resp $\widehat i$) we have
\begin{equation} \label{st}
{\rm Sp}_{i}\circ {\widehat \pi}_*=  {\overline \pi}_* \circ {\rm Sp}_{\widehat i}
\end{equation} 
This 
 follows from the construction of specialization 
via the nearby cycles functor $\Psi$,
 using the fact that $\Psi$ commutes with proper direct images.  Observe that   
$\widetilde {\cal F}^*_{\cal A, \cal B}$ is a direct summand of 
$\widehat \pi_*\widehat  \pi^* \widetilde {\cal F}^*_{\cal A, \cal B}$. 
 So ${\rm Sp}_{i}(\widetilde {\cal F}^*_{\cal A, \cal B})$ is a direct summand of 
$$
{\rm Sp}_{i}(\widehat \pi_*\widehat  \pi^* \widetilde {\cal F}^*_{\cal A, \cal B}) 
\stackrel{(\ref{st})}{=}
{\overline \pi}_* {\rm Sp}_{\widehat i}(\widehat  \pi^* \widetilde {\cal F}^*_{\cal A, \cal B})
$$
Let $\pi^{-1}s = \{s_i\}$. Then restriction of the right object  to $X \times s$ is a direct sum of the restrictions  of 
${\rm Sp}_{\widehat i}(\widehat  \pi^* \widetilde {\cal F}^*_{\cal A, \cal B})$ to  $X \times s_i$. 
Each of them  satisfies the condition of the theorem. The Hodge version of the theorem is proved.

{\it The l-adic setting}. It follows the same pattern as the proof of the Hodge version. 
There is a perverse $l$-adic sheaf $\widetilde {\cal F}^*_{\cal A, \cal B}$ over 
$X \times \Sigma^0$.  
Choose an algebraic closure $\overline K$ of the field of functions 
$K:= \Q(\Sigma^0)$, 
and an embedding $K \hra \overline K$ providing a generic geometric  point 
of $\Sigma^0$. Denote by $i_{\overline K}: X_{\overline K} \lra X \times \Sigma^0$ 
the corresponding map.

Let $y$ be a generic geometric  point of $\Sigma^0$ providing map 
$ i_y: X_{\overline K} \lra 
X \times \Sigma^0$. Then  for any 
perverse sheaf ${\cal G}$ on  $X \times \Sigma^0$ 
the (shifted) restriction $i_y^* {\cal G}[-1]$   is a perverse $l$-adic sheaf on 
$X_{\overline K}$. 

In particular let $\overline A := i_{\overline K}^* {\cal A}$, 
$\overline B:= i_{\overline K}^* {\cal B}$ and 
$\overline D_J:= i_{\overline K}^* {\cal D}_J$. Then there is 
a perverse $l$-adic sheaf  ${\cal F}^*_{\overline A, \overline B}$ 
on $X_{\overline K}$ which by 
proposition \ref{10.24.01.2} is glued from $\delta_{\overline D_J}$. 
Since all the standard functors commute 
with restriction to generic geometric point, there is an isomorphism   
$$
{\cal F}_{\overline K}:= i_{\overline K}^*\widetilde {\cal F}^*_{\cal A, \cal B}[-1]
\stackrel{\sim}{=}{\cal F}^*_{\overline A, \overline B}
$$ 
So ${\cal F}_{\overline K}$
 is glued from 
$\delta_{\overline D_J}$. 
Therefore one can find $J$ such that 
\begin{equation} \label{11.26.01.2}
{\rm Hom}_{\overline K}({\cal F}_{\overline K}, \delta_{\overline D_J}) \not = 0
\end{equation}
 Here 
{\rm Hom} is  in the category of perverse $l$-adic sheaves over $\overline K$. 
Both perverse sheaves are equipped with an action of the Galois group 
${\rm Gal}(\overline K/K)$. So the {\rm Hom}
has a structure of a finite dimensional 
$l$-adic 
${\rm Gal}(\overline K/K)$-module. 

Let us show that, twisting (\ref{11.26.01.2}) by $\Q_l(-m)$ for some 
$m \leq {\rm codim}D_J$,  we can find there a non zero 
${\rm Gal}(\overline K/K)$-invariant vector. Indeed, 
thanks to (\ref{11.26.01.5}) the ${\rm Gal}(\overline K/K)$-module  (\ref{11.26.01.2}) 
 is glued from 
the Tate modules $\Q_l(n)$. 
It follows that for a certain integer $m$ 
there is a ${\rm Gal}(\overline K/K)$-invariant vector in 
\begin{equation} \label{11.29.01.1}
{\rm Hom}_{\overline K}({\cal F}_{\overline K}, \delta_{\overline D_J})(-m)
\end{equation}
Recall that 
 ${\cal F}_{\overline K}$ and $\delta_{\overline D_J}$ 
are isomorphic to the perverse sheaves obtained by restriction 
$i_{\overline K}^*[-1]$ from 
 the perverse sheaves ${\cal F}^*_{\cal A, \cal B}$ and $\delta_{{\cal D}_J}$ 
defined over $K$. Therefore the subspace of ${\rm Gal}(\overline K/K)$-invariants 
in (\ref{11.29.01.1}) is identified with  
\begin{equation} \label{11.29.01.5}
{\rm Hom}_{K}({\cal F}^*_{\cal A, \cal B}, \delta_{{\cal D}_J}(-m))
\end{equation}
Here 
{\rm Hom} is  in the category of perverse $l$-adic sheaves over $K$. 
It follows that there is a non zero element in (\ref{11.29.01.5}). 
It remains to show that $m \leq {\rm codim}D_J$. Unfortunately so far 
there is no formalism of weights over $\overline K$. We will 
circumvent this problem by  
using the reduction 
to a finite field, where we can use the weights. 
 Choose a  finitely generated ring  ${\cal O}$ over $\Z$ 
such that $X,  \{{\cal D}_J\}$ are defined over ${\cal O}$. 
The restriction  to a generic geometric point $i_{\varepsilon}: 
{\rm Spec} F_q \lra {\rm Spec} {\cal O}$
 commutes 
with all the functors involved in the construction of the object 
${\cal F}^*_{\cal A, \cal B}$. So 
$i^*_{\varepsilon}{\cal F}^*_{\cal A, \cal B}$ is glued 
from the object $(i_{\varepsilon}^*\delta_{D_J})(-m)$ where $m \leq {\rm codim}D_J$. 
So by proposition \ref{10.24.01.2} there is  a non zero  element of 
\begin{equation} \label{11.29.01.2}
{\rm Hom}_{\overline F_q}
(i^*_{\varepsilon}{\cal F}^*_{\cal A, \cal B}\otimes_{F_q} \overline F_q, 
i^*_{\varepsilon}\delta_{{\cal D}_J}(-m)\otimes_{F_q}\overline F_q)
\end{equation}
invariant with respect to ${\rm Gal}(\overline F_q/F_q)$. 
Further, for generic $\varepsilon$ the vector spaces 
(\ref{11.29.01.1})  is isomorphic to the one 
 (\ref{11.29.01.2}), so that the isomorphism transforms one Galois 
action to the other 
via a surjective   projection
${\rm Gal}(\overline K/K) \lra {\rm Gal}(\overline F_q/F_q)$. So we identify 
the Tate twists in  (\ref{11.29.01.1})  and (\ref{11.29.01.2}). 
The l-adic version of the lemma \ref{11.1.01.4} is proved. 

The rest of the proof of the $l$-adic version of the theorem 
copies the end of the proof of the Hodge   version. 
 The theorem is proved.

{\it A specialization theorem data}. Suppose we have a family of divisors 
$D(\varepsilon) = A(\varepsilon) \cup B(\varepsilon)$ in $X$ parametrized by 
$\varepsilon \in \Sigma^0$ and satisfying all the conditions of theorem \ref{11.1.01.3}. 
It is given by a divisor ${\cal D} = {\cal A} \cup {\cal B} \subset X \times \Sigma^0$. 
Since 
$({\cal A}, {\cal B})$ is admissible, there are perverse sheaves 
$$
\widetilde {\cal F}_{\cal B}, \widetilde {\cal F}^!_{\cal A, \cal B}, 
\widetilde {\cal F}^*_{\cal A, \cal B}, \widetilde {\cal F}^*_{\cal A} 
$$
on $X \times \Sigma^0$ and morphisms
\begin{equation} \label{11.25.01.11}
\widetilde {\cal F}_{\cal B}  \stackrel{\widetilde \beta}{\lra }
\widetilde {\cal F}^!_{\cal A, \cal B}\stackrel{\widetilde c}{\lra }
\widetilde {\cal F}^*_{\cal A, \cal B}\stackrel{\widetilde \alpha}{\lra }
\widetilde {\cal F}^*_{\cal A} 
\end{equation}

Let us assume in addition that 
\begin{equation} \label{11.25.01.13}
{\rm Sp}_s\widetilde {\cal F}_{\cal A}   = {\cal F}_{A}; \qquad 
{\rm Sp}_s\widetilde {\cal F}_{\cal B}   = {\cal F}_{B}
\end{equation}

\begin{corollary} \label{11.25.01.12} Assuming (\ref{11.25.01.13}), 
 morphisms (\ref{11.25.01.11}) induce isomorphisms
\begin{equation} \label{11.25.01.14}
{\rm Gr}^{DW}_n({\rm Sp}_s\widetilde {\cal F}^!_{\cal A, \cal B}) 
\stackrel{\sim}{\lra}{\rm Gr}^{DW}_n({\rm Sp}_s\widetilde {\cal F}^*_{\cal A, \cal B})
\stackrel{\sim}{\lra}{\rm Gr}^{DW}_n({\cal F}_{A})
\end{equation}
\begin{equation} \label{11.25.01.15}
{\rm Gr}^{T}_0({\cal F}_{B}) 
\stackrel{\sim}{\lra}{\rm Gr}^{T}_0({\rm Sp}_s\widetilde {\cal F}^!_{\cal A, \cal B})
\stackrel{\sim}{\lra}{\rm Gr}^{T}_0({\rm Sp}_s\widetilde {\cal F}^*_{\cal A, \cal B})
\end{equation}
\end{corollary}

{\bf Proof}. Using proposition \ref{11.7.01.12} we get a similar isomorphism for the 
$\widetilde {\cal F}$-sheaves on $X \times \Sigma^0$. Since 
the specialization functor is exact we get them for 
${\rm Sp}_s\widetilde {\cal F}$-sheaves. It remains to use (\ref{11.25.01.13}).
The corollary is proved.

Suppose now that we are given sections of the following local systems on 
$\Sigma^0$:
\begin{equation} \label{11.25.01.1}
[\widetilde \omega_{\cal A}] 
\in \Gamma(\Sigma^0; {\rm Gr}^W_{n+{\rm dp}({\cal A})}(p_* \widetilde {\cal F}_{\cal A}))
\end{equation}
\begin{equation} \label{11.25.01.2}
[\widetilde \Delta_{\cal B}] \in \Gamma(\Sigma^0; {\rm Gr}^W_{n-{\rm dp}({\cal B})}
(p_* \widetilde {\cal F}_{\cal B}))
\end{equation}
In this situation there are two different framed objects:
\begin{equation} \label{11.25.01.4}
(H^0(X, {\cal F}^*_{A,B}), [\omega_A],  [\Delta_B])
\end{equation}
and 
\begin{equation} \label{11.25.01.5}
(H^0(X, {\rm Sp}_s\widetilde {\cal F}^*_{\cal A,\cal B}), {\rm Sp}_s
[\widetilde \omega_{\cal A}],  
{\rm Sp}_s[\widetilde \Delta_{\cal B}])
\end{equation}

We will assume that local systems staying on the right of 
(\ref{11.25.01.1})-(\ref{11.25.01.2}) extend to $0 \in \Sigma$, the sections 
$[\widetilde \omega_{\cal A}]$ and $[\widetilde \Delta_{\cal B}]$ 
extend to sections $[\overline \omega_{\cal A}]$ and 
$[\overline  \Delta_{\cal B}]$ over $\Sigma$. Then ${\rm Sp}_s
[\widetilde \omega_{\cal A}]$ and $  
{\rm Sp}_s[\widetilde \Delta_{\cal B}]$ do not depend on the choice of $v$ (and thus $s$) 
and can be identified with the values of the sections $[\widetilde \omega_{\cal A}]$ and $[\widetilde \Delta_{\cal B}]$ at zero. 

Therefore having (\ref{11.25.01.13})  it makes sense to ask that these values 
at $0$ coincide with the initially given  $[\omega_{A}]$ and 
$[\Delta_B]$:
\begin{equation} \label{11.25.01.3}
[\overline \omega_{\cal A}(0)] = [\omega_A];\quad [\overline  \Delta_{\cal B}(0)] = [\Delta_B]
\end{equation}
Below we will assume (\ref{11.25.01.3}).

Observe that (\ref{11.25.01.13}) implies that 
${\rm dp}({\cal A}) = {\rm dp}({A})$ and ${\rm dp}({\cal B}) = {\rm dp}({B})$.  Therefore 
using isomorphism (\ref{11.25.01.13})
and  equality (\ref{11.25.01.3}) we can identify the classes
$$
{\rm Sp}_s[\widetilde \omega_{\cal A}] \in 
{\rm Gr}^W_{n+{\rm dp}({\cal A})}H^0(X; {\rm Sp}_s\widetilde {\cal F}_{{\cal A}}) \quad \mbox{and} \quad [\omega_{A}] \in  
{\rm Gr}^W_{n+{\rm dp}({A})}H^0(X; {\cal F}_{{A}})
$$
and similarly 
$$
{\rm Sp}_s[\widetilde \Delta_{\cal B}] \in 
{\rm Gr}^W_{n-{\rm dp}({\cal B})}H^0(X; {\rm Sp}_s\widetilde {\cal F}_{B}) \quad \mbox{and} \quad [\omega_{B}] \in   
{\rm Gr}^W_{n-{\rm dp}({B})}H^0(X, {\cal F}_{B})
$$
Summarizing, we identified the framings of the mixed objects 
(\ref{11.25.01.4}) and (\ref{11.25.01.5}).

\begin{theorem} \label{10.24.01.1} 
a) The framed objects (\ref{11.25.01.4}) and (\ref{11.25.01.5}) are 
equivalent.

b) If the divisor $D(\varepsilon)$ provides  a Tate stratification for 
generic 
$\varepsilon \in \Sigma^0$, then these framed objects
 are mixed Tate objects, and they are equivalent as mixed Tate objects.
\end{theorem}

{\bf Proof}. Let us rewrite the definition of ${\cal F}^{\bullet}_{A,*,B}$ 
in a different way. Recall that 
$$
U^A_k := U(k-1) \cup \mbox{all codimension $k$ pure $A$-strata and mixed strata} 
$$
Consider the open embeddings
$$
q^A_k: U(k-1) \hra X; \qquad q^B_k: U_k^A \hra X
$$
Then, since $q_{1}^{A*} \delta_X = \delta_U$, one has 
$$
{\cal F}^{\bullet}_{A,*,B} = q_{n!}^Bq^{B!}_n q_{n*}^A q^{A*}_n ... 
  q_{1!}^Bq^{B!}_1 q_{1*}^A q_{1}^{A*} \delta_X
$$
Set
$$
Q_k^{A}:= q_{k*}^{A} q_{k}^{A*}; \qquad Q_k^B:= q_{k!}^B q_{k}^{B!}
$$
$$
Q_{k_1, ...,k_m}^{C_1, ..., C_m}:= q_{k_1}^{C_1} \circ ... \circ q_{k_m}^{C_m}; 
\quad C_i = \quad 
\mbox{$A$ or $B$}
$$
Let ${\cal G}^{\bullet} \in D_{\rm Sh}^b(X)$ be an object 
whose restriction to $U$ is isomorphic to $\delta_U$. 
Then   one obviously has 
an isomorphism 
\begin{equation} \label{11.2.01.12q}
{\cal F}^{\bullet}_{A,*,B} = Q_{n,n,...,1,1}^{B, A, ..., B, A}{\cal G}^{\bullet}
\end{equation}
Therefore using the adjunctions ${\rm Id}_X \lra Q_k^{A}$ and $Q_k^{B} \lra {\rm Id}_X$ 
 we get a sequence of 
morphisms:
$$
{\cal F}^{\bullet}_{A,*,B} = Q_{n,n,...,1,1}^{B, A, ..., B, A}{\cal G}^{\bullet} 
\longleftarrow
Q_{n,n,...,1}^{B, A, ..., B}{\cal G}^{\bullet} \lra  Q_{n,n, ...,2,2}^{B,A, ..., B,A }{\cal G}^{\bullet} \longleftarrow
$$
$$
... 
\lra 
Q_{n,n}^{B,A}{\cal G}^{\bullet}\longleftarrow Q_{n}^{B}{\cal G}^{\bullet}
\lra {\cal G}^{\bullet}
$$
Applying the functor 
$H_{\tau}^0(-)$ to this sequence we get a similar 
one
$$
{\cal F}^*_{A,B} = H_{\tau}^0\Bigl(Q_{n,n,...,1,1}^{B, A, ..., B, A}{\cal G}^{\bullet}\Bigr)
 \longleftarrow 
H^0_{\tau}\Bigl(Q_{n, n, ...,1}^{B, A, ...,  B}{\cal G}^{\bullet}\Bigr) \lra
$$
\begin{equation} \label{11.2.01.12}
H^0_{\tau}\Bigl(Q_{n,n,...,2,2}^{B, A, ..., B,A}{\cal G}^{\bullet}\Bigr)\longleftarrow ... 
\lra 
H^0_{\tau}\Bigl(Q_{n,n}^{B,A}{\cal G}^{\bullet}\Bigr)\longleftarrow 
H^0_{\tau}\Bigl(Q_{n}^{B}{\cal G}^{\bullet}\Bigr)\lra {\cal G}
\end{equation}
of perverse sheaves relating 
${\cal F}^*_{A,B}$ and ${\cal G}:= H_{\tau}^0{\cal G}^{\bullet}$.

A similar treatment for the perverse sheaves ${\cal F}^!_{A,B}$ looks as follows. 
There are  open embeddings
$$
p^A_k = q^A_k: U(k-1) \hra X; \qquad p^B_k: V_k^A \hra X
$$
where 
$$
V^A_k := U(k-1) \cup \mbox{all codimension $k$ pure $A$-strata}
$$
Then setting 
$P_k^{A}:= p_{k*}^{A} p_{k}^{A*}$ and $P_k^B:= p_{k!}^B p_{k}^{B!}$ we get 
$$
{\cal F}^{!}_{A,B} = H^0_{\tau}\Bigl(P_{n,n,...,1,1}^{B, A, ..., B, A}\delta_X\Bigr)
$$

Set
$$
{\cal G}_*^{\bullet}:= {\rm Sp}_s\widetilde {\cal F}^{\bullet}_{\cal A, *, \cal B}; \qquad
{\cal G}_!^{\bullet}:= {\rm Sp}_s\widetilde {\cal F}^{\bullet}_{\cal A, !, \cal B}
$$
Since specialization sends perverse sheaves to perverse sheaves we have 
\begin{equation} \label{12.1.01.3}
H^0_{\tau}({\cal G}_*^{\bullet}):= {\rm Sp}_s\widetilde {\cal F}^{*}_{\cal A, \cal B}; \qquad
H^0_{\tau}({\cal G}_!^{\bullet}):= {\rm Sp}_s\widetilde {\cal F}^{!}_{\cal A,  \cal B}
\end{equation} 
Observe that restrictions of these objects to $U$ are isomorphic to $\delta_U$, providing 
isomorphism (\ref{11.2.01.12q}).

\begin{proposition} \label{11.25.01.20} 
There are natural maps 
$$
{\cal F}_B \stackrel{}{\lra}  Q^{B, A, ..., C}_{n,n,  ..., m}
{\cal G}_!^{\bullet} \stackrel{}{\lra}  Q^{B,A,  ..., C}_{n,n,  ..., m}
{\cal G}_*^{\bullet} \stackrel{}{\lra} {\cal F}_{A} \quad 
\mbox{where $C=A$ or $B$}
$$
The first map induces isomorphism on ${\rm Gr}^{T}_{0}H^0_{\tau}(-)$, 
the second on both 
 ${\rm Gr}^{DW}_{n}H^0_{\tau}(-)$ and ${\rm Gr}^{T}_{0}H^0_{\tau}(-)$, and the last one 
induces isomorphism on 
${\rm Gr}^{T}_{0}H^0_{\tau}(-)$. 

These maps give rise to a commutative diagram
\begin{equation} \label{11.25.01.21}
\begin{array}{ccccccccc}
{\cal F}_{A}&=&{\cal F}_{A}&=&...&=&{\cal F}_{A}&=&{\cal F}_{A}\\
\alpha \uparrow &&\uparrow &&&&\uparrow && \uparrow\\
{\cal F}^*_{A, B}& \longleftarrow &H^0_{\tau}
\Bigl(Q^{B, A, ..., B}_{n,n,  ..., 1}{\cal G}_*^{\bullet}\Bigr)& 
\lra &...&\longleftarrow  &H^0_{\tau}
(Q^{B}_{n}{\cal G}_*^{\bullet}) &\lra & H^0_{\tau}({\cal G}_*^{\bullet}) \\
c \uparrow &&\uparrow &&&&\uparrow && \uparrow\\
{\cal F}^!_{A, B}& \longleftarrow &H^0_{\tau}
\Bigl(P^{B, A, ..., B}_{n,n,  ..., 1}{\cal G}_!^{\bullet}\Bigr)& 
\lra &...&\longleftarrow  &H^0_{\tau}
(P^{B}_{n}{\cal G}_!^{\bullet}) &\lra & H^0_{\tau}({\cal G}_!^{\bullet}) \\
\beta \uparrow &&\uparrow &&&&\uparrow && \uparrow\\\
{\cal F}_{B}&=&{\cal F}_{B}&=&...&=&{\cal F}_{B}&=&{\cal F}_{B}
\end{array}
\end{equation} 
\end{proposition}

Thanks to (\ref{11.25.01.13}) and (\ref{12.1.01.3}) 
the  right vertical maps are identified with 
$$
{\cal F}_{B} = {\rm Sp}_s \widetilde {\cal F}_{\cal B}\lra 
{\rm Sp}_s\widetilde {\cal F}^!_{\cal A, \cal B} \lra 
{\rm Sp}_s\widetilde {\cal F}^*_{\cal A, \cal B} \lra {\rm Sp}_s
\widetilde {\cal F}_{\cal A} = 
{\cal F}_{A} 
$$

{\bf Proof}. There is canonical morphism 
${\cal G}^{\bullet}_* \lra {\cal F}_A$ given by
\begin{equation} \label{12.24.01.1rt}
{\cal G}^{\bullet}_* = {\rm Sp}_s(\widetilde {\cal F}^*_{\cal A, \cal B}) 
\lra {\rm Sp}_s(\widetilde {\cal F}_{\cal A})  = {\cal F}_{A}
\end{equation} 
We define the 
maps 
$$
Q_{n, n, ..., m,m}^{B, A, ..., B,A}{\cal G}^{\bullet}_* = j_{n!}^B j_{n*}^A ... 
j_{m!}^B j_{m*}^A {\rm Res}_{U(m-1)}{\cal G}^{\bullet}_* \lra {\cal F}_{A}
$$
and 
$$
Q_{n, n, ..., m}^{B, A, ..., B}{\cal G}^{\bullet}_* = j_{n!}^B j_{n*}^A ... 
j_{m!}^B {\rm Res}_{U^A(m-1)}{\cal G}^{\bullet}_* \lra {\cal F}_{A}
$$
by induction using restrictions of the map (\ref{12.24.01.1rt}) to $U(m-1)$ and 
$U^A(m-1)$ and the two constructions i) and ii) 
introduced in the proof of lemma \ref{11.7.01.1}. This generalizes 
the procedure used in lemma \ref{11.7.01.1}. The starting point is the isomorphism
${\rm Res}_{U(0)}{\cal G}^{\bullet}_* = {\rm Res}_{U(0)}{\cal F}_{A}$. 
After this we apply the construction i) for $U(0)$, then ii) for $U^A(1)$, 
then i) for $U(1)$, after that ii) for $U^A(1)$, and so on. 

To show that the top rectangle in diagram (\ref{11.25.01.21}) is commutative 
we need to show that the following diagram is commutative
$$
\begin{array}{ccccc} \label{12.24.01.2}
{\cal F}_A & = & {\cal F}_A & = & {\cal F}_A \\
\uparrow &&\uparrow &&\uparrow \\
Q_{n, n, ..., m-1}^{B, A, ..., B}{\cal G}^{\bullet}_*&\lra &
Q_{n, n, ..., m,m}^{B, A, ..., B,A}{\cal G}^{\bullet}_*&\longleftarrow &
Q_{n, n, ..., m}^{B, A, ..., B}{\cal G}^{\bullet}_*
\end{array}
$$
To show that the left square is commutative it is sufficient to show that the diagram
$$
\begin{array}{ccc}
{\rm Res}_{U(m-1)}{\cal F}_A & = 
&{\rm Res}_{U(m-1)}{\cal F}_A  \\
\uparrow (\ref{12.24.01.12}) &&\uparrow   \\
(j_{m-1}^B)_!{\rm Res}_{U^A(m-2)}{\cal G}^{\bullet}_*  & \lra 
&{\rm Res}_{U(m-1)}{\cal G}^{\bullet}_* 
\end{array}
$$
is commutative. Since the left bottom object in this diagram is isomorphic to 
$$
(j^B_{m-1})_!j_{m-1}^{B!}{\rm Res}_{U(m-1)}{\cal G}^{\bullet}_*
$$
 this boils down to 
the fact that the adjunction $j_!j^! \lra {\rm Id}$ is a morphism of functors.

To show that the right square 
is commutative it is sufficient to show that the diagram
$$
\begin{array}{ccc}
{\rm Res}_{U^A(m-1)}{\cal F}_A & = 
&{\rm Res}_{U^A(m-1)}{\cal F}_A  \\
\uparrow (\ref{12.24.01.11}) &&\uparrow   \\
j_{m*}^B{\rm Res}_{U(m-1)}{\cal G}^{\bullet}_*  & \longleftarrow 
&{\rm Res}_{U^A(m-1)}{\cal G}^{\bullet}_* 
\end{array}
$$
is commutative. Since the left bottom object in this diagram is isomorphic to 
$$
j^A_{m*}j_{m}^{A*}{\rm Res}_{U^A(m-1)}{\cal G}^{\bullet}_*, 
$$
 This boils down to 
the fact that the adjunction ${\rm Id}\lra j_*j^* $ is a morphism of functors.

Therefore we defined the top rectangle of the diagram (\ref{11.25.01.21}) 
and proved that all its squares 
 are commutative. Applying the duality $\ast$ and interchanging the role of $A$ and $B$ we deduce 
from this a similar  statement about the bottom rectangle. 
The maps involved in the middle rectangle are constructed using the canonical map 
$j_! \lra j_*$ for the extensions to the mixed strata. 
It follows that the middle rectangle is also built from commutative squares. 
The commutative diagram (\ref{11.25.01.21}) is constructed.

It remains to prove that the vertical arrows in this diagram induce the 
isomorphisms on the appropriate pieces of the $WD$ or/and $T$ filtrations. 
This follows from lemma \ref{12.25.01.1}. 
The condition (\ref{12.26.01.1}) for $i=0$ is deduced from theorem 
\ref{11.1.01.3}.  For other $i$'s we use a similar result, see the remark 2 
after theorem \ref{11.1.01.3}.

Applying the duality $\ast$ and interchanging the role of $A$ and $B$ we deduce 
the statements about the bottom rectangle from the ones about the top rectangle. 

The maps involved in the middle rectangle are constructed using the canonical map 
$j_! \lra j_*$ for the extensions to the mixed strata. The proposition is proved.

The specialization theorem follows immediately from this proposition. 
Indeed, applying the functor $H^0(X,-)$ to the second horizontal sequence of morphisms in 
(\ref{11.25.01.21}) we get a sequence of mixed objects and maps between them connecting 
$H^0(X; {\cal F}^*_{A, B})$ and $H^0(X; {\rm Sp}_s\widetilde {\cal F}^*_{\cal A, \cal B})$. 
One needs to check that each of these objects is naturally 
framed, and each of the maps 
between them 
respects the frames. The $A$-parts of the frames on $H^0(X; -)$, where $-$ stays 
for the objects 
in the second horizontal line,  are provided by the vertical 
frame maps to ${\cal F}_A$. Since the top 
rectangle of the diagram is commutative and the top vertical arrows induce isomorphisms on 
${\rm Gr}^{DW}_{n}$,  the maps respect the $A$-parts of the frame by their very construction, see s. 2.4-2.5. 

To handle the $B$-parts of the frame we employ a similar argumentation 
for the bottom rectangle in the diagram, 
and in addition use the fact that the middle vertical maps induce isomorphisms on 
${\rm Gr}^{T}_{0}H^0_{\tau}(-)$. The part a) of the theorem is proved. 
Having a), the part b) is straitforward. The theorem is proved.

 \section{Applications to scissor congruence groups}
 
In this section we apply the results of 
the sections 2-4 when $D= A \cup B$ is an admissible configuration of hyperplanes in $P^n$. 
In particular when $A$ and $B$ are two simplices in $P^n$ 
we produce  a framed mixed Hodge structure of geometric origin whose 
period is given by the Aomoto polylogarithm. 
The $l$-adic version of this result provides  
 a framed mixed Tate l-adic representation 
of the Galois group ${\rm Gal}(\overline F/F)$. When $F$ is a number field we produce a mixed Tate 
motive over $F$ whose Hodge and $l$-adic realizations are 
the ones mentioned above.

Similar results are valid for the scissor congruence groups 
considered in [G9], including the classical scissor congruence groups 
in the hyperbolic and spherical spaces.

{\bf 1. The generalized scissor congruence groups and 
framed Hodge-Tate structures}. 
An algebraic  simplex in $P^n$ is a union of $n+1$ 
hyperplanes. It is called a nondegenerate algebraic  simplex if the intersection of 
these hyperplanes is empty. 
Let  
$$
A = A_0 \cup ...\cup A_n\quad \mbox{and} \quad B= B_0 \cup ...\cup B_n
$$ 
be two nondegenerate algebraic simplices  in $P^n$. Here $A_i$ and $B_j$ are hyperplanes. 
The pair $(A,B)$ is admissible 
if and only if $A_I \not = B_J$ if $|I| = |J|$. 
The stratification defined by any collection of hyperplanes in $P^n$ is always regular. 

Let us introduce a $\Z/2\Z$-torsor of orientations of a simplex $A$. Namely, 
an ordering of the hyperplanes $A_i$ provides an element of the torsor, and 
two orderings provide the same element if and only if they differ by an even permutation.

For a field $F$ the  generalized 
scissor congruence group $A_n(F)$ (see [BMS], [BGSV]) is generated by admissible pairs $$
(A; B) = 
(A_0, ..., A_n; B_0, ..., B_n)
$$
 of oriented simplices 
 in ${P}^n(F)$ subject to the following relations:

\vskip 3mm \noindent

1){\it Nondegeneracy}.  $(A; B) = 0$ if one of the simplices $A$ or $B$ is degenerate.

2) {\it Orientation}. $(A;B)$ changes the sign if we change orientation of $A$ or $B$.

3){\it Additivity}.  For any $n+2$ hyperplanes $A_0,...,A_{n+1}$ one has 
$$
\sum_{i=0}^{n+1}(-1)^i (A_0 ...,\hat A_i,...,A_{n+1}; B_0,..., B_n) =0
$$
if all  the terms are admissible (additivity in $A$). 

We 
impose a similar  additivity in $B$ condition.

4) {\it Projective invariance}. $(gA; gB) = (A; B)$ for any $g \in PGL_{n+1}(F)$.

Recall the group ${\cal H}_n$ of the Hodge-Tate structures framed by 
$\Q(0)$ and $\Q(-n)$. 
Its motivic version is the corresponding group   ${\cal A}_n(F)$ of framed mixed Tate 
motives over a number field $F$. 
The l-adic counterpart 
 ${\cal A}^{ et}_n(F)$ is the $\Q_l$-vector space of equivalence 
classes of mixed Tate l-adic ${\rm Gal}(\overline F/F)$-modules framed 
by $\Q_l(0)$ and $\Q_l(-n)$ ([BD1], see also [G7]). Here we have to assume that $F$ does not contain all 
$l^{\infty}$ roots of unity.

Let ${\cal F}^*_{A,B}$ be the perverse mixed Hodge sheaf defined 
in section 2, and by ${\cal F}^{et, *}_{A,B}$ its l-adic counterpart. 

\begin{theorem} \label{11.18.01.5} a) There is a canonical homomorphism 
$h_n: A_n(\C) \lra {\cal H}_n$ defined on the generators by 
\begin{equation} \label{11.18.01.6}
h^{\cal H}_n: (A,B) \lms H^0\Bigl(\C P^n; {\cal F}^*_{A,B} \Bigr)
\end{equation} 

b) Let $F$ be a field that does not contain all 
$l^{\infty}$ roots of unity. Then there is a canonical homomorphism 
$h^{et}_n: A_n(F) \lra {\cal A}^{ et}_n(F)$ defined  by 
\begin{equation} \label{11.18.01.6}
h^{et}_n: (A,B) \lms H^0_{et}\Bigl(P^n \otimes \overline F; {\cal F}^{et, *}_{A,B} \Bigr)
\end{equation} 

c) Let $F$ be a number field. 
 Then there is a canonical homomorphism 
$h_n: A_n(F) \lra {\cal A}_n(F)$ whose Hodge and l-adic realizations provided by a) and b). 
\end{theorem} 

{\bf Remark}. If we restrict our attention to the subgroup 
generated by pairs $(A;B)$ of simplices in generic position, i.e. $A \cup B$ 
is a normal crossing divisor, then a) and b) follow from  [BGSV], and c) from chapter 5 in [G9].

\begin{conjecture} The map constructed in theorem \ref{11.18.01.5}c) is an isomorphism. 
\end{conjecture}

{\bf Proof}. a) The right hand side in (\ref{11.18.01.6}) is obviously a Hodge-Tate structure. 
Indeed, the perverse sheaf ${\cal F}^*_{A,B}$ is glued from $\delta_C(-m)$ 
where $C$ are planes in $\C P^n$. 

Let $f_{i}$ be a rational function 
on $\C P^n$ with the  divisor is $A_i - A_0$. Set
$$
\omega_A = d\log f_1 \wedge ... \wedge  d\log f_n \in \Omega^n_{\rm log}(P^n - A)
$$
Let $\Delta_B$ be  a chain defining a class in 
$H_n(\C P^n, B(\C); \Z)$ corresponding to given orientation of the algebraic simplex $B$.
According to section 2.4  the 
classes
$$
[\omega_A] \in H^n_{\rm DR}(\C P^n - A) = {\rm Gr}^W_{2n}H^n_{\rm DR}(\C P^n - A)= \Z(-n)
$$
$$
[\Delta_B] \in H_n(\C P^n,
 B(\C); \Z) = {\rm Gr}^W_{0}H_n(\C P^n,
 B(\C); \Z) = \Z(0)
$$
 provide a framing of $h_n^{\cal H}(A,B)$.

It remains to check the relations. The relations 2) and 4) are obvious, and 1) 
is true by  definition. Let us check the additivity in $A$. 
We need the following simple general result.

\begin{lemma} \label{4.17.01.77} Let $M$ an object of a mixed 
Tate category ${\cal C}$ 
(e.g. the abelian 
category of mixed Tate motives over a number field $F$, or the category of 
Hodge-Tate structures). Then 

a) Given  non zero maps 
$$
v_0: {\rm Gr}^W_0M \lra \Q(0); \quad f_n^{\alpha}, 
f_n^{\beta}:  \Q(-n) \lra {\rm Gr}^W_{2n}M  
$$
such that  $\sum_{\alpha} f_n^{\alpha} = \sum_{\beta} f_n^{\beta} \not = 0$ one has an equality of  framed  objects
$$
\sum_{\alpha} (M, v_0, f_n^{\alpha}) = \sum_{\beta} (M, v_0, f_n^{\beta}) 
$$

b) Similarly if $\sum_{\alpha}v_0^{\alpha} = 
\sum_{\beta}v_0^{\beta} \not = 0$, then 
$$
\sum_{\alpha} (M, v_0^{\alpha}, f_n) = \sum_{\beta} (M, v_0^{\beta}, f_n) 
$$
\end{lemma}

{\bf Proof}. a) It is sufficient to prove that if 
$f_1, f_2: \Q(-n) \lra {\rm Gr}^W_{2n}M  $ 
are non zero maps whose sum is also non zero then 
$$
(M, v_0, f_1) + (M, v_0, f_2)  = (M, v_0, f_1+ f_2) 
$$
By definition the left element is represented by the framed object
$(M \oplus M, v_0 + v_0, (f_1,f_2))$. We claim that 
the natural projection $({\rm id},{\rm id}): M \oplus M \lra M$ 
induces an equivalence of the framed objects 
$$
(M \oplus M, v_0 + v_0, (f_1,f_2)) \stackrel{\sim }{\lra} (M, v_0, f_1+ f_2) 
$$
Indeed, there are commutative diagrams
$$
\begin{array}{ccc}
\Q(-n)&  =  & \Q(-n)\\
&&\\
(f_1, f_2) \downarrow &&\downarrow f_1 + f_2\\
&&\\
{\rm Gr}^W_{2n}M \oplus 
{\rm Gr}^W_{2n}M &\stackrel{({\rm id},{\rm id}) }{\lra} & {\rm Gr}^W_{2n}M
\end{array}
$$
and 
$$
\begin{array}{ccc}
{\rm Gr}^W_{0}M \oplus {\rm Gr}^W_{0}M & \stackrel{({\rm id},{\rm id}) }{\lra} & {\rm Gr}^W_{0}M\\
&&\\
v_0+v_0 \downarrow &&\downarrow v_0\\
&&\\
\Q(0)&=   & \Q(0)
\end{array}
$$
The proof of the second statement is similar. 
The lemma is proved. 

More generally, suppose that $A$ and $B$ are unions of hyperplanes in $P^n$, and the 
pair $A \cup B$ is admissible. 
Then there is a Hodge-Tate structure $h_n^{\cal H}(A,B)$ defined by the same formula 
(\ref{11.18.01.6}). 

Let $A':= A_0 \cup ...  \cup A_{n+1}$. Set $A^{(i)}:= 
A_0\cup ...  \widehat A_i ...  \cup A_{n+1}$.  
Then obviously
\begin{equation} \label{1.2.02.2}
\sum^{*}_{0 \leq i \leq n+1}(-1)^i\omega_{A^{(i)}} = 0
\end{equation}
Here the summation is over $0 \leq i \leq n+1$ such that 
$[\omega_{A^{(i)}}] \not = 0$.

\begin{lemma} \label{1.2.02.1}
Let $A_0, ...  , A_{n}$ be any collection of hyperplanes in $P^n$. 
Let $f_i$ be a rational function with the divisor $A_i - A_{0}$. Then  
$$
\{f_1, ..., f_n\} \in K_n^M(F(P^n))/ F^* \cdot K_{n-1}^M(F(P^n))
$$ 
is not zero if and only if 
the hyperplanes are in generic position. 
\end{lemma}

{\bf Proof}. If the hyperplanes are in generic position then 
they are projectively equivalent to the coordinate hyperplanes $z_i =0$, so our symbol is 
$\{z_1/z_0, ..., z_n/z_0\}$. It is non zero 
since the homomorphism
$$
d\log^n: K_n^M(F(P^n))/ F^* \cdot K_{n-1}^M(F(P^n)) \lra \Omega^n_{\rm log}((F(P^n));
$$
$$
\{f_1, ..., f_n\} \lms d\log (f_1) \wedge ... \wedge d\log (f_n)
$$
maps it to a non zero differential form. 

If $g_1, ..., g_m$ are elements of a field $k$ 
such that $g_1 + ... + g_m =1$ then it is easy to check 
by induction that 
$\{g_1, ..., g_m\}=0$ in $K^M_n(k)$. Let $F_i$ be a linear  equation of the hyperplane 
$A_i$. If these hyperplanes are not in generic position then, renumbering the 
hyperplanes, we can find $F_1, ..., F_m$ such that $\lambda_1 F_1 + ... + \lambda_m F_m =0$, 
$\lambda_i \not = 0$.
Since $\{F_1, ..., F_m\} \sim \{\lambda_1F_1, ..., \lambda_mF_m\}$ modulo 
$F^* \cdot K_{m-1}^M(F(P^n))$ the lemma follows. 
 
 This lemma implies that the terms omitted during the summation in (\ref{1.2.02.2}) 
are precisely the ones which are zero by the non degeneracy relation. 
 Therefore using lemma \ref{4.17.01.77} we see that
$$
\sum^{*}_{0 \leq i \leq n+1} (-1)^i(h_n^{\cal H}(A',B), [\omega_{A^{(i)}}], [\Delta_B]) = 0
$$
The canonical morphism of perverse sheaves 
${\cal F}^*_{A^{(i)},B} \lra {\cal F}^*_{A',B}$ provides, assuming $[\omega_{A^{(i)}}] \not = 0$,  
an equivalence of framed Hodge-Tate structures 
$$
(h_n^{\cal H}(A',B), [\omega_{A^{(i)}}], [\Delta_B]) \sim 
(h_n^{\cal H}(A^{(i)},B), [\omega_{A^{(i)}}], [\Delta_B])
$$
So the additivity in $A$ is proved. 
The additivity in $B$ is checked similarly. Or we can use the duality since
$\ast h_n^{\cal H}(A,B) = h_n^{\cal H}(B, A)$. The part a) of the theorem is proved. 

The part b) is completely similar to a). The only point requiring a comment 
is construction of the $A$-part of the framing. Namely, the class $[\omega_{A}]$ in the l-adic setting is the cohomology class given by the l-adic regulator 
applied to the symbol $\{f_1, ..., f_n\}$ where $f_i$ are as in lemma 
\ref{1.2.02.1}. Then for $n+2$ hyperplanes we have 
$\sum_i(-1)^i \{f_1, ... \widehat f_i, ... , f_{n+1}\} =0$ providing, together with 
lemma \ref{1.2.02.1}, identity (\ref{1.2.02.2})

The part c) follows from the results of chapter 3 and, say, a) 
using the injectivity of the regulators, or can be deduced directly 
using the same arguments as used in the proof of part a). 
The theorem is proved. 

{\bf Remark}. A similar result for the scissor congruence groups defined 
in [G9] is left to the reader as an easy exercise.

\section{Applications to  motivic torsors of path on curves}

In this section we apply the results of chapters 2-4 
when $A \cup B$ is a specific configuration of divisors 
on the $n$-th power of an arbitrary  regular curve $X$. 

Recall that 
in the Hodge or \'etale realization the torsor of path 
${\cal P}(X; v_x, v_y)$ between the 
tangential base points $v_x, v_y$ was constructed 
by Deligne [D] as a pro-object in the corresponding category. 
There is a weight filtration indexed by integers $n \leq 0$ on it, and 
${\rm Gr}^W_0{\cal P}(X; v_x, v_y) = \Q(0)$ (or $\Q_l(0)$). 
Let $v_0$ be the image of $1 \in \Q(0)$ under this isomorphism. 
Choose a non zero vector $f_k \in ({\rm Gr}^W_{-k}{\cal P}(X; v_x, v_y))^{\vee}$.
We  construct  
a framed object of geometric origin equivalent to the framed object
\begin{equation} \label{1qaq}
({\cal P}(X; v_x, v_y), v_0, f_k) 
\end{equation}
The period of its  Hodge realization is given by an iterated integral. This construction provides a  framed 
mixed Hodge structure of geometric origin corresponding 
to any iterated 
integral between tangential base points on $X$.  
In the case of the classical  base points this has been done by Beilinson, and our 
construction in this case leads to the same object. 

Let $F$ be a number field and $v_x, v_y$
   are non zero tangent vectors at the points $x, y \in {\Bbb P}^1(F)$. 
Suppose that  $\Lambda$ is a finite subset of $ {\Bbb P}^1(F)$,  and
 $v_x, v_y$ are  defined over $F$. 
It was proved in [DG] that there exists  the motivic torsor 
\begin{equation} \label{1}
{\cal P}^{\cal
  M}({\Bbb P}^1 -\Lambda; v_x, v_y)
\end{equation}
 of path on ${\Bbb P}^1 -\Lambda$ between the tangential
base points $v_x$ and $v_y$ understood as a pro-object in the abelian category 
${\cal M}_T(F)$ of mixed Tate motives over  $F$. Its Hodge and
$l$-adic realizations are isomorphic to the standard  Hodge and
$l$-adic 
realizations of the torsor of path.

Applying our general construction in the case $X = {\Bbb P}^1 -\Lambda$ we get
 an independent   geometric construction of  the motivic torsor 
of path (\ref{1}). Namely, we describe it by constructing all its matrix elements. 
We use the injectivity of the regulators to prove 
that it has all the needed properties.

{\bf 1. A framed mixed Hodge structures corresponding 
to iterated integrals on curves}.  Suppose that  $X$ is 
a regular projective curve over $\C$ and 
$\omega_i \in \Omega^1_{\log}(X)$. 
Let $p_i: X^n \lra X$ be the projection onto the $i$-th factor. 
 Choose a
path $\gamma: [0,1] \lra X(\C)$ 
connecting the (tangential) base points  $v_x,
v_y$  at $x$ and $y$. Then there is 
an iterated integral 
\begin{equation} \label{1.21.01.12}
\int_{\gamma; v_x, v_y} \omega_1 \circ ... \circ \omega_n:= 
\int_{\gamma(\Delta_n)} p_1^*\omega_1 \wedge ... \wedge p_n^*\omega_n
\end{equation}
Here
$
\Delta_n = \{0 \leq s_1 \leq
... \leq s_n \leq 1\}
$ is the standard $n$-dimensional simplex and 
\begin{equation} \label{1.28.01.10}
\gamma: \Delta_n \lra X^n(\C), \qquad \gamma(s_1, ..., s_n):=
(\gamma(s_1), ..., \gamma(s_n))
\end{equation}
Denote by ${\rm Sing}(\omega)$ the singular locus of a form $\omega$. 
Integral (\ref{1.21.01.12}) 
is convergent if and only if 
\begin{equation} \label{1.21.01.13}
x \not \in {\rm Sing}(\omega_1) \quad \mbox{and} \quad y \not \in {\rm
  Sing}(\omega_n) 
\end{equation}
Our goal is to provide an explicit construction 
of a framed mixed Hodge structure whose period is given by 
a convergent iterated integral (\ref{1.21.01.12}). 

Let $X$ be
a regular projective curve over a field $F$ and 
$\omega_i \in \Omega^1_{\log}(X)$. 
Set
$$
A_i := p_i^{-1}{\rm Sing}(\omega_i); \qquad   A := A_1 \cup ... \cup A_n
$$
So 
$$
X^n - A = (X- {\rm Sing}(\omega_1))\times ... \times (X- {\rm Sing}(\omega_n))
$$
Observe that ${\rm dp}(A)$ is the number of singular forms among the forms 
$\omega_i$. 
Let $(t_1, ..., t_n)$ be a point of $X^n$. 
Set 
$$
B:= \{x = t_1\}\cup  \cup_{i=1}^{n-1}\{t_i = t_{i+1}\}\cup  \{t_n = y\};
$$

\begin{lemma} \label{11.19.01.10} a) $A \cup B$ defines a regular stratification of $X^n$. 

b) The pair of divisors $(A,B)$ is admissible 
if and only if condition (\ref{1.21.01.13}) is satisfied. 
\end{lemma}

{\bf Proof}. a) If 
 $\{a^k_i\}$ are given points of $X$ then 
the stratification defined by the divisors $t_i = a^k_i$ and any collection of 
 diagonals 
$t_i = t_j$ in $X^n$ 
is obviously regular. 

b) Straitforward check. 
The lemma is proved. 

Below we assume (\ref{1.21.01.13}). 
So thanks to lemma \ref{11.19.01.10} the  construction of section 2 applied to the 
divisor $A \cup B$ in  $X^n$ provides 
a mixed Hodge structure
$H^0(X^n; {\cal F}_{A,B}^*)$.  
It is equipped with a framing  
provided by  s. 2.5 by   the  element 
$$
[p_1^*\omega_1 \wedge ... \wedge p_n^*\omega_n] \in {\rm Gr}^W_{n+{\rm dp}(A)}H^n(X^n -A),
$$
and the canonical element 
$$
[\Delta_n] \in {\rm Gr}^W_0H_n(X^n, B) = \Z(0),
$$
The sign of the generator here is determined by
 the natural ordering of irreducible components of divisor $B$. 
(If $x = y$ we modify slightly the definition of the group on
the right in a similar way as in  ch. 4 of [G7]). 
So we get a framed mixed Hodge structure 
\begin{equation} \label{1.22.01.1}
h(x; \omega_1\otimes ... \otimes \omega_n; y) := 
\Bigl(H^0(X^n; {\cal F}_{A,B}^*), [p_1^*\omega_1
\wedge ... \wedge p_n^*\omega_n],  [\Delta_n]\Bigr)
\end{equation}

If a more restrictive condition 
\begin{equation} \label{1.21.01.15}
x, y \not \in {\rm Sing}(\omega_i) \quad \mbox{for $1 \leq  i \leq n$}
\end{equation}
is valid then there is {\it another} framed mixed Hodge structure which has a  clear 
motivic origin, and whose 
  period is given 
by integral (\ref{1.21.01.12}). Its  construction 
 goes back to
 Beilinson. Namely set $B_A:= B - (B \cap A)$. Consider the mixed motive 
\begin{equation} \label{1.22.01.1q}
H^n(X^n -A, B_A) = H^0(X^n, \alpha_{*} j_{!}\delta_{
  U})
\end{equation}
where 
$$
j: U:= X^n - (A  \cup B) \hookrightarrow X^n - A;\quad   
\alpha: X^n - A  \hookrightarrow X^n; 
$$
Both inclusions here are open affine embeddings,  
so  $\alpha_{*} j_{!}\delta_{U}$ is a perverse sheaf. 

\begin{lemma} \label{1.21.01.151} Let us assume (\ref{1.21.01.15}). 
Then none of the pure $B$-strata is contained in a pure $A$-stratum. Therefore we have  
\begin{equation} \label{1.22.01.1p}
\alpha_{*} j_{!}\delta_{U} = {\cal F}^*_{A,B}
\end{equation}
\end{lemma}

{\bf Proof}. It is  sufficient to check that none of the $B$-vertices is contained in a pure $A$-stratum. The $B$-vertices are given by $t_1 = ... = t_k = x; t_{k+1} = ... = t_n =y$, 
so it follows immediately from (\ref{1.21.01.15}). Therefore $A\cup B$ 
is a union of mixed strata and pure $A$-strata. This implies that the inductive construction of the object ${\cal F}^{\bullet}_{A, *, B}$ can be replaced by a single 
 step construction given by $\alpha_{*} j_{!}\delta_{U}$, i.e. 
we have an isomorphism $\alpha_{*} j_{!}\delta_{U} = {\cal F}^{\bullet}_{A, *, B}$.
 In particular the latter object is a perverse sheaf. So we get (\ref{1.22.01.1p}). 
The lemma is proved. 

Therefore we have a natural framing on the mixed object (\ref{1.22.01.1q})

 \begin{lemma} \label{2.10.02.1} Assuming \ref{1.21.01.15},  
integral (\ref{1.21.01.12}) 
is a  period of the framed mixed Hodge structure (\ref{1.22.01.1q}). 
\end{lemma}

{\bf Proof}. The natural framing constructed in chapter 2 
admits in this case a more explicit description. Namely, the 
restriction of the $n$-form 
$p_1^*\omega_1 \wedge ... \wedge p_n^*\omega_n$  to $B$ is zero, and it
provides a well defined class
$$
[p_1^*\omega_1 \wedge ... \wedge p_n^*\omega_n] \in 
{\rm Gr}^W_{n+{\rm dp}(A)}H^n(X^n -A, B_A)
$$
Further, there is canonical isomorphism 
\begin{equation} \label{22113}
\Z(0)\stackrel{\sim}{=} {\rm Gr}^W_{0}H_n(X^n , B)    
= {\rm Gr}^W_{0}H_n(X^n -A, B_A) 
\end{equation}
providing the second component of the frame. 
(If $x = y$ we modify  the group on
the right as in ch. 4 of [G7]). 
We get  a framed mixed Hodge structure.  
\begin{equation} \label{1.22.01.1we}
h_B(x; \omega_1\otimes ... \otimes \omega_n; y) := \Bigl(H^n(X^n -A, B_A), [p_1^*\omega_1
\wedge ... \wedge p_n^*\omega_n],  [\Delta_n]\Bigr)
\end{equation}
A path $\gamma: [0,1] \lra X(\C)$ from $x$ to $y$ 
provides a lift of the generator of (\ref{22113}) to a relative homology class 
$$
[\gamma(\Delta_n)] \in H_n^{{\rm Betti}}(X^n -A, B_A; \Z) 
$$
It depends only on the homotopy class of $\gamma$. The lemma follows.

{\it Applications of the specialization theorem}. 
Let $\Sigma$ be a curve and $0 \in \Sigma$ is its distinguished  point.
 Let $\{\omega_i(\varepsilon)\}$ be  
$1$-forms on $X$ parametrized by $\varepsilon \in \Sigma$, and $\omega_i(0) = \omega_i $. 
Suppose that the base points $x(\varepsilon)$ and $y(\varepsilon)$ 
also depend regularly  on 
$\varepsilon \in \Sigma$, i.e. we are given two maps $x,y: \Sigma \lra X$ 
regular near $0$, and $x(0) =x, y(0)=y$. 
Suppose in addition to  condition (\ref{1.21.01.13}) that 
\begin{equation} \label{4.6.01.1.w}
x(\varepsilon), y(\varepsilon) \not \in {\rm Sing}(\omega_i(\varepsilon)) 
\quad \mbox{for} \quad \varepsilon \not = 0, \quad i = 1, ..., n
\end{equation}
Then for $\varepsilon \not = 0$ there are framed objects (\ref{1.22.01.1we}), while 
at $\varepsilon  = 0$ we constructed a framed object (\ref{1.22.01.1}). 
The specialization theorem \ref{10.24.01.1} in this particular case allows us to compare them:

\begin{theorem} \label{12.06.01.10} Let us assume conditions 
(\ref{1.21.01.13}) and (\ref{4.6.01.1.w}). Then 
one has 
$$
{\rm Sp}_{\varepsilon \to 0}h_B(x(\varepsilon); 
\omega_1(\varepsilon) \otimes ... \otimes \omega_n(\varepsilon); y(\varepsilon))   = 
h(x; \omega_1 \otimes ... \otimes \omega_n; y)
$$
\end{theorem}

{\bf Proof}. It follows immediately from theorem \ref{10.24.01.1} thanks to 
(\ref{1.22.01.1p}), (\ref{1.22.01.1q}) and the very definitions. 
The theorem is proved.

Here is a version of this result. Let $S:= \cup {\rm Sing}(\omega_i)$ and $X_S := X-S$.  
The family of the mixed Hodge structures $H^n(X_S^n, B_S)$
forms a unipotent 
variation ${\cal H}^n$ of mixed Hodge  structures over $X_S \times X_S $ 
(see chapter 4 in [G7]). 
Combining it with the rigidity theorem of Vologodsky [Vol] one
immediately sees that this variation is isomorphic to the one constructed by Hain
and Zucker [HZ]. 
The canonical class $[\Delta_n]$ and  the given forms $\omega_i$
provide a framing on it. We denote by 
${\cal H}^n(\omega_1\otimes  ...\otimes
  \omega_n)$ the corresponding variation of framed mixed Hodge 
  structures. Its fiber over a point $(x,y)\in X_S^2$ is precisely 
$h_B(x; \omega_1 \otimes ... \otimes \omega_n;y)$. 

\begin{corollary} \label{1.21.01.20} Let us assume 
  (\ref{1.21.01.13}). Then  specialization of the variation
  of framed 
Hodge-Tate structure ${\cal H}^n(\omega_1\otimes  ...\otimes
  \omega_n)$ 
at the tangent vector 
\begin{equation} \label{1.22.01.3}
(v_x, v_y) \in T_{(x,y)} (X \times X); \qquad v_x \not = 0, v_y \not = 0
\end{equation}
is equivalent to the framed 
Hodge-Tate structure $h(x; \omega_1\otimes  ...\otimes \omega_n;
  y)$. 
In particular it does not depend on the choice of vectors 
  $v_x, v_y$. 
\end{corollary}

{\bf Proof}. Apply twice theorem \ref{12.06.01.10}, using 
first specialization with respect to  $x \in X$, and then with respect to $y \in Y$. 
The corollary is proved.

If $x  \in {\rm Sing}(\omega_1)$ or $y \in {\rm Sing}(\omega_n)$ then 
the integral is divergent, and it has to be regularized. The
corresponding framed mixed Hodge  structure is {\it defined} 
by the specialization
of the variation ${\cal H}^n(\omega_1\otimes  ...\otimes
  \omega_n; \gamma)$ at the tangent vector (\ref{1.22.01.3}), and its
  period is, by definition,  the regularized value of the integral.

\begin{corollary} \label{11.19.01.10q} a) Let us assume (\ref{1.21.01.13}). Then 
iterated integral (\ref{1.21.01.12}) is the period of the 
 framed mixed Hodge structure (\ref{1.22.01.1}). 

b) If  $X= {\Bbb P}^1$ then (\ref{1.22.01.1}) it is a framed Hodge-Tate structure. 
\end{corollary}

{\bf Proof}. a) Follows, for instance from  lemma \ref{2.10.02.1} and 
theorem \ref{1.21.01.20}. Another way to deduce this result is provided by the results of 
Sections 3.2 and 3.4. The part 
b) is obvious. The corollary is proved.

{\bf 2. Framed mixed Tate motives corresponding to iterated integrals on ${\Bbb A}^1$}. 
Below $X = P^1$.

\begin{theorem} \label{2.8.02.1} Assume that
  $x$, $y$ and non zero forms $\omega_i$ are defined over a number field $F$ and 
(\ref{1.21.01.13}) holds. Then  
there is a framed mixed Tate motive 
$$
m(x; \omega_1\otimes ... \otimes \omega_n; y)\in {\cal A}_n(F)
$$
such that for a complex embedding $\sigma: F \hra \C$ its Hodge realization is provided by $
h(\sigma(x); \sigma(\omega_1) \otimes ... \otimes \sigma(\omega_n); \sigma(y))
$. 
\end{theorem} 

{\bf Proof}. We apply
 the construction of chapter 2 to the divisor $A\cup B \subset X^n$ 
defined in s. 6.1.  
Thanks to the results of chapter 3 we get a framed mixed Tate motive over $F$. 
Corollary \ref{1.21.01.20} provides the comparison of its Hodge 
realization with the construction used in [G7]. In particular the constructed 
element (\ref{1=wqa}) 
does not depend on the choice of tangent vectors $v_x, v_y$. 
The theorem is proved. 

Suppose that  $F$ is a number field, $v_x, v_y$
   are non zero tangent vectors at the points $x, y \in {\Bbb P}^1(F)$, 
  and
 $v_x, v_y$ are  defined over $F$. 

\begin{theorem} \label{2.8.02.11} In the   above situation  
for any  $a_1, ..., a_n \in F$ there is an element
\begin{equation} \label{1=wqa}
{\rm I}^{\cal M}(v_{x}; a_{1}, ..., a_{n}; v_{y}) \in 
{\cal A}_{n}(F)
\end{equation}
such that for a complex embedding $\sigma: F \hra \C$ 
its Hodge realization coincides with the defined in [G7]  
element
$$
{\rm I}^{\cal H}(\sigma(v_x); \sigma(a_{1}), ..., \sigma(a_{n}); \sigma(v_y)) \in 
{\cal H}_{n}
$$
\end{theorem}

{\bf Proof}. If $x \not = a_1$ and $y \not = a_2$ this is given by theorem 
\ref{2.8.02.1}. 

If $x  = a_1$ or $y  = a_2$ in the Hodge realization 
there are explicit formulas provided by  lemma 6.7 
(and also propositions 2.14, 2.15) 
in [G7] expressing the corresponding element as a product 
of the elements constructed on the first step. Using these explicit formulas as a 
definition in the motivic case we define element (\ref{1=wqa}) and at the same time 
prove that its Hodge realization is as needed. The theorem is proved.

{\bf 3. A geometric construction of the motivic torsor of path between tangential base points}. 
The fundamental Hopf algebra ${\cal A}_{\bullet}(F)$ can be define 
as a commutative Hopf algebra 
${\cal A}_{\bullet}(F)$
in the tensor category ${\cal P}_T$ of pure Tate motives. Indeed, the category 
${\cal P}_T$ is canonically equivalent to the  category  of 
finite dimensional graded $\Q$-vector spaces. 

The dual to the motivic torsor of path (\ref{1}) is  an ind-object 
in ${\cal M}_T(F)$. We define it as an ind-object in the category of comodules over 
${\cal A}_{\bullet}(F)$. For this one needs 

1) to 
construct a pro-object $V$ in the category of pure Tate motives, 

2) to define the coaction map $a: V \lra V \otimes  {\cal A}_{\bullet}(F)$, and 

3) to check the coaction axioms for $a$. 

We do this  as follows. 

1) The $\Q(-n)$-isotipic component $V_n$ of $V$ is given by 
$$
V_n:= \otimes^nH^1({\Bbb P}^1 - \Lambda)
$$
Here, barring the trivial case  $\Lambda = \emptyset$, the motive  
$H^1({\Bbb P}^1 - \Lambda)$ is a pure Tate motive of weight $2$, i.e. a direct sum 
of copies of $\Q(-1)$. Its 
rank is $|\Lambda |-1$. One may assume without loss of generality 
 that 
$\infty \in \Lambda $. Then there is a basis in 
$H^1({\Bbb P}^1 - \Lambda) $ given by the classes
\begin{equation} \label{wewewe}
\omega_a:= d\log (t-a): \Q(-1) \lra H^1({\Bbb P}^1 - \Lambda) ; \qquad a \in 
\Lambda  - \{\infty\}
\end{equation}
 
2) To define the coaction we,  
given a pair of non zero maps
$$
v: \Q(-m) \lra V_m; \quad f: V_n \lra \Q(-n),
$$
will define a map
\begin{equation} \label{1=}
a(v ,f): \Q(n-m) \lra  {\cal A}_{m-n}(F)
\end{equation}
This map is supposed to come as  the composition 
$$
\Q(-m) \stackrel{v}{\lra} V_m \stackrel{a}{\lra} 
V_n \otimes {\cal A}_{m-n}(F) \stackrel{f\otimes {\rm id}}{\lra} 
\Q(-n) \otimes {\cal A}_{m-n}(F)
$$

Choose basis elements
\begin{equation} \label{1=wq}
v = \omega_{b_1} \otimes ... 
\otimes \omega_{b_k};\qquad f = \omega_{a_1} \otimes ... 
\otimes \omega_{a_n};  \qquad a_i, b_j \in \Lambda -\{\infty\}
\end{equation}
 Consider an inclusion of ordered subsets 
\begin{equation} \label{1=wqg}
\beta: \{b_1, ..., b_k\} \hookrightarrow \{a_1, ..., a_n\}
\end{equation} 
provided by  a sequence $0=i_0 < i_1 < ... < i_k < i_{k+1}= n+1$ such that 
$b_p =  a_{i_p}$. 

\begin{definition} \label{2.8.02.5} Assuming (\ref{1=wq}) we set 
$$
a(v,f):= \sum_{\beta} \prod_{p=0}^k{\rm I}^{\cal M}(v_{a_{i_p}}; a_{i_{p}+1}, ..., a_{i_{p+1}-1}; 
v_{a_{i_{p+1}}}) \in 
{\cal A}_{n-k}(F)
$$
where the sum is  over all different 
inclusions (\ref{1=wqg}). 
\end{definition}

This definition has been suggested by theorem 6.4 in [G7]. 

{\bf Examples}. i)  $a(v,f) = 0$ if there is no such inclusion $\beta$. 

ii) Let $v_0: \Q(0) \stackrel{=}{\lra} V_0$ be 
the canonical isomorphism, and $f$ is as above. In this case 
 $\{b_1, ..., b_k\}$ is the empty set, so there is just one inclusion $\beta$. Thus 
$$
a(v_0, f) = {\rm I}^{\cal M}(v_{x}; a_1, ..., a_n ; v_{y})
$$

\begin{theorem} \label{8.1.01.1s} a) Let $F$ be a number field. Then 
the elements 
$a(v,f) \in A_n(F)$ provide an ind-object in ${\cal M}_T(F)$. Its dual is by definition the 
motivic torsor of path ${\cal P}^{\cal
  M}({\Bbb P}^1 -\Lambda; v_x, v_y)$. 

b) There is a morphism of pro-objects in ${\cal M}_T(F)$, the composition of path, 
$$
{\cal P}^{\cal
  M}({\Bbb P}^1 -\Lambda; v_x, v_y) \otimes {\cal P}^{\cal
  M}({\Bbb P}^1 -\Lambda; v_y, v_z) \lra {\cal P}^{\cal
  M}({\Bbb P}^1 -\Lambda; v_x, v_z)
$$

c) The Hodge and l-adic realizations of the motivic torsor of path 
are isomorphic to the ones defined in [D].
\end{theorem}

{\bf Proof}. We use the lemma 3.4 from [G7] based on the injectivity of regulators to deduce 
the parts a) and b)  from the corresponding statements in the Hodge setting. 
Then the Hodge part of c) is then given by the construction.

The Hodge realizations of the motive $H^1(P^1 - \Lambda)$
are parametrized by the complex embeddings 
$\sigma: F \hookrightarrow \C$. The different realizations 
$H^1(\C P^1 - \sigma_i(\Lambda))$ 
are canonically isomorphic to each other. The isomorphism identifies 
the basis elements $\omega_{\sigma_i(a)}$ and $\omega_{\sigma_j(a)}$. 
There is a 
collection of the Hodge-Tate structures 
\begin{equation} \label{weweweq}
{\cal P}^{\cal H}(\C P^1 - \sigma(\Lambda); \sigma(v_x), \sigma(v_y)); \qquad 
\sigma: F \hra \C
\end{equation}
 such that ${\rm Gr}^W_{\bullet}(-)$ 
are identified for different $\sigma$'s. 
Each $\sigma$ 
provides a homomorphism of the Hopf algebras 
$h_{\sigma}: {\cal A}_{\bullet}(F) 
\lra {\cal H}_{\bullet}$. 
The matrix elements $a(v,f)$ are compatible with the Hodge 
realization in the following sense. 

\begin{lemma} \label{8.1.01.1} 
For any complex embedding  
$\sigma$ the element 
$h_{\sigma}(a(v,f))$ coincides with the matrix element of (\ref{weweweq}) 
corresponding to $v$ and $f$ (i.e. to the 
elements of ${\rm Gr}^W_{\bullet}(\ref{weweweq})$  
identified with $v$ and $f$ by the Hodge realization functor). 
\end{lemma}

{\bf Proof}. Follows from theorem 6.4 in [G7]. 
The lemma is proved.

Using lemma \ref{8.1.01.1} and lemma 
3.4 from [G7] we deduce a) and b) 
from  the corresponding fact in the Hodge realization.  

Let us prove the l-adic part of c). In order to follow the described above scheme
 we need only to establish the l-adic version of theorem 6.4 in [G7], which 
is  interesting on its own. 

Let $(a_0; a_1, a_2, ..., a_m; a_{m+1})$ be an 
arbitrary configuration of $F$-points in ${\Bbb A}^1$, where $F$ is an arbitrary field 
with $\mu_{l^{\infty}}\not \in F^*$.  
We define the framed mixed Tate $l$-adic ${\rm Gal}(\overline F/F)$-module  
\begin{equation} \label{CPr2*}
{{\rm I}}^{\rm et}(v_{a_0}; a_1, ..., a_m; v_{a_{m+1}})
\end{equation}
 as the 
l-adic torsor of path ${\cal P}^{(l)}({\Bbb A}^1 - S; v_{a_0}, v_{a_{m+1}})$ where $S:= \cup a_i$, 
framed by $v_0$ and the element $f \in \otimes^mH_{\rm et}^1({\Bbb A}^1 - S)$ 
defined as the image under the l-adic 
regulator of the symbol
$$
\{t-a_1, ..., t-a_m\} \in K_m^M(F(t)) 
$$

\begin{proposition} \label{2.10.02.2}  One has 
\begin{equation} \label{CP2*}
\Delta {{\rm I}}^{\rm et}(v_{a_0}; a_1,  ..., a_m; v_{a_{m+1}})= 
\end{equation}
$$
\sum_{0 = i_0 < i_1 < ... < i_k < i_{k+1} = m} 
{{\rm I}}^{\rm et}(v_{a_0}; a_{i_1}, ..., a_{i_k}; v_{a_{m+1}}) \otimes \prod_{p =0 }^k
{{\rm I}}^{\rm et}(v_{a_{i_{p}}}; a_{i_{p}+1}, ..., a_{i_{p+1}-1}; v_{a_{i_{p+1}}})
$$
\end{proposition}

{\bf Proof}. If $a_0, a_{m+1} \not \in \{a_1, ..., a_m\}$ 
then the corresponding motivic torsor of path is given by, say, Beilinson's construction. 
Thanks to  lemma \ref{11.19.01.10}, if  $F$ is a number field 
is isomorphic to the motivic torsor of path we defined above. 
Since the motivic version of formula (\ref{CP2*}) has been established in this situation, 
and since the l-adic realization of $
{{\rm I}}^{\cal M}(a_0; a_1, ..., a_m; a_{m+1})
$ is isomorphic to its l-adic counterpart (\ref{CPr2*}), we have formula (\ref{CP2*}) 
in this situation. By the specialization theorem it is valid for any 
configuration $(a_0; a_1, a_2, ..., a_m; a_{m+1})$ still assuming that $F$ is a number field. 

\begin{lemma} \label{infr} Let $X$ be a non empty geometrically 
connected  curve and  $M\in {\cal A}^{\rm et}_{w}(X)$ is 
the equivalence class of a lisse l-adic framed mixed Tate object on $X$. 
Then if restriction of $M$ to infinitely many 
points of $X$ is zero then  $M =  0$. 
\end{lemma}

{\bf Proof}. We proceed by induction. For $w=1$ this is clear since 
${\cal A}^{\rm et}_{1}(X) = {\cal O}^*(X)\otimes \Q_l$. 
Recall  
the following rigidity property ([BD], s. 1.7): for $n>1$ 
restriction to 
any point $x: {\rm Spec} F \hra U$ 
provides an isomorphism
\begin{equation} \label{2.10.02.11}
{\rm Ext}^1_X(\Q_l(0), \Q_l(n)) \stackrel{\sim}{\lra} {\rm Ext}^1_{{\rm Spec} F}(\Q_l(0), \Q_l(n))
\end{equation}
 Recall that the kernel of the restricted coproduct $\Delta'$ on ${\cal A}^{\rm et}_n(X)$
is identified with ${\rm Ext}^1_X(\Q_l(0), \Q_l(n))$. 
The condition of the lemma implies by induction that $\Delta' (M) =0$, 
and hence $M=0$ by the rigidity. The lemma is proved. 

  It follows from this that formula (\ref{CP2*}) holds for any field $F$ as above. 
The proposition is proved. 

The theorem is proved.

\section{The motivic shuffle relations}

We say that  $(x_1, ..., x_m), (n_1, ..., n_m)$ are the parameters of  
${\rm Li}_{n_1, ..., n_{m}}
(x_1, ..., x_{m})$. The parameters are admissible if
\begin{equation} \label{3.7.01.1qwe}
x_m \not = 1 \quad \mbox{or} \quad n_m >1
\end{equation}
For admissible parameters we have defined in [G7] a framed Hodge-Tate structure 
corresponding to multiple polylogarithms by  
\begin{equation} \label{1.16.02.11}
{\rm Li}^{{\cal H}}_{n_1, ..., n_{m}}(x_1, ..., x_{m}) := (-1)^m \cdot
{\rm I}^{{\cal H}}_{n_1, ..., n_{m}}(a_1, a_2, ..., a_m)
\end{equation}
$$
a_1:= (x_1 ... x_m)^{-1}, a_2 := (x_2 ... x_m)^{-1}, ..., a_m:= x^{-1}_m
$$
where  the right hand side is the framed Hodge-Tate structure 
corresponding to the iterated integral (\ref{5*zx}). Thanks to the results 
of previous chapters, e.g. theorem \ref{12.06.01.10} and \ref{1.21.01.20}, 
in the case of admissible parameters 
 the Hodge-Tate structure ${\rm I}^{{\cal H}}$ defined in 
[G7] via the specialization is equivalent 
to the one defined in chapter 2. In particular it is equivalent to a 
Hodge-Tate structure of geometric origin. 

There are two different ways to extend the definition 
of  the ${\rm Li}^{{\cal H}}$- and ${\rm I}^{{\cal H}}$- 
Hodge-Tate structures to the set of  non admissible parameters, 
keeping the corresponding  shuffle relations. 
In each case the definition is provided by the fiber at the tangent vector 
$\partial/\partial \varepsilon$ of 
the specialization functor applied to an appropriate variation of 
the framed multiple polylogarithm 
Hodge-Tate structures over a small punctured disc 
with a natural coordinate $\varepsilon$. For the ${\rm I}^{\cal H}$-elements 
such a procedure has been worked out in [G7], and called there the canonical regularization. For the ${\rm Li}^{\cal H}$-elements this is done  in this section.

It is essentially obvious 
in the ${\rm I}^{\cal H}$-case, and follows from proposition \ref{1.7.02.1} 
in the ${\rm Li}^{\cal H}$-case,  
that defined this way ${\rm I}^{\cal H}$- and 
${\rm Li}^{\cal H}$-generators satisfy the 
corresponding shuffle relations for all parameters. 
However the basic formula (\ref{1.16.02.11}) comparing the  
${\rm I}^{\cal H}$- and ${\rm Li}^{\cal H}$-generators for admissible  parameters 
has to be replaced by a more sophisticated relation - 
see theorem \ref{1.13.02.5}.  The proof shows that 
the comparison formula from theorem \ref{1.13.02.5} is uniquely 
determined if we want to keep the shuffle product formula 
for all parameters.

Our approach to regularization of the  ${\rm Li}$-generators 
is contained  in the last page  of [G3], see also s. 2.10 in [G7]. 
We show below that, using proposition \ref{1.7.02.1}, it is   easily 
transformed to the Hodge or \'etale setting.

For the multiple $\zeta$-numbers a different regularization of the ${\rm Li}$-power series
was considered by Zagier  and documented by Ihara-Kaneko [IK]. 
A  similar regularization was  
used by  Boute de Monville. 
However it is not clear how to  make this regularization 
  Hodge theoretic or motivic.

For multiple polylogarithms, working 
on the level of numbers, it is not quite clear how even to formulate 
an explicit version of 
the comparison problem for two regularizations. Indeed, 
the values of multiple polylogarithms are not quite well defined numbers 
since the corresponding functions are multivalued.

{\bf 1. The power series shuffle product formula}. Let us define a  
generalized shuffle of the two ordered sets 
\begin{equation} \label{11.22.00.2}
\{x_1,...,x_{m}\} \quad \mbox{and} \quad \{y_{1},...,y_{n}\}
\end{equation}
Take a string  of points on the real line, called slots,  and mark 
every slot either by $x_i$, or by $y_j$, or by $x_i$ and $y_j$,  in 
such a way that 
$x_a$ (resp. $y_a$) is on the left of $x_b$ (resp. $y_b$) if $a<b$. 
The generalized shuffles 
are the  combinatorial types  of configurations 
of labelled slots obtained this way. An example is presented on the picture. 

\begin{center}
\hspace{4.0cm}
\epsffile{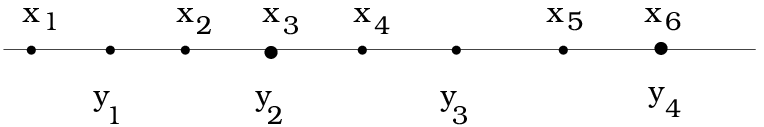}
\end{center}

For instance if $m=n=1$ there are three generalized shuffles: 
$(x_1; y_1)$, $(y_1; x_1)$ and  $(x_1, y_1)$. The last one  
has  just one slot where both $x_1$ and $y_1$ sit, while in 
the first two there are two different slots.

Let $\overline \Sigma_{p,q}$ be the set of all generalized shuffles 
of the ordered sets $\{1, ..., p\}$ and $\{p+1, ..., p+q\}$. 
Set $$
\Z_{++}^p:= \{(k_1, ..., k_p) \in \Z^p_+ \quad | \quad 0 < k_1 < ... < k_p\}, 
$$ 
Then there is a natural decomposition of the product of cones 
$\Z_{++}^{p} \times \Z_{++}^{q}$ into a union of the cones $Z_{++}^{\sigma}$ 
parametrized by the generalized shuffles:
\begin{equation} \label{4.16.01.2}
\Z_{++}^{p} \times \Z_{++}^{q} = \cup_{\sigma \in \overline \Sigma_{p,q}} Z_{++}^{\sigma}
\end{equation}
For example for $p=q=1$ we have 
$$
\{k_1 >0\} \times \{k_2 >0\}  = \{0 < k_1<  k_2 \} \cup \{0 < k_1=  k_2 \}
\cup \{k_1> k_2 >0\}
$$

For a generalized shuffle $\sigma \in \overline \Sigma_{p+q}$ consider the formal power series, convergent if $|x_i|<1$: 
$$
{\rm Li}^{\sigma}_{n_1, ..., n_{p+q}}(x_1, ..., x_{p+q}):= 
\sum_{(k_1, ..., k_{p+q}) \in Z^{\sigma}_{++}}\frac{x^{k_1}_1 ... 
x^{k_{p+q}}_{p+q}}{k_1^{n_1}... k_{p+q}^{n_{p+q}}}
$$
Then it follows from (\ref{4.16.01.2}) 
that there is an equality of formal power series, 
the power series shuffle product formula:

\begin{equation} \label{4.16.01.13}
{\rm Li}_{n_1, ..., n_{p}}(x_1, ..., x_{p}) \cdot 
{\rm Li}_{n_{p+1}, ..., n_{p+q}}(x_{p+1}, ..., x_{p+q}) = 
\end{equation}
$$
\sum_{\sigma \in \overline \Sigma_{p,q}}{\rm Li}^{\sigma}_{n_1, ..., n_{p+q}}(x_1, ..., x_{p+q})
$$

Let $\sigma \in \overline \Sigma_{p,q}$ be a generalized shuffle. Denote by 
$\sigma(n_1, ..., n_{p+q})$ the sequence of integers obtained from 
$n_1, ..., n_{p+q}$ using 
$\sigma$ and the following convention: the contribution of a slot 
of $\sigma$ into 
$\sigma(n_1, ..., n_{p+q})$ is the sum of the integers $n_i$ 
sitting at this slot. We define  $\sigma(x_1, ..., x_{p+q})$ similarly 
using the  convention that a slot 
of $\sigma$ contributes the product of $x_i$'s sitting at this slot. 
For example if $\sigma$ is the generalized shuffle of the sets $\{1\}$ and 
$\{2\}$ with just one slot 
then 
${\rm Li}_{\sigma(1, 1)}(\sigma(x_1, x_2)) = 
{\rm Li}_{2}(x_1 x_2)$. 

Using this convention we can write the right hand side of (\ref{4.16.01.13}) 
as 
\begin{equation} \label{1.20.02.1}
\sum_{\sigma \in \overline \Sigma_{p,q}}
{\rm Li}_{\sigma(n_1, ..., n_{p+q})}(\sigma(x_1, ..., x_{p+q})) 
\end{equation}

{\bf 2. The category of unipotent variations of Hodge-Tate structures 
on  ${\C}^*$}. 
Recall that it  is canonically equivalent to the category of graded 
finite dimensional comodules over the corresponding fundamental  Hopf algebra, denoted  
$
{\cal A}^{\cal H}_{\bullet}({\Bbb A}^1 - \{0\})
$, (see for example [G7] ch. 3). There is  the canonical element 
$$
\log^{\cal H}\varepsilon \in  {\cal A}^{\cal H}_{1}({\Bbb A}^1 - \{0\}) =
{\rm Ext}_{MHS({\Bbb A}^1 - \{0\})}^1(\Q(0), \Q(1)) 
$$

\begin{lemma} \label{1.13.02.2}
 There is an isomorphism of graded commutative algebras  
$$
{\cal A}^{\cal H}_{\bullet}({\Bbb A}^1 - \{0\}) = {\cal H}_{\bullet} \otimes_{\Q} 
\Q[\log^{\cal H}\varepsilon]
$$
\end{lemma}

{\bf Proof}. Indeed, let $X$ be a smooth complex algebraic variety and  $L^{\cal H}_{\bullet}(X(\C))$ 
 the fundamental Lie algebra of the category of unipotent variations of Hodge-Tate structures 
on $X(\C)$. Then there is an exact sequence of Lie algebras 
$$
0 \lra L^{\rm geom}(X(\C)) \lra L^{\cal H}_{\bullet}(X(\C)) \lra L^{\cal H}_{\bullet}({\rm Spec}(\C)) \lra 0
$$
where $L^{\rm geom}(X(\C))$ is the geometric fundamental Lie algebra. It is isomorphic to the 
Lie algebra provided by the 
pronilpotent completion  of the topological fundamental group. Therefore in the case $X = {\Bbb G}_m$ 
we have $L^{\rm geom}(\C^*) = \Q(1)$. 
Applying the Poincare-Birkhoff-Witt theorem we get the lemma.

{\bf 3. The framed Hodge-Tate structures corresponding to multiple polylogarithms}. 
We want to  define the framed Hodge-Tate structures  
$  {\rm Li}^{\cal H}_{n_1, ..., n_m}
(x_1, ..., x_m)$ for all $x_i$ and all $n_i >0$ so that 
they coincide with (\ref{1.16.02.11}) for admissible parameters, and satisfy 
the  shuffle relations (\ref{4.16.01.13})-(\ref{1.20.02.1}) for all all parameters. 

\begin{lemma} \label{1.7.01.10}
Let $U^*$ be a punctured  at zero disc. Suppose that 
$x_i(\varepsilon)$ are holomorphic for $\varepsilon \in U^*$, and $n_m >1$ or 
$x_m(\varepsilon) \not = 1$ for small non zero $\varepsilon$. 
Then for a sufficiently small disc $U^*$ the 
framed Hodge-Tate structures 
\begin{equation} \label{1.7.01.11}
{\rm Li}^{\cal H}_{n_1, ..., n_m}
(x_1(\varepsilon ), ...,  x_m(\varepsilon)); \quad \varepsilon \in U^*
\end{equation}
are equivalent to fibers of a certain unipotent variation of Hodge-Tate structures 
over $U^*$.  
\end{lemma}

{\bf Proof}. This is a standard general fact. In our case it 
follows immediately from theorem 5.5 in [G7]. 
Indeed, translating the assertion of that theorem from the ${\rm I}-$ 
to the ${\rm Li}-$ notation we get the  following. Say that 
$(x_1, ..., x_m)$ and $(x'_1, ..., x'_m)$ are of the same combinatorial type if 
$x_i ... x_j = 1$ if and only if $x'_i ... x_j' =1$. Then the canonical 
(in the sense of [G7])
Hodge-Tate structures  ${\rm Li}^{\cal H}_{n_1, ..., n_m}
(x_1, ...,  x_m)$ form a unipotent variation over the space of all configurations 
of points $(x_1, ..., x_m)$ of given combinatorial type. Making 
 $U^*$  smaller we can assume that 
 $(x_1(\varepsilon ), ...,  x_m(\varepsilon))$ has the same combinatorial type 
for all $\varepsilon \in U^*$. The lemma is proved.

If $x_1 ... x_m \not = 0$, by lemma \ref{1.7.01.10} there is a 
 unipotent variation framed Hodge-Tate structures on  a 
little punctured complex disc   with a natural 
 parameter $\varepsilon$ whose fibers at points $\varepsilon$ are 
\begin{equation} \label{1.12.02.6}
{\rm Li}^{{\cal H}}_{n_1, ..., n_{m}}(x_1(1-\varepsilon), ..., x_{m}(1-\varepsilon))
\end{equation} 
Observe that the parameters of these elements for small $\varepsilon \not = 0$ are admissible. 
Its specialization at $\varepsilon =0$ 
is a unipotent variation of framed Hodge-Tate structures on 
the punctured tangent space. 
It provides the  element  
$$
\widehat {\rm Li}^{{\cal H}}_{n_1, ..., n_{m}}(x_1, ..., x_{m}) \in 
{\cal A}^{\cal H}_{\bullet}({\Bbb A}^1 - \{0\})  
$$
 The tangent space is  
 equipped  with  coordinate $\varepsilon$. By lemma \ref{1.13.02.2} we have 
\begin{equation} \label{1.12.02.5}
\widehat {\rm Li}^{{\cal H}}_{n_1, ..., n_{m}}(x_1, ..., x_{m}) = \sum_{k\geq 0}
{\rm Li}^{{\cal H}}_{n_1, ..., n_{m}}(x_1, ..., x_{m})_{(k)} \cdot (\log^{{\cal H}}
\varepsilon)^k
\end{equation} 
 We define 
${\rm Li}^{{\cal H}}_{n_1, ..., n_{m}}(x_1, ..., x_{m})$ as 
 the constant term of this expression.   We set 
$$
\widehat {\rm Li}^{{\cal H}}_{n_1, ..., n_{m}}(x_1, ..., x_{m})=0 \quad \mbox{if $x_1 ... x_m =0$}
$$


\begin{proposition} \label{1.8.02.1} For admissible parameters we have
$$
\widehat {\rm Li}^{\cal H}_{n_1, ..., n_m}
(x_1, ..., x_m) =  {\rm Li}^{\cal H}_{n_1, ..., n_m}
(x_1, ..., x_m)
$$
\end{proposition}

{\bf Proof}. Follows immediately from the specialization theorem \ref{10.24.01.1}.

\begin{theorem} \label{4.16.01.13saq} For any $x_i \in \C$ and any positive integers $n_i$ 
then there is a
 shuffle product formula 
\begin{equation} \label{3.7.01.1132pp}
{\rm Li}^{{\cal H}}_{n_1, ..., n_{p}}(x_1, ..., x_{p}) 
 \cdot {\rm Li}^{{\cal H}}_{n_{p+1}, ..., n_{p+q}}
(x_{p+1}, ..., x_{p+q}) = 
\end{equation}
$$
\sum_{\sigma \in \overline \Sigma_{p,q}}
{\rm Li}^{{\cal H}}_{\sigma(n_1, ..., n_{p+q})}(\sigma(x_1, ..., x_{p+q})) 
$$
There is a similar formula with $ {\rm Li}^{{\cal H}}$ 
replaced by  $\widehat {\rm Li}^{{\cal H}}$. 
\end{theorem}

{\bf 4. The strategy of the proof of theorem \ref{4.16.01.13saq}}. 
The proof  consists of two steps of different nature. 

{\bf Step 1}. 
\begin{theorem} \label{4.17.01.77z} 
 Suppose that the parameters $(x_1, ..., x_m)$ of each of the terms in 
 (\ref{3.7.01.1132pp}) satisfy the following condition: 
\begin{equation} \label{4.11.01.3} 
x_i \not = 0, \quad 1 \not =  \{x_m^{-1}, \quad (x_{m-1}x_m)^{-1},\quad  ... \quad 
(x_1 ... x_{m})^{-1}\}
\end{equation}
Then the 
 shuffle product formula (\ref{3.7.01.1132pp}) is valid. 
\end{theorem}

We will prove this theorem in chapters 9  below. 
An  essentially different proof will be outlined in chapter 11. 
Now let us assume theorem \ref{4.17.01.77z} and proceed to the proof of theorem 
\ref{4.16.01.13saq}.

\begin{corollary}\label{1.20.02.2} The shuffle relations (\ref{3.7.01.1132pp}) are valid if 
all the terms in (\ref{3.7.01.1132pp}) 
are admissible. 
\end{corollary}

{\bf Proof}. Suppose that parameters $(x_1, ..., x_m)$ satisfy 
 the condition of the corollary. One can find a 
little complex curve $U$ passing through the point $(x_1, ..., x_m)$ 
 such that the corresponding punctured curve $U^*$
lies in  the domain determined by conditions (\ref{4.11.01.3}).  
For sufficiently small $U$ there are unipotent variations over 
$U^*$ corresponding to the terms of (\ref{3.7.01.1132pp}), and by theorem \ref{4.17.01.77z} 
we have identity 
(\ref{3.7.01.1132pp}) for the corresponding framed variations. Taking specialization, and using the specialization theorem, we get the corollary. The corollary is proved.

{\bf Remark}. To get  corollary \ref{1.20.02.2} 
it is sufficient to have formula (\ref{3.7.01.1132pp}) 
for any non empty Zariski open subset of the parameters set. 

{\bf Step 2}.  

\begin{proposition} \label{1.7.02.1}
Suppose that $x_i(\varepsilon) \in \C^*$ are holomorphic for  small $\varepsilon$, 
$x_{l+1}(0) = ... = x_m(0) = 1$, and $n_l>1$ or $x_l(0) \not = 1$.  
Assume in addition that $x_m(\varepsilon) \not = 1$ for small non zero $\varepsilon$, 
or $n_m >1$. 
Then 
$$
[{\rm Sp}_{\varepsilon = 0}]{\rm Li}^{\cal H}_{n_1, ..., n_l, 1, ..., 1}
(x_1(\varepsilon ), ..., x_l(\varepsilon ), x_{l+1}(\varepsilon), ..., x_m(\varepsilon )) = 
$$
$$
[{\rm Sp}_{\varepsilon = 0}]
{\rm Li}^{\cal H}_{n_1, ..., n_l, 1, ..., 1}
(x_1(0), ..., x_l(0), x_{l+1}(\varepsilon), ..., x_m(\varepsilon )) 
$$
\end{proposition}

{\bf Proof}. Thanks to  lemma \ref{1.7.01.10} 
 the use of  the  functor ${\rm Sp}_{\varepsilon = 0}$ is legitimate.

We prove the proposition by induction on $k:= m-l$. 
If $k=0$ it follows from the specialization theorem \ref{10.24.01.1}. 
The induction step is deduced from theorem \ref{4.17.01.77z}. 
Namely, consider the shuffle formula for the product 
\begin{equation} \label{1.1.02.100}
{\rm Li}^{\cal H}_{n_1, ..., n_l}
(x_1(\varepsilon ), ..., x_l(\varepsilon )) \cdot {\rm Li}^{\cal H}_{1, ..., 1}
( x_{l+1}(\varepsilon), ..., x_m(\varepsilon )), \quad \varepsilon  \not = 0
\end{equation}
Thanks to the condition in the proposition all the terms in the shuffle product formula are admissible for $\varepsilon  \not = 0$, so the shuffle product formula 
is available by  theorem \ref{4.17.01.77z}. 
This formula is a sum of several terms. One of them is  
\begin{equation} \label{1.8.02.21}
{\rm Li}^{\cal H}_{n_1, ..., n_l, 1, ..., 1}
(x_1(\varepsilon), ..., x_l(\varepsilon),  x_{l+1}(\varepsilon), ..., x_m(\varepsilon))
\end{equation}
In the rest of the terms  $x_{l+1}(\varepsilon)$ stays strictly 
to the left of 
$x_l(\varepsilon)$, or appears in the 
variable $x_l(\varepsilon) x_{l+1}(\varepsilon)$ 
coupled to the index $n_{l+1}+1$. 
For each of these terms our induction  invariant $k$ is down at least by one. Thus  
by the induction assumption applying ${\rm Sp}_{\varepsilon = 0}$ 
we can replace all the variables 
on the left of $x_l(\varepsilon)$ by their values at $\varepsilon =0$. 

On the other hand,  thanks to the condition $x_l(0) \not = 1$ or $n_l >1$, 
the specialization theorem implies 
$$
[{\rm Sp}_{\varepsilon = 0}]{\rm Li}^{\cal H}_{n_1, ..., n_l}
(x_1(\varepsilon ), ..., x_l(\varepsilon )) = {\rm Li}^{\cal H}_{n_1, ..., n_l, }
(x_1(0), ..., x_l(0))
$$
So the product (\ref{1.1.02.100}) equals to 
$$
{\rm Li}^{\cal H}_{n_1, ..., n_l}
(x_1(0 ), ..., x_l(0 )){\rm Li}^{\cal H}_{1, ..., 1}
( x_{l+1}(\varepsilon), ..., x_m(\varepsilon ))
$$
One of the terms in the shuffle formula for this product is 
\begin{equation} \label{1.8.02.22}
{\rm Li}^{\cal H}_{n_1, ..., n_l, 1, ..., 1}
(x_1(0), ..., x_l(0),  x_{l+1}(\varepsilon), ..., x_m(\varepsilon))
\end{equation}
By induction we already know that the specialization of 
the other terms equal to the specialization of the ones 
in the shuffle formula for (\ref{1.1.02.100}). Thus 
$$
[{\rm Sp}_{\varepsilon = 0}](\ref{1.8.02.21})= [{\rm Sp}_{\varepsilon = 0}](\ref{1.8.02.22})
$$ 
The proposition is proved. 

Now we can finish the proof of theorem \ref{4.16.01.13saq}. 
The case when $x_1 ... x_m =0$ is trivial since all the terms in 
(\ref{3.7.01.1132pp}) are zero. 
So we may assume $x_1 ... x_m \not = 0$. Then by theorem \ref{4.17.01.77z} we have 
the shuffle product formula 
$$
 {\rm Li}^{{\cal H}}_{n_1, ..., n_{p}}(x_1(1-\varepsilon), ..., x_{p}(1-\varepsilon)) 
 \cdot  {\rm Li}^{{\cal H}}_{n_{p+1}, ..., n_{p+q}}
(x_{p+1}(1-\varepsilon), ..., x_{p+q}(1-\varepsilon)) = 
$$
\begin{equation} \label{3.7.01.1132q}
\sum_{\sigma \in \overline \Sigma_{p,q}}
{\rm Li}^{{\cal H}}_{\sigma(n_1, ..., n_{p+q})}
(\sigma(x_1(1-\varepsilon), ..., x_{p+q}(1-\varepsilon))) 
\end{equation}
Applying to it the functor $[{\rm Sp}_{\varepsilon = 0}]$ we get a valid formula. 
It remains to match the terms of that formula with the ones  
of  (\ref{3.7.01.1132pp}). For the left hand side of (\ref{3.7.01.1132q}),
 as well as for the terms in  the right hand side of depth $p+q$,  
this is clear from the very definition. 
For the depth $<p+q$ terms in the right hand side this is deduced from 
proposition \ref{1.7.02.1}. Indeed, there exists a slot of a given generalized shuffle 
where two indices sit: $n_i$ and $n_j$. Take the very right such a slot. 
Since $n_i+n_j >1$ 
we may replace the corresponding variables $x_s(1-\varepsilon)^{a(s)}$ 
at this slot and all the slots to the left of it 
by $x_s$,  as well as by $x_s(1-\varepsilon)$. 
Theorem \ref{4.16.01.13saq} is proved assuming 
theorem \ref{4.17.01.77z}. 

{\bf 5. The canonical specialization and the comparison theorem}. 
Consider the 
unipotent variation framed Hodge-Tate structures over a small 
punctured complex disc with the coordinate 
$\varepsilon$, whose fiber at $\varepsilon$ is 
\begin{equation} \label{1.13.02.3}
{\rm Li}^{{\cal H}}_{n_1, ..., n_{m}}(x_1, ..., x_{m-1}, 
x_{m}(1-\varepsilon))
\end{equation} 
Observe that 
\begin{equation} \label{1.15.02.1}
{\rm Li}^{{\cal H}}_{n_1, ..., n_{m}}(x_1, ..., x_{m-1}, 
x_{m}(1-\varepsilon))= 
\end{equation}
$$
 {\rm I}^{{\cal H}}_{n_1, ..., n_{m}}(0; (x_1 ... x_{m})^{-1}, 
..., 
x^{-1}_{m}; (1-\varepsilon))
$$

So according to [G7] the fiber of its specialization at 
$\partial/\partial \varepsilon$  serves as  
the definition of the Hodge-Tate structure 
${\rm Li}^{{\cal H}}_{n_1, ..., n_{m}}(x_1, ..., x_{m})$ 
for arbitrary set of parameters. 
This procedure was called  in [G7] the canonical regularization.

The  specialization of the family 
(\ref{1.13.02.3}) provides  the element
$$
\widehat {\rm I}^{{\cal H}}_{n_1, ..., n_{m}}(a_1, ..., a_{m})\in {\cal A}^{\cal H}_{\bullet}({\Bbb A}^1 - \{0\})  
$$

\begin{theorem} \label{1.13.02.5} For any $x_i \in \C$ one has 
$$
\widehat {\rm Li}^{{\cal H}}_{n_1, ..., n_{m}}(x_1, ..., x_{m}) = {\Bbb L} \circ 
\widehat {\rm I}^{{\cal H}}_{n_1, ..., n_{m}}(a_1, ..., a_{m})
$$
\end{theorem} 

{\bf Proof}. Let us establish first the case when $n_i =1$ and $x_i =1$ 
for all $i$. 
Observe that (see [G7]) 
$$
\widehat {\rm Li}^{{\cal H}}_{1, ..., 1, 1}
(\underbrace{1, ...,1, 1-\varepsilon}_{n } ) = \frac{(\log^{\cal H}(\varepsilon))^n}{n!}
$$

\begin{lemma} \label{1.14.02.1}
\begin{equation}\label{1.13.02.1}
\sum_{n \geq 0}\widehat {\rm Li}^{{\cal H}}_{1, ..., 1}
(\underbrace{1, ...,1}_{n } ) \cdot u^n =  {\rm exp}\left( - \sum_{n=1}^{\infty} 
(-1)^n \frac{\widehat \zeta^{\cal H}(n)}{n}u^n\right)
\end{equation}
\end{lemma}

{\bf Proof}. Applying $\frac{d}{du}$ to 
 both parts of 
 to (\ref{1.13.02.1}) and replacing the exponential 
factor by the left hand side of (\ref{1.13.02.1}) we get the following corollary of 
(\ref{1.13.02.1})
\begin{equation}\label{1.13.02.1d}
\sum_{k,l \geq 0}(-1)^k \widehat {\rm Li}^{{\cal H}}_{k+1}(1)\cdot 
\widehat {\rm Li}^{{\cal H}}_{1, ..., 1}
(\underbrace{1, ...,1}_{l} ) u^{k+l}\stackrel{?}{=}
 \sum_{m \geq 0}(m+1) \widehat {\rm Li}^{{\cal H}}_{1, ..., 1}
(\underbrace{1, ...,1}_{m+1} )\cdot u^m
\end{equation}
Arguing by induction on $n$ we see that it  is equivalent to 
the original identity (\ref{1.13.02.1}). Indeed,  
all the terms in the exponential have degree at least $1$, so we 
can use the induction assumption to replace the exp-term in the formula 
for $\frac{d}{du}$(the right hand side of (\ref{1.13.02.1d}))  by 
the left hand side of (\ref{1.13.02.1}).

Formula (\ref{1.13.02.1d})  
is equivalent to collection  of identities, one for each $m$: 
$$
\sum_{k= 0}^m(-1)^k \widehat {\rm Li}^{{\cal H}}_{k+1}(1)\cdot 
\widehat {\rm Li}^{{\cal H}}_{1, ..., 1}
(\underbrace{1, ...,1}_{m-k} ) = (m+1) \widehat {\rm Li}^{{\cal H}}_{1, ..., 1}
(\underbrace{1, ...,1}_{m+1} )
$$
Writing the shuffle relations (available by 
theorem \ref{4.16.01.13saq}),  for each term in the sum, and 
taking the sum, we get the last identity. The lemma is proved. 

{\bf Examples}. We have 
$$
\widehat {\rm Li}^{{\cal H}}_{1}(1)\cdot 
\widehat {\rm Li}^{{\cal H}}_{1}
(1)  - \widehat{\rm Li}^{{\cal H}}_{2}(1)= 2 \widehat {\rm Li}^{{\cal H}}_{1, 1}(1,1)
$$
$$
\widehat {\rm Li}^{{\cal H}}_{1}(1)\cdot 
\widehat {\rm Li}^{{\cal H}}_{1,1}
(1,1)  - \widehat{\rm Li}^{{\cal H}}_{2}(1)\widehat {\rm Li}^{{\cal H}}_{1}
(1) + \widehat {\rm Li}^{{\cal H}}_{3}(1)= 3 \widehat {\rm Li}^{{\cal H}}_{1, 1, 1}(1,1, 1)
$$
So
$$
\widehat {\rm Li}^{{\cal H}}_{1, 1}(1,1) = -\frac{{\rm Li}^{{\cal H}}_{2}(1)}{2} 
+ 
\frac{(\log \varepsilon)^2}{2}
$$
$$
\widehat {\rm Li}^{{\cal H}}_{1, 1, 1}(1,1, 1) = \frac{{\rm Li}^{{\cal H}}_{3}(1) }{3}
+ \frac{{\rm Li}^{{\cal H}}_{2}(1) \log \varepsilon}{2} - 
\frac{(\log \varepsilon)^3}{6}
$$

Now let us treat the general case. We present 
${\rm Li}^{{\cal H}}_{n_1, ..., n_{m}}(x_1, ..., x_m)$ as 
$$
{\rm Li}^{{\cal H}}_{n_1, ..., n_{l}, 1, ..., 1}(x_1, ..., x_l,1,...,1)
$$
 where $n_l>1$ or $x_l \not = 1$. 
We will prove the theorem by induction on  $m-l$. 

By proposition \ref{1.7.02.1} $\widehat 
{\rm Li}^{{\cal H}}_{n_1, ..., n_{m}}(x_1, ..., x_m)$ is the specialization of the 
variation 
$$
{\rm Li}^{{\cal H}}_{n_1, ..., n_{l}, 1, ..., 1}
(x_1, ..., x_l,1 - \varepsilon, ... , 1 - \varepsilon)
$$
In particular this settles the  case $m-l=0$. Lemma 
\ref{1.14.02.1} just means that 
the theorem is true for $ 
{\rm Li}^{{\cal H}}_{1, ..., 1}(1, ..., 1 )$. 
Consider the shuffle  product formulas
\begin{equation} \label{1.14.02.2}
{\rm Li}^{{\cal H}}_{n_1, ..., n_{l}}
(x_1, ..., x_l,)\cdot 
{\rm Li}^{{\cal H}}_{1, ..., 1}
(1 - \varepsilon, ... , 1 - \varepsilon) = 
\end{equation}
\begin{equation} \label{1.14.02.4}
{\rm Li}^{{\cal H}}_{n_1, ..., n_{l}, 1, ..., 1}
(x_1, ..., x_l,1 - \varepsilon, ... , 1 - \varepsilon)+ \mbox{the rest of the terms}
\end{equation}
as well as 
\begin{equation} \label{1.14.02.3}
{\rm Li}^{{\cal H}}_{n_1, ..., n_{l}}
(x_1, ..., x_l,)\cdot 
{\rm Li}^{{\cal H}}_{1, ..., 1}
(1, ..., 1, 1 - \varepsilon)
\end{equation}
\begin{equation} \label{1.14.02.5}
{\rm Li}^{{\cal H}}_{n_1, ..., n_{l}, 1, ..., 1}
(x_1, ..., x_l,1, ... , 1, 1 - \varepsilon)+ \mbox{the rest of the terms}
\end{equation}
Let us apply the $<{\rm Sp}_{\varepsilon = 0}>$ functor to these 
identities. 
Thanks to the lemma we have $(\ref{1.14.02.2}) = {\Bbb L}\circ (\ref{1.14.02.3})$.
By the induction assumption we have 
$$
[{\rm Sp}_{\varepsilon=0}]\mbox{(the rest of the terms in (\ref{1.14.02.4}))} 
=
$$
$$ {\Bbb L}\circ
[{\rm Sp}_{\varepsilon=0}]\mbox{(the rest of the terms in (\ref{1.14.02.5}))}
$$
It follows that 
$$
[{\rm Sp}_{\varepsilon=0}]{\rm Li}^{{\cal H}}_{n_1, ..., n_{l}, 1, ..., 1}
(x_1, ..., x_l,1 - \varepsilon, ... , 1 - \varepsilon)= 
$$
$$
{\Bbb L}\circ [{\rm Sp}_{\varepsilon=0}]{\rm Li}^{{\cal H}}_{n_1, ..., n_{l}, 1, ..., 1}
(x_1, ..., x_l,1 , ... , 1, 1 - \varepsilon)
$$
The theorem is proved. 

\begin{corollary} \label{1.17.02.1}
There exists an explicit formula expressing 
${\rm Li}^{{\cal H}}_{n_1, ..., n_{m}}(x_1, ..., x_{m})$ via 
${\rm Li}^{{\cal H}}_{*}(-)$ with admissible parameters.
\end{corollary}

{\bf Proof}. Theorem \ref{1.13.02.5} (and a shuffle product formula) 
provides a formula expressing 
${\rm Li}^{{\cal H}}_{n_1, ..., n_{m}}(x_1, ..., x_{m})$ 
as $\Q$-linear combinations 
of  ${\rm I}^{{\cal H}}_{*}(-)$ with not necessarily admissible parameters. 
Non admissible  ${\rm I}^{{\cal H}}_{*}(-)$'s are expressed 
as linear combinations of admissible ones using explicit formulas 
from lemma 6.7 (or formulas from proposition 2.14 and 2.15) in [G7]. 
It remains to use formula  (\ref{1.16.02.11}). 
The corollary follows. 

{\bf 6. The shuffle relations on the motivic level}. 
\begin{theorem} \label{4.16.01.13sa1} Let $F$ be a number field. Suppose that 
$x_i \in F$. Then

a) There exists a framed mixed Tate motive over $F$ 
\begin{equation} \label{1.8.01.33}
{\rm Li}^{{\cal M}}_{n_1, ..., n_{m}}(x_1, ..., x_{m}) \in {\cal A}_w(F)
\end{equation}
so that for any embedding $\sigma: F \hra \C$ 
its Hodge realization coincides with 
$$
 {\rm Li}^{{\cal H}}_{n_1, ..., n_{m}}(\sigma(x_1), ..., \sigma(x_m)) \in {\cal H}_w
$$

b) 
There  is a
 shuffle product formula 
\begin{equation} \label{3.7.01.1132}
 {\rm Li}^{{\cal M}}_{n_1, ..., n_{p}}(x_1, ..., x_{p}) 
 \cdot  {\rm Li}^{{\cal M}}_{n_{p+1}, ..., n_{p+q}}
(x_{p+1}, ..., x_{p+q}) = 
\end{equation}
$$
\sum_{\sigma \in \overline \Sigma_{p,q}}
{\rm Li}^{{\cal M}}_{\sigma(n_1, ..., n_{p+q})}(\sigma(x_1, ..., x_{p+q})) 
$$
\end{theorem}

{\bf Proof}. a) If the parameters of (\ref{1.8.01.33}) are admissible the
 construction is given by the results of chapters 2 and 3. In particular its Hodge realization 
is the same as the one defined in [G7] thanks  to the specialization theorem. 

If the parameters are not admissible, we use  the explicit formulas provided by 
corollary \ref{1.17.02.1} (with 
${\rm Li}^{{\cal H}}$ changed to ${\rm Li}^{{\cal M}}$),  to define 
(\ref{1.8.01.33}). The compatibility with the Hodge realization is then obvious. 

 b) By theorem \ref{4.16.01.13saq} 
the shuffle product formula is valid on the Hodge level. 
Therefore a) and lemma 3.4 from [G7] provide (\ref{3.7.01.1132}). 
Theorem \ref{4.16.01.13sa1} is proved. 

{\bf 7. The shuffle relations in the $l$-adic setting}. 
\begin{theorem} \label{4.16.01.13sa} 
Let $F$ be a field and $\mu_{l^{\infty}} \not \in F$. Suppose that 
$x_i \in F$. Then

a) There exists a framed l-adic mixed Tate ${\rm Gal}(\overline F/F)$-module 
\begin{equation} \label{1.8.01.34}
{\rm Li}^{{\rm et}}_{n_1, ..., n_{m}}(x_1, ..., x_{m}) \in {\cal A}^{{\rm et}}_w(F)
\end{equation}

b) If $F$ is a number field,  the l-adic realization 
of (\ref{1.8.01.33}) is given by  (\ref{1.8.01.34}). 

c) The elements (\ref{1.8.01.34}) satisfy the 
 shuffle product formula (\ref{3.7.01.1132}).
\end{theorem}

{\bf Proof}.  If $F$ is a number field we can  
apply the l-adic realization functor to the motivic objects from theorem 
\ref{4.16.01.13sa1}. In particular this way we can settle the important case 
when $x_i$ are roots of unity. 

In  general we proceed just as in the Hodge case, since the proofs of all 
the results we used there work also in the l-adic situation. 
More specifically:

a) If the parameters are  admissible  the part a) is given 
by the constructions of chapters 2 and 3. In general 
we define (\ref{1.8.01.34}) by specialization, 
via the l-adic analog of the definition of $\widehat {\rm Li}$. 
The l-adic version of the specialization theorem \ref{10.24.01.1} guarantees that 
these two approaches give the same result in the case of admissible parameters. 

b) In the case of admissible parameters this is true  thanks to the results of 
chapters 2 and 3. For non admissible parameters we work 
out the l-adic version of corollary \ref{1.17.02.1}, which follows from 
the l-adic version of theorem \ref{1.13.02.5} and the l-adic version of 
lemma 6.7 in [G7], and use it just as we did in the Hodge case. 

c) Just as in the Hodge case, it is deduced from the fact that it is true for generic parameters 
(see the part d) of theorem \ref{4.17.01.77qa}, which is 
the l-adic version of theorem \ref{4.16.01.13saq}) and specialization theorem \ref{10.24.01.1}. 
The theorem is proved.

\section{How to prove identities between periods}


Below we discuss two general approaches to prove identities 
between the framed objects. Then we apply each of them to prove  theorem \ref{4.17.01.77z} 
in  the  double logarithm case

The first   approach suggested below allows 
to reduce a proof of any functional equation between 
the framed mixed Tate objects, e.g. Hodge-Tate structures, 
to a routine calculation. This  method 
is of  
algebraic nature, it proceeds by induction, and relies 
 only on the formula for the coproduct of the corresponding 
framed objects,
plus the fact that the identity is true at one particular 
point. 

The second  is the well known  direct approach to the problem. 
Any relation between  the framed mixed Tate structures of geometric origin 
is supposed to have a proof of this kind. However  such a proof does require
 an inspiring guess. So, unlike the first method, it 
is not an algorithmic procedure.

Both methods work equally well in the Hodge or l-adic setting.

{\bf 1. The ${\rm Li}$-shuffle relations for the double logarithm mixed Tate objects}. 
For 
$
xy \not = 0;  y \not = 1
$ we follow [G7] and define the double logarithm framed 
Hodge-Tate structures by 
$${\rm Li}^{\cal H}_{1,1}(x,y) := 
{\rm I}^{\cal H}_{1,1}((xy)^{-1}, y^{-1})
$$ 
They are fibers of  a unipotent variation of framed Hodge-Tate structures 
over a bit smaller domain (see [G10] or [G7]):
 \begin{equation} \label{4.1.02.12} 
\{(x,y) \subset {\Bbb C}^* \times {\Bbb C}^*| x \not =1; y \not =1; xy \not = 1\}
\end{equation}

\begin{theorem} \label{1.6.02.1} a) For any $(x,y)$ from the domain (\ref{4.1.02.12})  one has 
\begin{equation} \label{4.1.02.1tt} 
 {\rm Li}^{\cal H}_{1}(x)
 {\rm Li}^{\cal H}_{1}(y) - ( {\rm Li}^{\cal H}_{1,1}(x,y) + 
 {\rm Li}^{\cal H}_{1,1}(y,x) +  {\rm Li}^{\cal H}_{2}(xy))  = 0
\end{equation} 

b) The l-adic version of this identity holds. 
\end{theorem}

{\bf 2.  Two methods  to prove identities between the periods}.  
The first requires  to introduce parameters, i.e.  
interpret the periods as special values of 
period functions arising from  (unipotent) 
variations of mixed Hodge structures of geometric nature (e.g. multiple polylogarithms). 
Then we proceed as  follows:  

1). {\it Prove that the differential of the suspected identity is zero. 
So it is 
valid up to a constant. Specializing to a (degenerate) point check that it is zero}. 

This method always works, but sometimes requires a lot of routine labor. 
The reason  it is so efficient is provided by the following basic fact: 
$$
d(\mbox{weight $w$ unipotent period 
function }) = 
$$
$$
\sum (\mbox{weight $w-1$ unipotent period 
functions}) \cdot d\log f_i
$$ 
where $f_i \in \C(X)^*$ are some rational functions. So taking a basis in $\C(X)^*$ 
and decomposing $f_i$ in this basis 
we find out that one needs to prove some identities for the weight $w-1$ period functions. 
So we proceed by induction on the weight.

{\bf Example}. To prove formula (\ref{2.1.02.4}) we use 
$$
d {\rm Li}_{1,1}(x,y) = \log (1-xy) \cdot d\log\frac{x(1-y)}{1-x} + 
\log(1-y) \cdot d\log (1-x)
$$ 
Then an easy algebraic calculation shows that the differential 
of (\ref{2.1.02.4}) is zero. It remains to notice that formula (\ref{2.1.02.4}) is obvious 
for $x=y=0$. 

{\bf Remark}. Of course  formula (\ref{2.1.02.4}) is obvious by the power series expansion, 
but as far as we know this approach has no apparent Hodge/motivic incarnation. 

The second approach is this:

2). {\it Prove an identity between the periods by using  
identities of algebraic-geometric nature 
between 
between the corresponding differential forms and cycles, and  the Stokes formula}. 

This method does not require  to interpret the   periods 
as special values of some period functions. 

The Hodge/motivic version of  the first  method  is explained in section 8.3. 
The second method is manifestly motivic.

{\bf 3. How  to prove identities between the framed objects}. 
The motivic version of the first method 
 is based on the following rigidity lemma. We  spell it first the Hodge case. 
Recall that $\Delta$ is the coproduct in the Hopf algebra ${\cal H}_{\bullet}$ of the framed Hodge-Tate structures. 
Similarly one can consider the Hopf algebra of unipotent 
variations of framed Hodge-Tate structures
over a base $S$.

\begin{lemma}\label{1.4.02.2}
Suppose that $H_{S}$ is a unipotent variation of Hodge-Tate structures 
over a connected manifold $S$ framed 
by $\Q(0)$ and $\Q(n)$, where $n>1$. 
Then if $\Delta $ kills $H_{S}$  then 
$H_{S}$ is equivalent to the  constant variation over $S$ representing an element of 
${\rm Ext}_{MHS/S}^1(\Q(0), \Q(n)) $. 
\end{lemma}

In particular if $H(s)$ is the fiber at a point $s \in S$ 
then it represents an element of  ${\rm Ext}_{MHS}^1(\Q(0), \Q(n)) $ 
which does not depend on $s$. 

{\bf Proof}. Recall that  ${\rm Ker} \Delta \subset {\cal H}_n$ is identified with 
$$
{\rm Ext}_{MHS}^1(\Q(0), \Q(n)) = 
\frac{\C}{(2\pi i \Q)^n}\subset {\cal H}_n
$$
The Griffiths transversality condition implies that for $n>1$, and 
a connected smooth base $S$, 
one has 
$$
{\rm Ext}_{MHS}^1(\Q(0), \Q(n)) =  {\rm Ext}_{MHS/S}^1(\Q(0), \Q(n))
$$
The lemma is proved. 

{\it A general method to prove that 
an element $\sum_i H_i(s) \in {\cal H}_n$ is zero}:

1) Find a variation of framed Hodge-Tate 
structures over a smooth connected base $S$ whose fiber at $s \in S$ is equivalent to 
$\sum_i H_i(s)$.

2) Check that  $\Delta(\sum_i H_i(s)) =0$ for all $s \in S$.

3) Prove that it is zero at a single point $s_0 \in S$, or that 
specialization of the underlying variation 
at a certain tangential base point at infinity is zero.

Thanks to lemma \ref{1.4.02.2} the conditions 1) and 2) guarantee that our 
variation is constant, and  
3)  implies that it is zero. 

In the l-adic case there is a similar rigidity lemma provided by  (\ref{2.10.02.11}). 
Thus we can repeat the described  above scheme  in the l-adic case.

  {\bf 4. The first proof of theorem \ref{1.6.02.1}}. We apply the method above 
to  the left hand side of (\ref{4.1.02.1tt}).  
Recall the following formula for the coproduct, see proposition 2.3 in [G10] or example 1 in 
ch. 6 in [G7]:
$$
\Delta {\rm Li}^{\cal H}_{1,1}(x,y) = (1-xy) \otimes \frac{x(1-y)}{1-x} + 
(1-y) \otimes (1-x) \in   \C^*_{\Q} \otimes \C^*_{\Q} = {\cal H}_1^{\otimes 2}
$$
Here $\C^*_{\Q} =\C^*\otimes {\Q}$, and  the last isomorphism is provided by the canonical isomorphisms
 ${\cal H}_1 = {\rm Ext}_{MHS}^1(\Q(0), \Q(1)) = \C^*_{\Q}$. 
Therefore  
$$
\Delta( {\rm Li}^{\cal H}_{1,1}(x,y) + {\rm Li}^{\cal H}_{1,1}(y,x)) = 
(1-xy) \otimes xy  + (1-y) \otimes (1-x)+ (1-x) \otimes (1-y)
$$
and, since ${\rm Li}^{\cal H}_{1}(z) = -(1-z) \in \C^*_{\Q} = {\cal H}_1$,
$$
\Delta( -{\rm Li}^{\cal H}_{2}(xy) + {\rm Li}^{\cal H}_{1}(x){\rm Li}^{\cal H}_{1}(y)) = 
(1-xy) \otimes xy  + (1-y) \otimes (1-x)+ (1-x) \otimes (1-y)
$$
Thus  $\Delta$ kills the left hand side of (\ref{4.1.02.1tt}).  
Each of the objects from the left hand side of (\ref{4.1.02.1tt}) 
is a fiber of a certain 
unipotent variation of Hodge-Tate structures over the   space 
(\ref{4.1.02.12}),  
framed by $\Q(0)$ and $\Q(2)$. Taking the specialization of 
the left hand side of (\ref{4.1.02.1tt}) at $x =0$, and 
after this  at  $y=0$ we see that  each term in (\ref{4.1.02.1tt}) 
specializes to zero.  This implies the Hodge version of the theorem. 
Therefore we have the motivic version of this result when $x,y$ are from  a number field. 
The l-adic case now follows from lemma \ref{infr}. 
Theorem \ref{1.6.02.1} is proved.

{\bf 5. The second proof of theorem \ref{1.6.02.1}}. 
A peculiar property 
 of  formula (\ref{4.1.02.1tt}) is  this. Even for generic $x,y$, and even 
on the level of framed  
Hodge-Tate structures,  one of the terms in the formula, ${\rm Li}^{\cal H}_{2}(xy)$,
 was defined in [G7] 
not 
directly, but using the specialization (= canonical regularization), 
which is a rather sofisticated functor from the  algebraic geometric point of view. 
So the first thing we might want  to do is to use a model where a geometric definition 
of ${\rm Li}^{\cal H}_{2}(xy)$ is available. Here is how it works.

Let us blow up the point $(0,0)$ at the $(t_1, t_2)$ plane. The 
natural coordinates on the blow up are $(u_1, u_2)$ such that 
$t_1 = u_1u_2$ and $t_2 = u_2$. The iterated integral 
for ${\rm Li}_{1,1}(x,y)$ on the blow up looks as follows. 
$$
{\rm Li}_{1,1}(x,y) = \int_{0 \leq t_1 \leq t_2 \leq 1} \frac{dt_1}
{t_1 - (x_1x_2)^{-1}}\wedge \frac{dt_2}{t_2 - x_2^{-1}} = 
$$
$$
\int_{0 \leq u_i \leq 1} 
\frac{d(u_1u_2)}
{u_1u_2 - (x_1x_2)^{-1}}\wedge \frac{du_2}
{u_2 - x_2^{-1}}  
$$
Introduce new variables $v_i = x_iu_i$  we arrive to the identity
$$
{\rm Li}_{1,1}(x,y) = \int_{0 \leq v_i \leq x_i} 
\frac{d(v_1v_2)}
{v_1v_2 - 1}\wedge \frac{dv_2}
{v_2 - 1}  
$$
Similarly 
$$
{\rm Li}_{1,1}(y,x) = \int_{0 \leq v_i \leq x_i} 
\frac{d(v_1v_2)}
{v_1v_2 - 1}\wedge \frac{dv_1}
{v_1 - 1} 
$$
and 
$$
{\rm Li}_{2}(xy) =  \int_{0 \leq v_i \leq x_i} 
\frac{d(v_1v_2)}
{v_1v_2 - 1}\wedge \frac{dv_2}
{v_2}  =  \int_{0 \leq v_i \leq x_i} 
\frac{dv_1\wedge  dv_2}
{v_1v_2 - 1} 
$$

On  the other hand 
$$
{\rm Li}_{1}(x ) {\rm Li}_{1}(y)=   \int_{0 \leq v_i \leq x_i} 
\frac{dv_1}
{v_1  - 1} \wedge \frac{dv_2}
{v_2  - 1}
$$
One has an equality of rational differential forms 
\begin{equation} \label{1.11.02.100}
\frac{dv_1}
{1 - v_1} \wedge \frac{dv_2}
{1 - v_2} \quad = \quad \frac{d(v_1v_2)}
{ 1 - v_1v_2 }\wedge \frac{dv_2}
{1  - v_2} - \frac{d(v_1v_2)}
{1 - v_1v_2 }\wedge \frac{dv_1}
{1 - v_1}  + \frac{dv_1\wedge  dv_2}
{1 - v_1v_2}
\end{equation}
To check it one can  develop
 the denominators in the left and right hand sides 
into the geometric progression 
$$
\sum_{0 < k_1, k_2}  v_1^{k_1}v_2^{k_2} \cdot dv_1 \wedge dv_2 = 
\Bigl(\sum_{0 < k_1 < k_2} + \sum_{k_1 > k_2 >0} + 
\sum_{k_1 = k_2 >0}\Bigr) v_1^{k_1}v_2^{k_2} \cdot dv_1 \wedge dv_2
$$

Let us interpret 
 each of the integrals above as periods of appropriate framed Hodge-Tate structures.

Lifting the real triangle $ 0 \leq t_1\leq t_2  \leq 1$ on the blow up 
we get a square  $0 \leq u_i \leq 1$ denoted $C_2$. 
It is just the same as the rectangle $0 \leq v_i \leq x_i$. 

Let ${\Bbb U}_{(2)}$ be the affine chart on the blow up 
where the coordinates $(u_1, u_2)$ are defined. The square $C_2$  lies inside of 
${\Bbb U}_{(2)}(\R)$. 
Each of the four double integrals above has the form 
\begin{equation} \label{1.12.02.1}
\int_{C_2}\Omega
\end{equation} 
 where 
$\Omega$ is one of the 2-forms from (\ref{1.11.02.100}). 

Denote by $C^*_{(2)}$ 
the restriction to ${\Bbb U}_{(2)}$ of 
the preimage of the algebraic triangle 
$
B:= \{0=t_1\}\cup \{t_1=t_2\}\cup \{t_2=1\}
$. 
Thus the square provides a generator  $$
[C_{2}] \in 
H_2({\Bbb U}_{(2)}(\C), C^*_{(2)}(\C))
$$
Let ${\rm Sing}(\Omega)$ be the singularity divisor of $\Omega$.
Then if $y \not = 1$ and $\Omega$ is one of the forms in (\ref{1.11.02.100}) 
then ${\rm Sing}(\Omega) \cup C^*_{(2)}$ is a normal crossing divisor on ${\Bbb U}_2$. 

{\bf Remark}. For generic $x,y$ a similar pair of divisors on the $(t_1, t_2)$ plane 
is  a normal crossing divisor for ${\rm Li}_{1,1}(x,y)$, but not for ${\rm Li}_{2}(xy)$. 
We had to make the blow up to get a normal crossing divisor for ${\rm Li}_{2}(xy)$.

So setting 
$$
C^*_{(2), \Omega}:= C^*_{(2)} - C^*_{(2)}\cap {\rm Sing}(\Omega)
$$
 we see that 
 integral (\ref{1.12.02.1})  is a period of the 
framed Hodge-Tate structure of geometric origin 
$$
\Bigl( H^2({\Bbb U}_{(2)} - {\rm Sing}(\Omega), C^*_{(2), \Omega}), [C_{2}], [\Omega]\Bigr)
$$
It is equivalent (and in fact isomorphic) to the corresponding term in \ref{1.6.02.1}. 
Replacing  ${\rm Sing}(\Omega)$ in 
$H^2({\Bbb U}_{(2)} - {\rm Sing}(\Omega), C^*_{(2), \Omega})$ by the  union of such divisors for 
all four 2-forms appearing in (\ref{1.11.02.100}), and cutting down 
$ C^*_{(2), \Omega}$ accordingly,  we get a framed object equivalent 
to the one above. Applying to these objects identity 
(\ref{1.11.02.100}) and lemma \ref{4.17.01.77} we get theorem \ref{1.6.02.1} in the 
Hodge setting, and hence, as explained in the end of s. 8.3, in the motivic and 
finally l-adic settings.

{\bf 6. Calculation of  ${\rm Li}^{\cal H}_{1,1}(x,1)$}. 
By proposition \ref{1.7.02.1} for $x \not = 1$ we have 
$$
{\rm Li}^{\cal H}_{1,1}(x,1):= 
{\rm Sp}_{\varepsilon \to 0}{\rm Li}^{\cal H}_{1,1}(x, 1 -\varepsilon)
$$
If $x \not = 1$ applying the ${\rm Sp}_{\varepsilon \to 0}$ functor 
to 
$$
 {\rm Li}^{\cal H}_{1}(1 - \varepsilon)
 {\rm Li}^{\cal H}_{1}(x) = 
 {\rm Li}^{\cal H}_{1,1}(1 - \varepsilon, x) 
+ {\rm Li}^{\cal H}_{1,1}(x, 1 - \varepsilon) 
+  {\rm Li}^{\cal H}_{2}((1 - \varepsilon)x)
$$
we get
$$
\widetilde {\rm Li}^{\cal H}_{1,1}(x, 1) = - 
{\rm Li}^{\cal H}_{1,1}(1, x) + {\rm Li}^{\cal H}_{2}( x) 
$$
This identity (as we already know) remains valid fop $x=1$.

\section{A  proof of theorem \ref{4.17.01.77z}} 

In this chapter we give an algebraic-geometric proof  of the shuffle product formula 
from theorem \ref{4.17.01.77z}. The subtlety of this 
 problem  is explained by the following. 
Even for generic $x_i$, and even 
on the 
Hodge level,  the depth $< p+q$ terms  in  (\ref{3.7.01.1132pp}) 
 were defined in [G7] by  using the specialization functor,  
which is a rather sofisticated functor from the  motivic point of view. 
So we wanted to find  a model where a direct geometric definition 
of all terms is available. Below we use a  sequence of blow ups  which transforms 
the algebraic simplex to the algebraic cube. 
After this we use arguments similar to the one used in chapter 5.

{\bf 1.  The set up}. 
Let  $t$ be the canonical coordinate on $X:= {\Bbb A}^1$, and  
and $t_1, ..., t_n$ the  coordinates 
on $X^n$. Recall the standard algebraic simplex $B$ in $X^n$:
\begin{equation} \label{4.16.01.21}
B:= \quad \{0 = t_1\} \cup \{t_1 = t_2\}\cup  \{t_2 = t_3\} \cup  ... \cup  
\{t_{n-1} = t_n\}\cup   \{t_n = 1\}
\end{equation}
Let us use the flag of subvarieties 
$$
\{0 = t_1 = ... = t_n\} \subset \{0 = t_1 = ... = t_{n-1}\}
\subset ... \subset \{0 = t_1 = t_2\}
$$
of dimensions $0, 1, ..., n-2$ to construct a sequence of blow ups of 
$X^n$. Namely, blow up the zero dimensional subvariety on $X^n$, then blow up 
the strict preimage of the one  dimensional subvariety of the flag,
 then blow up  the strict preimage of the two dimensional subvariety of the flag,  
and so on. On the final stage we get a variety denoted $Y_{(n)}$. 
The preimage of the divisor $B$ in $Y_{(n)}$
is denoted $C_{(n)}$. 

\begin{lemma} The divisor $C_{(n)}$ has a 
shape of the $n$-dimensional cube. 
\end{lemma}
  
{\bf Proof}. After the blow up of the point $\{0 = t_1 = ... = t_n\}$ 
the preimage of $D$ has the shape of $\Delta^1 \times \Delta^{n-1}$. 
Then one proceeds by induction. 

Another way to see the cube is this. By the very definition of the blow up 
the divisor 
$C_{(n)}$  has a shape of the polytop obtained as follows. 
Take the standard simplex 
\begin{equation} \label{4.7.01.1}
\Delta^{n} = \{0 \leq t_1 \leq ... \leq t_n \leq 1\} \subset \R^n
\end{equation} 
and cut out 
the vertex 
$(0, ..., 0)$ by a hyperplane $t_n = \varepsilon$, $\varepsilon > 0$, then  
(the remaining of) the edge $(0, ..., 0, t)$ by $t_{n-1} = \varepsilon$, 
then cut out the 2-face 
$(0, ..., 0, t, t)$ and so on. The polytop we get is a cube. 
The lemma is proved. 

There is a  natural coordinate system $(u_1, ..., u_n)$ on $Y_{(n)}$ 
 such that 
\begin{equation} \label{1.2.3.4.5.}
t_1 = u_1 ... u_n; \quad t_2 = u_2 ... u_n; \quad ... \quad; 
t_{n-1} = u_{n-1} u_n; 
\quad t_n = u_n
\end{equation}
so that 
$$
u_1 = \frac{t_1}{t_2}, \quad u_2 = \frac{t_2}{t_3}, \quad ...\quad  
u_{n-1} = \frac{t_{n-1}}{t_n}, \quad u_n = t_n
$$

Then the standard simplex (\ref{4.7.01.1}) is transformed to the standard cube
$$
C_n = \{(u_1, ..., u_n) \subset \R^n \quad | \quad 0 \leq u_i \leq 1\}
$$

Let $w := n_1 + ... + n_m$. 
Let us write   
 integral representation (\ref{5*}) on the blow up $Y_{(w)}$ using  the variables $u_i$. Then   introduce another variables  $v_j$:
$$
v_j:= \left\{ \begin{array}{ll}
u_j  & \mbox{if $j \not = n_1 + ... + n_k $ \quad for some $k$}
\\ 
x_p u_{n_1+ ... + n_k} & \mbox{if $j = n_1 + ... + n_k$ \quad 
for some $k$}\end{array} \right. 
$$ 
The cube $C_w$ looks in this coordinates as the cube 
\begin{equation} \label{7.11.01.27}
C_{n_1, ..., n_m}(x_1, ..., x_m):= 
\end{equation}
$$
\{ (v_1, ..., v_w) \quad | \quad 0 \leq 
v_{n_1 + ... + n_p} \leq x_p, \quad 0 \leq v_j \leq 1 \quad \mbox{otherwise}\}
$$

Let us split the sequence of variables $v_1, ..., v_w$ into $m$ consecutive segments of lengths $n_1, ..., n_m$:
\begin{equation} \label{7.11.01.7}
S_1:= \underbrace {v_1, ..., v_{n_1}}_{\mbox{$n_1$  {\rm terms}}}; \quad 
S_2:= \underbrace {v_{n_1+1}, ..., v_{n_1+n_2}}_{\mbox{$n_2$  {\rm terms}}}; \quad ... \quad ;
S_w:= \underbrace {v_{w- n_{m}+1}, ..., v_w}_{\mbox{$n_m$  {\rm terms}}} 
\end{equation}
 and, using this splitting, define the following differential form 
on $Y_{(w)}$:
$$
\Omega_{n_1, ..., n_{m}}(x_1, ..., x_m):= 
$$
\begin{equation} \label{1.2.3.4}
\underbrace {\frac{d(v_1 ... v_w)}{1-v_1 ... v_w}\wedge 
\frac{d(v_2 ... v_w)}{ v_2 ... v_w}\wedge ...\wedge  
\frac{d(v_{n_1} ... v_w)}{ v_{n_1} ... v_w}}_{\mbox{$n_1$ {\rm factors}}} \wedge ... \wedge 
\end{equation}
$$
  \underbrace {\frac{d(v_{w-n_{m} + 1}  \cdot ... \cdot  v_w)}{1-v_{w-n_{m} + 1}\cdot  ... \cdot v_w}\wedge  ...\wedge \frac{d(v_{w-1} v_w)}{v_{w-1} v_w}  \wedge
\frac{dv_w}{ v_w}}_{\mbox{$n_m$  {\rm factors}}} $$
For example
$$
\Omega_{2,1}(x_1, x_2):= \frac{d(v_1  v_2 v_3)}{1-v_1  v_2 v_3}\wedge 
\frac{d( v_2 v_3)}{ v_2 v_3}\wedge \frac{d v_3}{ v_3}
$$
 Let us stress that it indeed depends on the parameters 
$x_i$ since the variables $v_j$ were defined using $x_j$'s and  the 
natural coordinates 
$u_j$ on $Y_{(w)}$. On the other hand the cube (\ref{7.11.01.27}) on $Y_{(w)}$ 
does not depend on $x_i$'s. 

To clarify the algebraic nature of this differential form 
recall the map
\begin{equation} \label{1.10.02.1}
d\log^{\wedge n}: \Lambda^n F(Y)^* \lra \Omega^n_{\log}(F(Y))
\end{equation}
$$ f_1 
\wedge ... \wedge f_n \lms d\log (f_1) \wedge ... \wedge d\log (f_n)
$$
The form (\ref{1.2.3.4}) is obtained by applying such a map 
to an element 
$$
F_{n_1, ..., n_m}(x_1, ..., x_m) \in \Lambda^w F(Y_{(w)})^*; \qquad F = 
\Q(x_1, ..., x_m)
$$
Precisely, $F_{n_1, ..., n_m}(x_1, ..., x_m) = (-1)^wf_1\wedge ... \wedge f_w$ where 

\begin{equation} \label{4.16.01.1}
f_j:= \left\{ \begin{array}{ll}
\prod_{i \geq j}v_i  & \mbox{if $j \not = n_1 + ... + n_k +1$ 
\quad for some $k$}
\\ 
1-\prod_{i \geq j}v_i & \mbox{if $j = n_1 + ... + n_k +1$ \quad 
for some $k$}\end{array}\right.
\end{equation}

 For example
$$
F_{2,1}(x_1, x_2) = - (1-v_1v_2v_3) \wedge (v_2v_3) \wedge v_3
$$

\begin{lemma} \label{4.9.01.1} One has 
$$
{\rm Li}_{n_1, ..., n_m}(x_1, ..., x_m) = 
\int_{C_w}  \Omega_{n_1, ..., n_{m}}(x_1, ..., x_m)
$$
 \end{lemma}

 {\bf Proof}. Here is how it works in the simplest cases. 
For ${\rm Li}_{1,1}(x,y) $ see ch. 8.

{\bf Example 1}. 
In the case of  ${\rm Li}_{2,1}(x,y) $ the computation looks as follows: 
$$
{\rm Li}_{2,1}(x,y) = (-1)^2\int\frac{d (u_1u_2u_3)}{u_1u_2u_3 - (xy)^{-1}} 
\wedge \frac{d (u_2u_3)}{u_2u_3 } \wedge 
\frac{d u_3}{u_3 - y^{-1}}; \qquad 0 \leq u_i \leq 1
$$
Changing the variables $v_1 = x u_1, v_2 = u_2, v_3 = y u_3$ we get 
$$
\int_{C_{2,1}(x,y)}\frac{d (v_1v_2v_3)}{1-v_1v_2v_3} \wedge 
\frac{d (v_2v_3)}{v_2v_3}\wedge \frac{d v_3}{1-v_3}  \quad = \quad 
\int_{C_{2,1}(x,y)}\frac{v_3d v_1 \wedge dv_2 \wedge dv_3}
{(1-v_1v_2v_3)(1-v_3)} 
$$
 
{\bf Example 2}. For ${\rm Li}_{3}(x) $ one has:
$$
{\rm Li}_{3}(x) =- \int\frac{d (u_1u_2u_3)}{u_1u_2u_3 - x^{-1}} 
\wedge \frac{d (u_2u_3)}{u_2u_3} \wedge 
\frac{d u_3}{u_3}; \qquad 0 \leq u_i \leq 1
$$
Changing the variables $v_1 = u_1 x, v_2 = u_2, v_3 = u_3$ we get 
$$
\int_{C_{3}(x )} \frac{d (v_1v_2v_3)}{1-v_1v_2v_3} 
\wedge \frac{d (v_2v_3)}{v_2v_3} \wedge 
\frac{d v_3}{v_3}  \quad = \quad \int_{C_{3}(x )}\frac{d v_1 \wedge dv_2 \wedge dv_3}{1- v_1v_2v_3}
$$

In general the proof follows the same pattern. The lemma is proved.

{\bf 2. A framed mixed Tate motive corresponding 
to multiple polylogarithms 
with general arguments}. Let ${\Bbb U}_{(w)}$ be the affine chart  of ${Y}_{(w)}$ were 
the coordinates $u_1, ..., u_w$ are defined. So there is an isomorphism 
$(u_1, ..., u_w): {\Bbb U}_{(w)} \stackrel{=}{\lra} {\Bbb A}^w$. 
Observe that  $C_w \subset {\Bbb U}_{(w)}(\R)$.

Let $\Omega$ be a differential  $w$-form  on  ${\Bbb U}_{w}$ with logarithmic singularities 
at the divisor  ${\rm Sing}(\Omega) \subset {\Bbb U}_{w}$. 
We will assume that  $\Omega$ is in the image $d\log^{\wedge n}$ map, 
see (\ref{1.10.02.1}). 
In particular it is of weight $2w$. 
Let us suppose in addition that 
\begin{equation} \label{4.11.01.1}
\mbox{${\rm Sing}(\Omega)$ does not contain any 
face of the algebraic cube $C_{(w)}$}
\end{equation}
We set 
$$
C^*_{(w)}:=  C_{(w)} \cap {\Bbb U}_{w}; \quad  C^*_{(w), \Omega}:=  
C^*_{(w)} - (C^*_{(w)} \cap {\rm Sing}(\Omega))
$$
Then there is a mixed motive
\begin{equation} \label{1.11.02.13}
H^w({\Bbb U}_{(w)} - {\rm Sing}(\Omega), C^*_{(w), \Omega} )
\end{equation}
as well as its Hodge and l-adic realizations.

\begin{lemma} \label{1.11.02.10}
Let us assume (\ref{4.11.01.3}).  Then the  form 
\begin{equation} \label{4.16.01.14256}
\Omega:= \Omega_{n_1, ..., n_m}(x_1, ..., x_m)
\end{equation}
satisfies condition (\ref{4.11.01.1}). 
\end{lemma}

{\bf Proof}. It is sufficient to prove that 
${\rm Sing}(\Omega)$ does not contain any vertex of the cube, which is clear from a 
computation in the coordinates $u_i$. 

Thanks to lemma \ref{1.11.02.10} the object (\ref{1.11.02.13}) is canonically isomorphic to 
$H^0({\Bbb U}_{(w)}, {\cal F}^*_{A,B})$ for $A:= {\rm Sing}(\Omega)$ and $B:= C_{(w)}$. 
If $x_i$ are elements of a number field satisfying (\ref{4.11.01.3})  then 
(\ref{1.11.02.13}) 
is defined as a mixed Tate motive. 

The constructions of  chapter 2 provide a natural frame on 
the object (\ref{1.11.02.13}).  Its weight $2w$ component  
is provided by the cohomology class 
$[\Omega]$. The weight $0$ component  is provided by 
the  relative homology class  
\begin{equation} \label{4.16.01.1425}
[C_w] \in H_w( {\Bbb U}_{(w)},  
C^*_{(w)})
 \end{equation}
of the cube $C_w$. 
\begin{proposition} \label{chupa} Let us assume (\ref{4.11.01.3}). Suppose that 
 $\Omega$ is given by (\ref{4.16.01.14256}). Then  
there is an equivalence of framed Hodge-Tate structures 
$$
{\rm Li}^{\cal H}_{n_1, ..., n_m}(x_1, ..., x_m) \sim 
\Bigl(H^w({\Bbb U}_{(w)} - {\rm Sing}(\Omega),  C^*_{(w), \Omega}), [C_w], [\Omega] \Bigr)
$$
If $x_i$ are elements of  a number field then there is a similar equivalence of framed 
mixed Tate motives. If $x_i \in F$ and $\mu_{l^{\infty}} \not= F$ 
then there is a similar statement in the l-adic case. 
\end{proposition}

{\bf Proof}.  This is a very particular case of the situation considered in chapters 2 and 3. 
The proposition is proved. 

{\bf 3. Some algebraic identities between the differential forms}. 
Set 
\begin{equation} \label{4.11.01.1}
s_j := \prod_{v_i \in S_j}v_i = \prod_{n_j \leq i\leq  n_{j+1}-1}v_i \quad = \quad x_j 
\prod_{n_j \leq i\leq  n_{j+1}-1}u_i
\end{equation}
\begin{lemma} One has an identity of formal power series
$$
\Omega_{n_1, ..., n_{m}}(x_1, ..., x_m):= \Bigl(\sum_{0 < k_1 < ... < k_m} s_1^{k_1} ... s_m^{k_m} \Bigr)\cdot dv_1 \wedge ... \wedge dv_w
$$
The right hand side is convergent if $x_i \in \C$ and $|x_i| <1$. 
\end{lemma} 

{\bf Proof}. One has 
$$
\Omega_{n_1, ..., n_m}(x_1, ..., x_m) = 
$$
$$
\frac{\prod_{j=2}^{m}s_j}{1- \prod_{j=1}^{m}s_j} \cdot 
\frac{\prod_{j=3}^{m}s_j}{1- \prod_{j=2}^{m}s_j} \cdot ... \cdot 
\frac{1}{1- s_j}\cdot dv_1 \wedge ... \wedge dv_w = 
$$
$$
\Bigl(\sum_{0 < k_1 < ... < k_m} s_1^{k_1} ... s_m^{k_m} \Bigr)\cdot dv_1 \wedge ... \wedge dv_w
$$
The  power series are 
convergent on the cube $C_w$ if $|x_i| <1$. The lemma is proved.

Let $n_1, ..., n_{p+q}$ be positive integers and 
$w = n_1 + ... + n_{p+q}$. 
For a 
generalized shuffle $\sigma \in \overline \Sigma_{p,q}$ 
consider the form with  formal power series coefficients 
\begin{equation} \label{4.11.01.12}
\Omega^{\sigma}_{n_1, ..., n_{p+q}}(x_1, ..., x_{p+q}):= 
\Bigl(\sum_{(k_1, ..., k_{p+q}) \in Z^{\sigma}_{++}} s_1^{k_1} ... s_{p+q}^{k_{p+q}} \Bigr)
\cdot dv_1 \wedge ... \wedge dv_w
\end{equation}
This power series are convergent on the cube $C_w$ for $|x_i| <1$. 

Let $w_1 = n_1+...+n_p$ and $w_2 = n_{p+1}+...+n_{p+q}$. Consider the natural projections $p_i: {\Bbb U}_{(w_1+w_2)} \lra {\Bbb U}_{(w_i)} $ given by  $$
p_1: (v_1, ..., v_{w_1+w_2}) \lra (v_1, ..., v_{w_1}); \quad 
p_2: (v_1, ..., v_{w_1+w_2}) \lra (v_{w_1+1}, ..., v_{w_1 +w_2})
$$   
Then 
\begin{equation} \label{4.11.01.112}
  p_1^*\Omega_{n_1, ..., n_{p}}(x_1, ..., x_p) = 
\Bigl(\sum_{0< k_1 <  ... < k_{p}} s_1^{k_1} ... s_{p}^{k_{p}} \Bigr)
\cdot dv_1 \wedge ... \wedge dv_{w_1}
\end{equation}
\begin{equation} \label{4.11.01.212} 
p_2^*\Omega_{n_{p+1}, ..., n_{p+q}}(x_{p+1}, ..., x_{p+q}) = 
\end{equation}
and 
$$
\Bigl(\sum_{0 < k_{p+1} <  ... < k_{p+q}} s_{p+1}^{k_{p+1}} ... s_{p+q}^{k_{p+q}} \Bigr)
\cdot dv_{w_1+1} \wedge ... \wedge dv_{w_1+w_2}
$$

\begin{lemma} \label{7.11.01.18}
$$
p_1^*\Omega_{n_1, ..., n_{p}}(x_1, ..., x_p)  \wedge 
p_2^*\Omega_{n_{p+1}, ..., n_{p+q}}(x_{p+1}, ..., x_{p+q})  = 
$$
$$
\sum_{\sigma \in \overline \Sigma_{p,q}}\Omega^{\sigma}_{n_1, ..., n_{p+q}}(x_1, ..., x_{p+q})
$$
\end{lemma}

{\bf Proof}. This is obviously true 
on the formal power series level thanks to decomposition (\ref{4.16.01.2}).  
The lemma is proved.

Let us define  forms (\ref{4.11.01.12}) 
as rational forms with logarithmic singularities on $Y_{(w)}$. 
Just as in the definition of the 
form (\ref{7.11.01.7}), given positive integers $n_1, ..., n_{p+q}$ 
we split the sequence of variables $v_1, ..., v_w$ into $p+q$ consecutive 
segments of lengths $n_1, ..., n_{p+q}$:
\begin{equation} \label{4.11.01.6}
S_1:= \underbrace {v_1, ..., v_{n_1}}_{\mbox{$n_1$  {\rm terms}}};  \quad ... \quad ;\quad 
S_p:= \underbrace {v_{n_1+...+ n_{p-1}+1}, ...,v_{n_1+...+ n_{p}} }_{\mbox{$n_{p}$  {\rm terms}}} 
\end{equation}
\begin{equation} \label{4.11.01.11}
S_{p+1}:= \underbrace {v_{n_1 + ... + n_p+1}, ..., v_{n_1+...+n_{p+1}}}_{\mbox{$n_{p+1}$  {\rm terms}}}; \quad ... \quad ;\quad 
S_{p+q}:= 
\underbrace {v_{w- n_{p+q}+1}, ..., v_w}_{\mbox{$n_{p+q}$  {\rm terms}}} 
\end{equation}
Given  a generalized shuffle $\sigma\in \overline \Sigma_{p,q}$ 
we define a new sequence of 
segments by shuffling, according to $\sigma$, the segments 
(\ref{4.11.01.6}) and (\ref{4.11.01.11}). If a segment $S_i$ from 
(\ref{4.11.01.6}) and 
a segment $S_j$ from  
(\ref{4.11.01.11}) correspond to the same slot of 
 $\sigma$ we put $S_iS_j$ as a new segment 
corresponding to this slot. 
The new sequence of segments determines a new ordered sequence 
$(v_1', ..., v_w')$ of $v_i$'s. Denote by 
$\widetilde \sigma $ the permutation  
$(v_1, ..., v_w) \lra (v'_1, ..., v'_w)$. 

Let us cook up the differential form corresponding to the 
new  sequence of 
segments following the rule (\ref{7.11.01.7}), and multiply it by 
$(-1)^{|\widetilde \sigma|}$.  It is easy to see that 
the form we get admits a power series decomposition identical with the right 
hand side of (\ref{4.11.01.12}). 
Thus 
$$
\Omega^{\sigma}_{n_1, ..., n_{p+q}}(x_1, ..., x_{p+q}) = 
d\log^{\wedge w}(F^{\sigma}_{n_1, ..., n_{p+q}}(x_1, ..., x_{p+q}))
$$
for an element 
$F^{\sigma}_{n_1, ..., n_{p+q}}(x_1, ..., x_{p+q}) = 
(-1)^{|\widetilde \sigma|} f^{\sigma}_1 \wedge ... \wedge
 f^{\sigma}_w$ defined via the rule (\ref{4.16.01.1}) applied to 
the new sequence of segments. 

\begin{lemma} \label{4.17.01.110}
The formula in  lemma \ref{7.11.01.18} remains valid for the corresponding 
rational differential forms.
\end{lemma}

{\bf Proof}. Clear. 

{\bf Remark}. One can prove a stronger statement: the 
following equality holds in the Milnor $K$-group $K^M_w(F(Y_{(w)}))$
$$
p_1^*F_{n_1, ..., n_{p}}(x_1, ..., x_p)  \wedge 
p_2^*F_{n_{p+1}, ..., n_{p+q}}(x_{p+1}, ..., x_{p+q})
 $$
$$
\sum_{\sigma \in \overline \Sigma_{p,q}} 
F^{\sigma}_{n_1, ..., n_{p+q}}(x_1, ..., x_{p+q}) = 
$$

 A path $\gamma$ be between $0$ and $1$ provides a relative cycle 
 modulo the algebraic simplex $B$. Let  
$\gamma(C_w)$ be the closure of its lifting   
on  $Y_{(w)}(\C)$. 

\begin{lemma} \label{4.17.01.41z} 
Define ${\rm Li}^{\sigma}_{n_1, ..., n_{p+q}}(x_1, ..., x_{p+q})$  
using an integral presentation like 
(\ref{1.16.02.1z}) 
written for an arbitrary path $\gamma$ between $0$ and $1$ instead of the path $[0,1]$. Then 
\begin{equation} \label{4.17.01.41}
\int_{\gamma(C_w)} \Omega^{\sigma}_{n_1, ..., n_{p+q}}(x_1, ..., x_{p+q}) = 
{\rm Li}^{\sigma}_{n_1, ..., n_{p+q}}(x_1, ..., x_{p+q})
\end{equation}
\end{lemma} 

{\bf Proof}. We work out an example when $p=q=1$. The general case is completely similar. 
Consider the shuffle $\sigma$ such that 
$$
{\rm Li}^{\sigma}_{n_1, n_{2}}(x_1, x_{2}) = 
{\rm Li}_{n_2, n_{1}}(x_2, x_{1})
$$
Then if $\Omega_{n_1, n_{2}}(x_1, x_{2})$ is defined using the sequence 
$$
S_1 = v_1, ..., v_{n_1}; \quad S_2 = v_{n_1+1}, ..., v_{n_1+n_2}
$$
the form  $\Omega^{\sigma}_{n_1, n_{2}}(x_1, x_{2}) $ is defined using the sequence
$$
S'_1 = v_{n_1+1}, ..., v_{n_1+n_2}; \quad S'_2 = v_1, ..., v_{n_1}
$$
So if we make change of variables
$$
v_1' = v_{n_1+1}, ..., v_{n_2}' = v_{n_1+n_2}; \qquad 
v_{n_2+1}' = v_{1}, ..., v_{n_1+n_2}' = v_{n_1}
$$
then
$$
\int_{C_w}\Omega^{\sigma}_{n_1, n_{2}}(x_1, x_{2}) = 
\int_{C_w'}\Omega_{n_2, n_{1}}(x_2, x_{1})
$$
where $C_w' = C_{n_2, n_{1}}(x_2, x_{1})$ is the corresponding cube in the 
$v'_i$ space.  

The lemma is proved. 

 Observe that 
\begin{equation} \label{4.18.01.1}
p_1^*C_{n_1, ..., n_p}(x_1, ..., x_p) \times 
p_2^*C_{n_{p+1}, ..., n_{p+q}}(x_{p+1}, ..., x_{p+q}) = 
\end{equation}
$$ C_{n_1, ..., n_{p+q}}(x_1, ..., x_{p+q})
$$ 
It follows from (\ref{4.17.01.41}) that  integrating the left and right parts of the formula in  lemma 
\ref{7.11.01.18} over 
the cycles located in the left and right parts of the  
identity (\ref{4.18.01.1})
 we get the ${\rm Li}$-shuffle relation (\ref{4.16.01.13}) where all the terms are understood 
as iterated integrals along a path $\gamma$. A 
proof of the identity (\ref{4.16.01.13}) for the case of 
multiple $\zeta$-values using similar ideas  was independently found by P. Cartier (unpublished).

{\bf 4. Proof of theorem \ref{4.17.01.77z}}.
Let 
$x_1, ..., x_{p+q} \in \C^*$ are as in (\ref{4.11.01.3}). Denote by 
$$
{\rm Li}^{{\cal H}, \sigma}_{n_1, ..., n_{p+q}}
(x_1, ..., x_{p+q}) 
$$
 the equivalence class of the framed mixed Tate motive 
$$
 \Bigl({\rm Li}^{{\cal H}}_{n_1, ..., n_{p+q}}(x_1, ..., x_{p+q}), [C_w], 
[\Omega^{\sigma}_{n_1, ..., n_{p+q}}(x_1, ..., x_{p+q})] \Bigr)
$$ 
where the object itself and the $[C_w]$-component of the frame were defined 
in proposition \ref{chupa}.

\begin{theorem} \label{4.17.01.77qa} a) Suppose that $x_i \in \C^*$ satisfy condition 
(\ref{4.11.01.3}) for all terms of the formula (\ref{4.17.01.73}) below. Then  
\begin{equation} \label{4.17.01.73}
{\rm Li}^{{\cal H}}_{n_1, ..., n_{p}}(x_1, ..., x_{p}) 
 \cdot {\rm Li}^{{\cal H}}_{n_{p+1}, ..., n_{p+q}}
(x_{p+1}, ..., x_{p+q}) = 
\end{equation}
$$
\sum_{\sigma \in \overline \Sigma_{p,q}}
{\rm Li}^{{\cal H}, \sigma}_{n_1, ..., n_{p+q}}(x_1, ..., x_{p+q}) 
$$

b) For a generalized shuffle $\sigma \in \overline \Sigma_{p,q}$ one has 
$$
{\rm Li}^{{\cal H}, \sigma}_{n_1, ..., n_{p+q}}(x_1, ..., x_{p+q}) 
  = {\rm Li}^{{\cal H}}_{\sigma(n_1, ..., n_{p+q})}
(\sigma(x_1, ..., x_{p+q})) 
$$

c) Suppose that $F$ is a number field and $x_i \in F^*$. Then there are
motivic versions of the parts a) and b) where ${\rm Li}^{{\cal H}}$ is replaced by 
${\rm Li}^{{\cal M}}$ 

d) If $\mu_{l^{\infty}} \not \in F^*$ then there is the  l-adic versions of a) and b).
\end{theorem}

{\bf Proof}. a) 
Proposition \ref{chupa} allows to interpret 
the framed object  in the left hand side of (\ref{4.17.01.73}) as 
 the object (\ref{1.11.02.13})  (where $w = n_1 + ... + n_{p+q}$)
with the frame provided by the cycle staying on the left of (\ref{4.18.01.1}), 
and by the product of forms (\ref{4.11.01.112}) and (\ref{4.11.01.212}). 
Each of the framed objects appearing in the sum on the right hand side of (\ref{4.17.01.73}) 
has a similar interpretation as the same object  (\ref{1.11.02.13})  with the frame given 
by the cycle class on the right hand side of (\ref{4.18.01.1}) and the forms $\Omega^{\sigma}$. 
Thus the statement follows from (\ref{4.18.01.1}) and lemmas \ref{7.11.01.18} and 
\ref{4.17.01.77}.

The part b) is a Hodge version of (\ref{4.17.01.41}). To prove it 
one needs to do a change of variables 
similar to the one which was done in the example following 
(\ref{4.17.01.41}). 

c) Follows from a) and b) using the injectivity of regulators, 
or can be deduced the same way we did a) and b). 

d) Is completely similarly to a) and b). The  theorem  \ref{4.17.01.77qa} is proved. 

Theorem \ref{4.17.01.77z} follows immediately from the parts a) and b) 
of theorem \ref{4.17.01.77qa}.  

The l-adic and motivic versions of theorem \ref{4.17.01.77z} follow immediately from the parts d) and c).

\section{Proofs of the theorems from  Introduction}


{\bf 1. Multiple polylogarithms Hopf and Lie algebras}. 
Let $F$ be a number field and $G \subset F^*$ be a subgroup. Then 
$
{\cal Z}_{w}^{\cal M}(G)
$ 
is the $\Q$-vector subspace of ${\cal A}_w(F)$ generated 
by the weight $w$ objects 
$$
 {\rm Li}^{{\cal M}}_{n_1, ..., n_{m}}(x_1, ..., x_{m}) \in {\cal A}_w(F); \quad  x_i \in G
$$   and 
$
{\cal Z}_{\bullet}^{\cal M}(G):= \oplus_{w \geq 0}{\cal Z}_{w}^{\cal M}(G)
$.

\begin{theorem} \label{4.16.01.s} ${\cal Z}_{\bullet}^{\cal M}(G)$ is a graded Hopf algebra.
\end{theorem}

{\bf Proof}. A similar result in the Hodge setting has been proved 
in theorem 6.12 of  [G7]. Observe that we use the $ {\rm Li}$-objects, 
while we have used the ${\rm I}$-objects in [G7]. However thanks to theorem 
\ref{1.13.02.5} or corollary \ref{1.17.02.1} 
this does not make any difference.  After that  lemma 3.4 from [G7] allows 
to deduce the motivic version from the Hodge one. The theorem is proved.

Similarly for any subgroup $G \subset F^*$ for any field $F$ such that $\mu_{l^{\infty}} \not \in F^*$ 
there is an l-adic Hopf algebra ${\cal Z}_{\bullet}^{\rm et}(G)$. If $F$ is a number field it can be 
defined as the $l$-adic realization of the corresponding motivic Hopf algebra. 

{\it The  depth filtration}. It is  an increasing filtration  $F^{D}_{\bullet}$ 
on the vector space ${\cal Z}_{\bullet}^{\cal M}(G)$ indexed by non negative integers. 
   $F^{D}_{m}{\cal Z}_{\bullet}^{\cal M}(G) $ is generated by the elements ${\rm Li}_{n_1, ..., n_k}^{\cal M}(x_1, ..., x_k)$ 
with  $k \leq m$. 
In particular $F^{D}_{0}{\cal Z}_{\bullet}^{\cal M}(G) = \Q$. 

\begin{lemma} \label{1.20.02.6}
a) The  depth filtration is compatible with the Hopf algebra 
structure on ${\cal Z}_{\bullet}^{\cal M}(G)$. 

b) A similar result hold for ${\cal Z}_{\bullet}^{\rm et}(G)$. 
\end{lemma}

{\bf Proof}. a) The Hodge version of this result is theorem 6.12 from [G7]. The motivic version is deduced from it using lemma 3.4 in [G7]. 

b) If $F$ is a number field it is deduced directly from a). The general case 
follows from proposition \ref{2.10.02.2}. 
The lemma is proved. 

Just like in section 1.5 we define for an arbitrary subgroup $G \subset F^*$ the  Lie coalgebra 
${\cal C}_{\bullet}^{\cal M}(G)$ and the Lie  algebra ${\rm C}_{\bullet}^{\cal M}(G)$. 
We can view them as  Lie (co)algebras  
in the category of pure Tate motives, or 
as a graded Lie (co)algebra over $\Q$. The grading is provided by the weight. 

The depth filtration is compatible with the weight grading. Taking
 the associate graded for the depth filtration we get the corresponding bigraded Lie coalgebra 
and Lie algebra  over $\Q$:
$$
{\cal C}^{\cal M}_{\bullet, \bullet}(G):= {\rm Gr}^D\left({\cal C}_{\bullet}^{\cal M}(G)\right); \qquad 
{\rm C}^{\cal M}_{\bullet, \bullet}(G):= {\rm Gr}^D\left({\rm C}_{\bullet}^{\cal M}(G)\right)
$$
 Equivalently, we can think about   them  as of a Lie (co)algebras 
in the category of pure Tate motives  graded by the depth.

For a subgroup $G \subset \C^*$ the corresponding Hodge Lie 
algebras ${\rm C}_{\bullet}^{\cal H}(G)$ and 
${\rm C}_{\bullet, \bullet}^{\cal H}(G)$ have been defined in chapter 6  of [G7]. 
There are their l-adic analogs 
${\rm C}_{\bullet}^{\rm et}(G)$ and 
${\rm C}_{\bullet, \bullet}^{\rm et}(G)$.

{\bf 2 Proof of theorem \ref{4.16.01.qa} }. We need to show that if $x_i \in \mu_N$ 
then the framed mixed Tate motive 
(\ref{1.8.01.33}) lies in ${\cal A}_w(S_N)$. A priori we know that   it lies in 
${\cal A}_w(\Q(\zeta_N))$. We employ definition 3.4 of ${\cal A}_w(S_N)$ 
from [G7]. It is equivalent to the one provided by  [DG]. 
According to that definition an element  $x \in {\cal A}_w(\Q(\zeta_N))$ belongs 
to the subspace ${\cal A}_w(S_N)$ if and only if 
$$
\Delta'(x):= \Delta(x) - (x \otimes 1 + 1 \otimes x) \in \oplus_{0 < k < w}{\cal A}_{k}(S_N) \otimes {\cal A}_{w-k}(S_N)
$$
Then the statement of the theorem follows 
by induction. Indeed,  the isomorphism (\ref{1.20.02.5}) 
provides the base of induction, and  theorem \ref{4.16.01.s} the induction step. 
The theorem is proved.

{\it A scetch of an alternative proof of theorem \ref{4.16.01.qa}}. 
Recall that a  mixed Tate motive over $\Q$ is  
unramified outside $N$ if and only if 
its $l$-adic realization for a prime   $l$ such that   $(l,N) =1$  
is unramified outside $lN$. The elements of the Hopf algebra 
${\cal Z}^{\cal M}(\mu_N)$ are the matrix elements 
(in the sence of section 3.2 in [G7]) 
of the pro-object 
${\cal P}^{\cal M}({\Bbb G}_m - \mu_N; v_{0}, v_{1})$. So the claim follows from the fact that its l-adic realization 
${\cal P}^{(l)}({\Bbb G}_m - \mu_N; v_{0}, v_{1})$ is unramified outside of $lN$. 
More directly, we can use 
corollary \ref{2.24.02.6} and the geometric construction 
of the motivic torsor of path given in s. 6.3 to show that 
${\cal P}^{\cal M}({\Bbb G}_m - \mu_N; v_{0}, v_{1})$ is unramified 
over $S_N$.

{\bf 3. Mixed Tate categories and their Galois groups}. To prove theorem 
\ref{2.6.02.1} we need to recall below some basic  
material about the mixed Tate categories and their Galois groups, 
see section 3.1-3.2 in [G7] 
for more details. 

Let $({\cal M}, K(1))$ be a mixed Tate $K$--category, where $K$ is a characteristic zero 
field. 
It is equipped with canonical fiber functor
$$
\Psi: {\cal M} \lra {\rm Vect}_{\bullet}, \quad X \lra \oplus \Psi_n(X) = 
\oplus_{n} {\rm Hom}_{\cal M}(K(-n), {\rm Gr}^W_{2n}X)
$$
Forgetting the grading we get the fiber functor $\widetilde \Psi$. 
Let $X$ be an object of ${\cal M}$. A choice of non zero vectors 
$v^{(p)} \in \Psi_{2p}(X)$ and  $f^{(q)}\in \Psi_{2q}(X)$   provide a framed 
object  $(X, v^{(p)}, f^{(q)})$. 

\begin{definition} \label{420.02.6} ${\cal A}_{\bullet}(X)$ is the commutative 
$K$--algebra generated 
by the equivalence classes of framed objects
\begin{equation} \label{4.20.02.2}
[X, v^{(p)}, f^{(q)}] \in {\cal A}_{p-q}(\cal M); \quad  
\end{equation}
 for all pairs of non zero vectors 
$v^{(p)} \in \Psi_p(X)$ and  $f^{(q)}\in \Psi_q(X)$, for all $p \geq q$. 
\end{definition}

Let us choose for every integer $n$ 
 a basis 
in  $\Psi_n(X)$. It provides  the  dual basis 
in $\Psi_n(X)^{\vee}$. Obviously it is sufficient to consider in 
(\ref{4.20.02.2}) the framed objects corresponding to all basis/cobasis vectors only.

Recall that according to the Tannakian formalism the category ${\cal M}$ is 
canonically equivalent to the category of graded finite dimensional $K$--modules 
over a proalgebraic 
unipotent group scheme $U({\cal M}):= {\rm Aut}^{\otimes } \widetilde \Psi$ over $K$. 
The equivalence is given by the fiber functor $\Psi$. In particular 
the group scheme $U({\cal M})$ acts on the graded vector space $\Psi(X)$, 
and this action factorizes through a quotient $U^X$ of $U({\cal M})$, i.e. 
$U^X$ is the image of the group scheme $U({\cal M})$ acting on $\Psi(X)$. Let 
${\cal X}$ be the mixed Tate category  generated by the object $X$. Then 
$U^X = U(\cal X)$. 

\begin{lemma} \label{4.20.02.1} a) ${\cal A}_{\bullet}(X)$ is a graded Hopf algebra. 

b) It is canonically isomorphic to the Hopf 
algebra of regular functions on $U^X$, so ${\rm Spec}{\cal A}_{\bullet}(X) = U^X$.
\end{lemma}

{\bf Proof}. a) The graded $K$-vector space generated by the elements (\ref{4.20.02.2})
is closed under the coproduct by its very definition. Since ${\cal A}_{\bullet}(\cal M)$ 
is a Hopf algebra, the statement follows. 

b) Recall that the equivalence between the category ${\cal M}$ and the category of graded 
${\cal A}_{\bullet}(\cal M)$--comodules assignes to an object $X$ the graded comodule 
$\Psi(X)$ with the ${\cal A}_{\bullet}(\cal M)$--coaction 
$\Psi(X) \otimes \Psi(X^{\vee}) \lra {\cal A}_{\bullet}(\cal M)$ given by 
$$
v \otimes f \lms \quad \mbox{the equivalence class of the framed object $[M, v, f]$}
$$
This immediately implies the part b) of the lemma. The lemma is proved. 

{\bf 4. The motivic torsor of path on ${\Bbb G}_m - \mu_N$ and theorems
 \ref{2.6.02.1} and \ref{2.6.02.wqw1} }. 
Recall ([DG] or chapter 6 above)  the motivic torsor of path  
\begin{equation} \label{4.20.09}
{\cal P}^{\cal M}({\Bbb G}_m - \mu_N; v_0, v_1)
\end{equation}
  between the canonical 
tangential base points $v_0$ and $v_1$ at $0$ and $1$. It is   a proobject in 
${\cal M}_T(S_N)$.

\begin{theorem} \label{4.16.01.wq} There is canonical isomorphism
$$
{\cal A}_{\bullet}\left({\cal P}^{\cal M}({\Bbb G}_m - \mu_N; v_0, v_1)\right) = 
{\cal Z}_{\bullet}^{\cal M}(\mu_N)
$$
\end{theorem}

{\bf Proof}. Follows immediately from section 6.3 and 
lemma \ref{4.20.02.1}. The theorem is proved.

\begin{proposition} \label{4.20.2.7} 
For any $\varepsilon \in \{0,1,\infty\}$ there is an isomorphism 
$$
{\cal A}_{\bullet}\left({\cal P}^{\cal M}(X_N; v_0, v_1)\right) = 
{\cal A}_{\bullet}\left(\pi_1^{\cal M}(X_N; v_{\varepsilon})\right)
$$
\end{proposition}

{\bf Proof}.  Let us prove the proposition  for $\varepsilon = 0$. 
We will use a shorthand $X_N:= {\Bbb G}_m - \mu_N$. 

Applying the fiber functor $\Psi$ to the left action of 
$\pi_1^{\cal M}(X_N; v_{1})$  on 
${\cal P}^{\cal M}(X_N; v_0, v_1)$ we get  a map of graded 
$L_{\bullet}(S_N)$-modules 
$$
\Psi(\pi_1^{\cal M}(X_N; v_{0})) \otimes 
\Psi({\cal P}^{\cal M}(X_N; v_0, v_1)) 
\lra 
\Psi({\cal P}^{\cal M}(X_N; v_0, v_1))
$$ 
Let 
$L_{\bullet}(S_N)$ be the Lie algebra of the pro-group scheme $U({\cal M}_T(S_N))$
 and 
 $I_{\bullet}(X)$ its ideal annihilating 
$\Psi(X)$. A non zero vector $\gamma_0 \in \Psi_0({\cal P}^{\cal M})$
 provides an isomorphism
of $I_{\bullet}({\cal P}^{\cal M})$--modules $
\Psi(\pi_1^{\cal M}(X_N, v_0)) \otimes \gamma_0
 \lra \Psi({\cal P}^{\cal M}(X_N; v_0, v_1))
$. Thus 
$$
I_{\bullet}({\cal P}^{\cal M}(X_N; v_0, v_1)) \subset I_{\bullet}
(\pi_1^{\cal M}(X_N, v_0))
$$

Let us prove the opposite inclusion. Since ${\cal P}^{\cal M}
(X_N; v_0, v_1)$ is a right
$\pi_1^{\cal M}(X_N; v_{1})$--torsor,  a choice of an element 
$\gamma \in \Psi({\cal P}^{\cal M})$ provides a map 
from $L_{\bullet}(S_N)$ to $\pi_1^{\cal M}(X_N; v_1)$. 
Namely, if $l \in L_{\bullet}(S_N)$ and  $l(\gamma) = \gamma s_l$ for 
$s_l \in \pi_1^{\cal M}(X_N; v_{1})$, then $l \lms s_l$. 
Its restriction to the ideal 
$I_{\bullet}(\pi_1^{\cal M}(X_N, v_0))$ is independent on the choice of $\gamma$, and 
hence provide a   
graded Lie algebra homomorphism 
\begin{equation} \label{4.21.02.1} 
I_{\bullet}\left(\pi_1^{\cal M}(X_N; v_0)\right) \lra 
\Psi\left(\pi_1^{\cal M}(X_N; v_1)\right)
\end{equation}By its very definition it induces an injective map
$$
\frac{I_{\bullet}\left(\pi_1^{\cal M}(X_N; v_0)\right)}
{I_{\bullet}\left({\cal P}^{\cal M}(X_N; v_0)\right)} \hookrightarrow 
\Psi\left(\pi_1^{\cal M}(X_N; v_1)\right)
$$
Thus to prove the opposite inclusion we have to show that this map is zero. 

 Recall the canonical  
map, ``the motivic loop around $1$ based at $v_1$'':  
$$
\alpha^{\cal M}_1: \Q(1) \lra \pi_1^{\cal M}(X_N; v_1)
$$ 
Observe that the ideal 
 $I_{\bullet}(\pi_1^{\cal M}(X_N; v_0))$ kills  
\begin{equation} \label{4.21.02.2}
\gamma\circ \Psi(\alpha^{\cal M}_1(\Q(1))\circ\gamma^{-1} \subset
 \pi_1^{\cal M}(X_N; v_0)
\end{equation}
Since $L_{\bullet}(S_N)$ kills the loop $\alpha^{\cal M}_1$, using this 
it is easy to see that 
the image of (\ref{4.21.02.1}) centralizes $\Psi(\alpha^{\cal M}_1)$.
Since $\pi_1^{\cal M}(X_N; v_1)$ is a free Lie algebra, this 
implies that 
the image of (\ref{4.21.02.1}) is contained  in $\Psi(\alpha^{\cal M}_1(\Q(1)))$.  
Therefore it lies in the one dimensional subspace of weight $-2$. 
So our our claim is reduced to the following elementary
 fact. Let $W_{[p, q]}(X)$ is the subquotient 
of $X$ of weights $2p \leq w \leq 2q$). Then 
$$
{\cal A}_{-1}\left(W_{[-2, -1]}\left(\pi_1^{\cal M}(X_N; v_1)\right)\right) =
{\cal A}_{-1}\left(W_{[-1, 0]}\left({\cal P}^{\cal M}(X_N; v_0, v_1)\right)\right)
$$
Indeed, we know that the right hand side is genertaed as a $\Q$--vector space 
by the elements $\log^{\cal M}(1-\zeta_N^{\alpha})$, where $0< \alpha <N$, and 
 the left hand side is contained in the right hand side. So it remains 
to check that every element 
$\log^{\cal M}(1-\zeta_N^{\alpha})$ appears in the left hand side.  
This follows, for instance, from the results of [D], where the polylogarithmic 
quotient of the Hodge fundamental group  $\pi_1^{\cal H}(P^1 - \{0, 1, \infty\}; x)$ 
was explicitly described. In particular it was shown that 
$\log^{\cal H}(1-x)$ appears as a matrix element for an appropriate 
$(\Q(1), \Q(2))$--framing on $\pi_1^{\cal H}(P^1 - \{0, 1, \infty\}, x)$. 
Then the 
map $t \lms t/x$ provides an isomorphism  
$$
\pi_1^{\cal H}(P^1 - \{0, 1, \infty\}; x) = 
\pi_1^{\cal H}(P^1 - \{0, \zeta_N^{\alpha}, \infty\}; 1)
$$ 
The right object is a subobject of 
$\pi_1^{\cal H}(X_N; v_1)$. Thus 
setting  $x := \zeta_N^{-\alpha}$ we get the statement. 

To prove the  proposition for other $\varepsilon$ we proceed as follows. 
There is an isomorphism of $L_{\bullet}(S_N)$--modules (``reversing the 
path'')
$$
{\cal P}^{\cal M}(X_N; v_0, v_1) = {\cal P}^{\cal M}(X_N; v_1, v_0)
$$
So interchanging the role of $0$ and $1$ in the proof above we get the statement for 
$\varepsilon =1$.  
The involution $x \lms x^{-1}$ on 
$P^1$ provides an isomorphism of mixed Tate motives 
$$
\pi_1^{\cal M}(X_N; v_0) \stackrel{\sim}{\lra} \pi_1^{\cal M}(X_N; v_{\infty})
$$
Therefore the statement for  $\varepsilon = \infty$ follows from the one 
for $\varepsilon = 0$. The proposition is proved.


Recall the pro-algebraic group scheme ${\rm Spec}({\cal Z}_{\bullet}^{\cal M}
(\mu_N)$ over $\Q$. 
The grading on ${\cal Z}_{\bullet}^{\cal M}(\mu_N))$ provides an action of 
the group ${\Bbb G}_m$ 
on this group. According to lemma \ref{4.20.02.1}, theorem \ref{4.16.01.wq} just means that the  motivic Galois group  acts on  (\ref{4.20.09}) via its quotient given 
by the semidirect product of 
${\Bbb G}_m $ and ${\rm Spec}({\cal Z}_{\bullet}^{\cal M}(\mu_N))$.

{\bf Proof of theorem \ref{2.6.02.wqw1}}. 
Follows immeduately from theorem \ref{4.16.01.wq} and proposition \ref{4.20.2.7}.

{\bf Remark}. Below we will use only the l-adic version of  proposition \ref{4.20.2.7}. 

{\bf 5. Proof of theorem \ref{5.2.02.2}}.  Let $X$ be a mixed Tate motive over a number field $F$. Denote by  $L_{\bullet}(X)$ the image of the motivic 
Lie algebra acting on  $\Psi(X)$. Let $X_{\Q_l}$ 
be the l-adic realization of the motive $X$. Denote by ${\cal G}^{(l)}_{X}$ 
the Lie algebra of the image of ${\rm Gal}(\overline F/F(\zeta_{l^{\infty}}))$ 
acting on $X_{\Q_l}$. 

\begin{lemma} \label{4.22.02.q1}
There is an isomorphism 
$
{\cal G}^{(l)}_{X} = L_{\bullet}(X)\otimes \Q_l
$, 
and  canonical isomorphism of graded Lie algebras
$$
{\rm Gr}^W_{\bullet}{\cal G}^{(l)}_{X} = L_{\bullet}(X)\otimes \Q_l
$$
\end{lemma}

{\bf Proof}. Follows from the content of section 3.7 in [G4], see especially the 
second half of page 424 loc. cit..

\begin{corollary} \label{4.22.02.11}
There is canonical isomorphism of the graded l-adic pro-Lie algebras
$$
{\rm Gr}_{\bullet}^W{\cal G}^{(l)}_N = L_{\bullet}(\pi_1^{\cal M}({\Bbb G}_m - \mu_N; 
v_{\infty})\otimes \Q_l
$$
\end{corollary}

{\bf Proof}.   Follows immediately by going to the projective limit 
from lemma \ref{4.22.02.q1} and proposition \ref{4.20.2.7}.

The part a) of theorem \ref{5.2.02.2} follows from corollary \ref{4.22.02.11} and theorem \ref{4.16.01.wq}. The part b) is also straitforward. 

{\bf 6. Proof of theorem \ref{4.16.01.pi}: the first part.} 
We prove the  Hodge case. By lemma 3.4 in [G7] it implies the motivic one. 
The l-adic version follows from the motivic one. So we study the map 
$$
\nu^{\cal H}_{\bullet, \bullet}(\mu_N): {\cal D}_{\bullet, \bullet}(\mu_N) \lra 
{\cal C}^{\cal H}_{\bullet, \bullet}(\mu_N)
$$ 
We have to check that this map sends the double shuffle and distribution relations to zero, 
and is compatible with the cobracket.

{\it The relations}. We use the list of the relations is given in section 4.1 of [G4]. 
Observe that we used in [G4] the homogeneous notations 
$$
{\rm I}_{n_1, ..., n_m}(g_1: ... :g_{m}:g_{m+1}) := 
{\rm I}_{n_1, ..., n_m}(g_{m+1}^{-1}g_1: ... :g_{m+1}^{-1}g_{m}) 
$$
So the relation (i) from [G4] is respected by the definition. 

The ${\rm I}$-shuffle relations, that is relation (63) in [G4], are valid by the first part of 
lemma 6.6 in [G7]. 

To get the ${\rm Li}$-shuffle relation, that is relation (62) in [G4], 
we notice that thanks to 
 theorem \ref{1.13.02.5} we have for any $x_i \in \C^*$ (with $a_i$ are defined by (\ref{1.16.02.1zsd})
$$
{\rm Li}_{n_1, ..., n_m}^{\cal H}(x_1, ..., x_m)  \sim  {\rm I}_{n_1, ..., n_m}^{\cal H}(a_1, ..., a_m)
$$ 
modulo  depth $<m$ terms and products. So 
using  theorem \ref{4.16.01.13saq} we get relation (62) in [G4]. 

The distribution relations  in [G4] are easily checked either using the general method of 
section 8.2, or directly using the realization described in section 9 and the specialization theorem. 
The distribution relations are proved in the l-adic setting in [G4], but the proof works in 
the Hodge  setting as well. This gives the third proof. 

The normalizing relation ${\rm I}_1(e:e)$ is true since ${\rm I}^{\cal H}(1) =0$.  
So we conclude that the map $\nu^{\cal H}_{\bullet, \bullet}(\mu_N)$ 
is a well defined map of the  bigraded vector spaces. 

 {\bf  7. The map $\nu^{\cal H}_{\bullet, \bullet}(\mu_N)$ 
commutes with the cobrackets}. The formula for 
the coproduct $\Delta$ in  the 
Hopf algebra ${\cal Z}^{\cal H}(\C^*)$ has been computed explicitly 
in  theorem 6.5 in [G7], and more specifically in proposition 6.8 loc. cit. 
The coproduct $\Delta$ induces a cobracket in the corresponding Lie coalgebra. We denote it 
 by  $\delta_{\cal H}$. 
On the other hand the cobracket $\delta$ in the dihedral Lie coalgebra was defined in 
section 4.4 of [G4]. Let us first  recall the definitions of $\Delta$ and $\delta$.

{\it The formula for $\Delta$}. For  convenience of the reader we 
recall some  definitions and results borrowing directly from chapter 6 in [G7]. 

We package the framed Hodge-Tate structures 
${\rm I}^{\cal H}_{ n_1, ..., n_m}(a_1, ... , a_m)$ 
 into the generating series
\begin{equation} \label{CO1e}
{\rm I}^{\cal H}( a_1: ... : a_m: a_{m+1}|  t_1, ... ,t_m) := 
\end{equation}
$$
\sum_{n_i \geq 1}{\rm I}^{\cal H}_{n_1, ..., n_m}(\frac{a_1}{a_{m+1}}, ... , \frac{a_m}{a_{m+1}}) 
t_1^{n_1-1} ... t_m^{n_m-1} 
 \in \quad {\cal H}_{\bullet}[[t_1, ..., t_m]] 
$$
Following lemma 6.7 in [G7], it is convenient to introduce 
the  homogeneous  variables $(t_0:t_1: ...: t_m)$ replacing  the variables $(t_1, ..., t_m)$:  
\begin{equation}  \label{3.7.01.19}
{\rm I}^{\cal H}(a_1: ...: a_m: a_{m+1}| t_0: t_1: ... :t_m) := 
\end{equation}
$$
a_{m+1}^{t_0}  {\rm I}^{\cal H}(a_1: ... : a_m: a_{m+1}| t_1-t_0, ..., t_m-t_0)
$$

By proposition 6.8 in  [G7] we have
$$
\Delta {\rm I}^{\cal H}(a_1: ... : a_{m+1}| t_0: ... :t_{m}) 
$$
$$
\sum {\rm I}^{\cal H}( a_{i_1}: ...: a_{i_{k}}: a_{m+1}| t_{j_0}: t_{j_1}: ... :
  t_{j_{k}}) \otimes
$$
\begin{equation}  \label{ur100100}
\prod_{\alpha =0}^k \Bigl(
(-1)^{j_{\alpha}-i_{\alpha}} {\rm I}^{\cal H}(a_{j_{\alpha}}: a_{j_{\alpha}-1}: ... : 
 a_{i_{\alpha}}| 
-t_{j_{\alpha}}: -t_{j_{\alpha}-1}: ... : - t_{i_{\alpha}}) \cdot 
\end{equation}
$$
   {\rm I}^{\cal H}(a_{j_{\alpha}+1}: ... : a_{i_{\alpha+1}-1}: a_{i_{\alpha+1}}|
t_{j_{\alpha}}: t_{j_{\alpha}+1}: ... : t_{i_{\alpha+1}-1} )\Bigr)
$$
where the sum is over all special marked colored segment, i.e. over all 
sequences $\{i_{\alpha}\}$ and $\{j_{\alpha}\}$ 
such that 
\begin{equation}  \label{ur10*}
i_{\alpha } \leq j_{\alpha } < i_{\alpha +1} \quad \mbox{ for any 
$0\leq \alpha \leq k$}, \qquad j_0 = i_0 = 0, i_{k+1} = i_{m+1} 
\end{equation}

{\it A  geometric  interpretation of the formula (\ref{ur100100})} . We picture  generating series 
(\ref{3.7.01.19}) by a colored segment where  $a_0 =0$:
\begin{center}
\hspace{4.0cm}
\epsffile{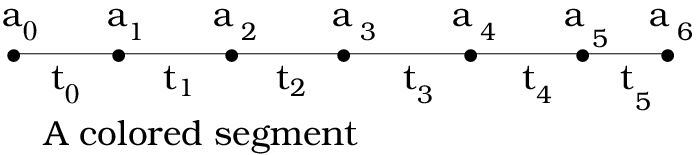}
\end{center}
The terms of this formula correspond  to the so-called  special marked colored segments: 

a) we mark (by making them boldface) the 
points $a_0; a_{i_1}, ..., a_{i_{k}}; 
a_{m+1}$.  

b) mark (by cross) segments $t_{j_0}, ..., t_{j_k}$ 
such that there is just one marked segment between any two 
neighboring marked points.

\begin{center}
\hspace{4.0cm}
\epsffile{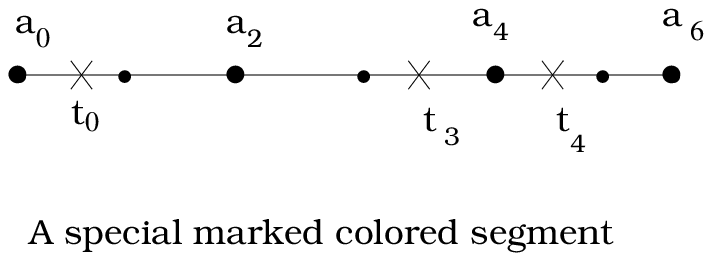}
\end{center}

{\it The cobracket 
$ 
\delta:     {\cal D}_{\bullet \bullet}(G)  \longrightarrow
  \Lambda^2{\cal D}_{\bullet \bullet}(G)$}. For the convenience of the reader we reproduce 
the definition given on the page 434 of 
[G4]. Consider the formal generating series similar to (\ref{CO1e})
\begin{equation} \label{ginf11}
\{g_1: ... : g_{m+1}|t_1:....:t_{m+1}\}:=  
\end{equation}
$$
\sum_{n_i >0} {\rm I}_{n_1,...,n_m}(g_1: ... : g_{m+1})
(t_1-t_{m+1})^{n_1-1}...(t_m -t_{m+1})^{n_m-1}
$$
We picture them on the oriented circle. Namely, the circle has slots, where the $g$'s sit, and in between the consecutive slots, 
dual slots, where $t$'s sit, see the picture below. The slots are marked by black points, and the 
dual slots by little circles:
\begin{center}
\hspace{4.0cm}
\epsffile{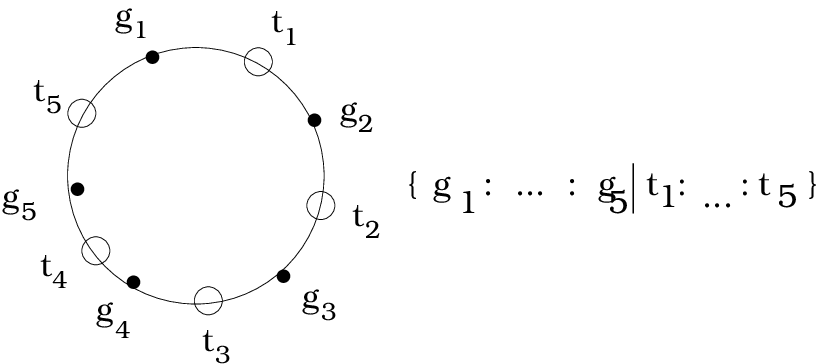}
\end{center}
Set
\begin{equation} \label{ccc3}
\delta \{g_1:  ... : g_{m+1}| t_1: ... :t_{m+1}\} = 
\end{equation}
$$
-\sum_{k=2}^{m} {\rm Cycle}_{m+1}\Bigl(\{g_{1}:... :g_{k-1}:g_k| t_{1}: ... : t_{k-1}: t_{m+1}    \} 
\wedge 
$$
$$
\{g_{k}  :   ... : g_{m+1}| t_k: ... : t_{m+1}  \} \Bigr)
$$ 
where the indices are modulo $m+1$ and  
$$
 {\rm Cycle}_{m+1} f(x_1,...,x_{m+1}) := \sum_{i = 1}^{m+1} f(x_i,...,x_{m+i})
$$  
Each term of the formula corresponds to the 
following procedure: choose a slot and 
a dual slot on the circle. 
Cut the circle 
 at the chosen  slot and dual slot  and make two   oriented circles      
with a dihedral words on each of them   out of the initial data. It is useful to 
think about the slots and dual slots as of little arcs, not points, so cutting one of them we 
get the 
arcs on each of the two new circles marked by the corresponding letters. 
The formula reads as follows: 
$$
\delta((\ref{ccc3})) = - \sum_{{\rm cuts}} 
\mbox{(start at the dual slot)} \quad  \wedge\quad 
\mbox{(start at the slot)}
$$ 
The only asymmetry between $g$'s and $t$'s is the order of factors.

\begin{center}
\hspace{4.0cm}
\epsffile{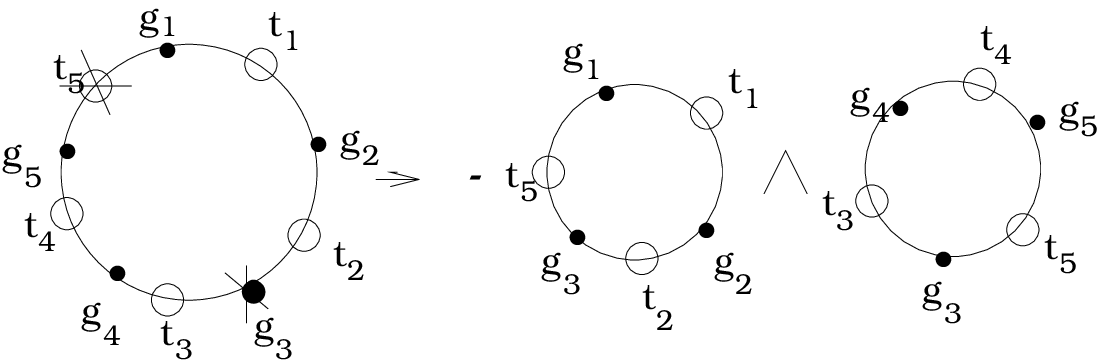}
\end{center}

{\it Comparing $\delta^{\cal H}$ and $\delta$}. We claim that the only terms in the sum  which may 
survive after we kill the products of the weight $\geq 1$ terms are the ones 
corresponding to  special marked colored segments shown on the two pictures below
\begin{center}
\hspace{4.0cm}
\epsffile{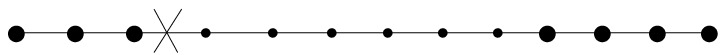}
\end{center}

\begin{center}
\hspace{4.0cm}
\epsffile{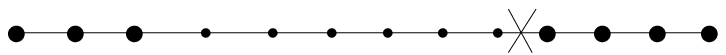}
\end{center}
More precisely, these are the special marked colored segments where we  mark every point from 
$a_0=0$ to  
some point $a_p$, then  jump on the right and mark a point $a_q$  and  every 
point after, till we hit the right endpoint $a_{m+1}$. 
The cross between $a_p$ and $a_q$ can be either at the very left segment 
$t_p$, or at the very right segment $t_{q-1}$. 

To check the claim we observe the following. The factors in the product (\ref{ur100100}) correspond 
to  ``gaps'' between the consequetive  marked points. We say that such a gap is wide, if 
it has more then   one segment inside, i.e. if there is unmarked point inside  the gap, as shown on 
 the picture
\begin{center}
\hspace{4.0cm}
\epsffile{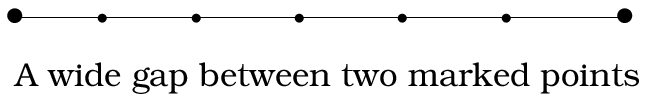}
\end{center}

The factor of the product (\ref{ur100100}) corresponding to a wide gap is of weight at least one. 
Thus if there are two such wide gaps on a marked colored segment, then 
the corresponding term dies after we kill the products. 
Therefore we allow to have 
no more then  one wide gap. 
\begin{center}
\hspace{4.0cm}
\epsffile{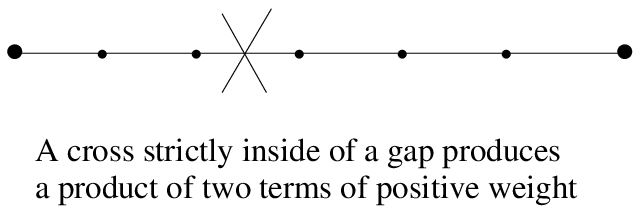}
\end{center}
Now let us make a close look at the cross located inside of the wide gap. If this cross is not 
in one of the corner segments, we  will pick up two terms of weight at least one 
in the term corresponding to this gap. 
So the cross must be in one of the corner segments of the gap, and we arrive precisely to the 
marked colored segments shown on the two pictures above.  

Let us show that the contribution of all special marked colored segments having a wide gap 
matches the formula for $\delta$. Observe that 
$
{\rm I}^{\cal H}(a|t) = a^t
$. 
 If there is a wide gap on a special marked colored segment then every other gap 
contributes just $1$  to the product formula. Indeed, its contribution is 
$$
{\rm I}^{\cal H}(a|t) = a^t = {\rm e}^{a \cdot t} =1 + \log^{\cal H}a\cdot t + 
(\log^{\cal H}a)^2\cdot t^2/2 
+ ... 
$$
 but any term  
of weight bigger then zero will multiply the term corresponding to the wide gap, 
and thus their product will be killed by going to the Lie coalgebra. 
 Similarly the contribution of the wide gap with the left marked cross is 
\begin{equation} \label{1.26.02.1}
(-1)^{i_{\alpha +1}-1 - i_{\alpha }} {\rm I}^{\cal H}(a_{i_{\alpha+1}-1}: a_{i_{\alpha+1}-2}: ... : 
 a_{i_{\alpha}}| 
-t_{i_{\alpha+1}-1}: -t_{i_{\alpha+1}-2}: ... : - t_{i_{\alpha}})
\end{equation}
It corresponds to the case $j_{\alpha} = i_{\alpha+1}-1$ in (\ref{ur100100}). 
The contribution of the wide gap with the right marked cross is
\begin{equation} \label{1.26.02.r}
{\rm I}^{\cal H}(a_{i_{\alpha}+1}: a_{i_{\alpha}+2}: ... : 
 a_{i_{\alpha+1}}| 
t_{i_{\alpha}}: t_{i_{\alpha}+1}: ... : t_{i_{\alpha+1}-1})
\end{equation}
It follows from  theorem 4.1 and formula (65) in [G4] that, 
modulo the double shuffle relations, lower depth terms and products the element 
(\ref{1.26.02.1}) is equivalent to the element 
\begin{equation} \label{1.26.02.l}
 -{\rm I}^{\cal H}(a_{i_{\alpha}}: ... :  a_{i_{\alpha+1}-2}
 : a_{i_{\alpha+1}-1}| 
t_{i_{\alpha}+1}: t_{i_{\alpha}+2}: ... :  t_{i_{\alpha}})
\end{equation}

Observe now that the map $\nu^{\cal H}_{\bullet, \bullet}(\mu_N)$ is described geometrically by the picture
\begin{center}
\hspace{4.0cm}
\epsffile{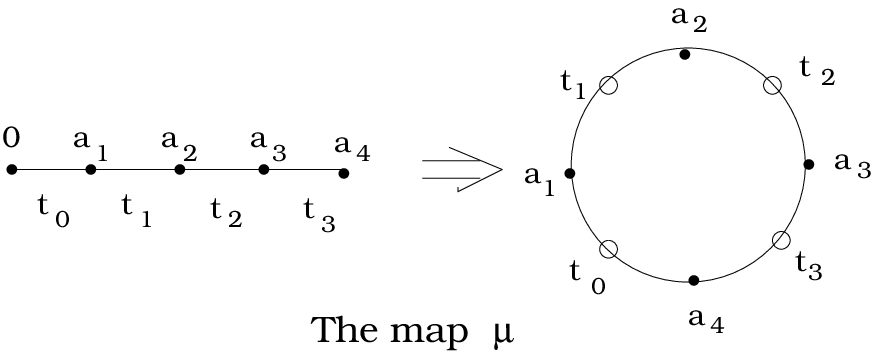}
\end{center}

Now it is easy to check  that the contribution of all special 
marked colored segments having a wide gap 
matches the formula for $\delta$. Namely, the two ways to mark by cross the wide gap 
produce the following two pictures. The first one corresponds to (\ref{1.26.02.l}):

\begin{center}
\hspace{4.0cm}
\epsffile{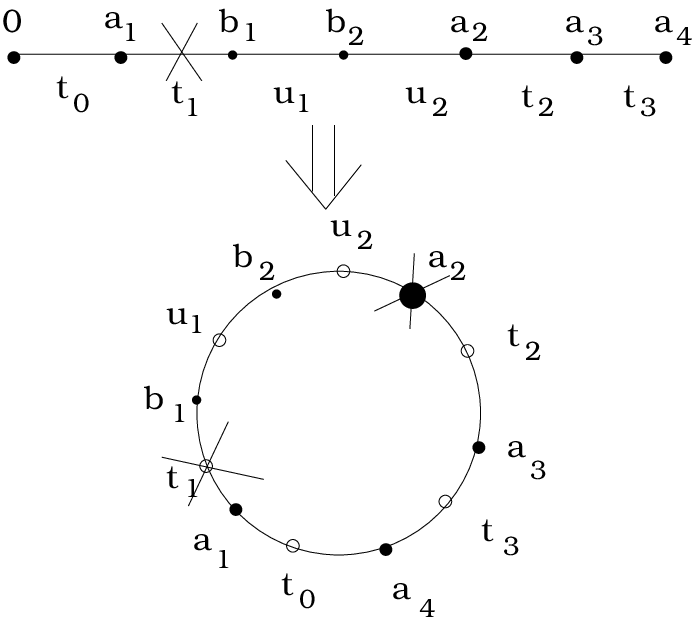}
\end{center}

and the second one to (\ref{1.26.02.r}):

\begin{center}
\hspace{4.0cm}
\epsffile{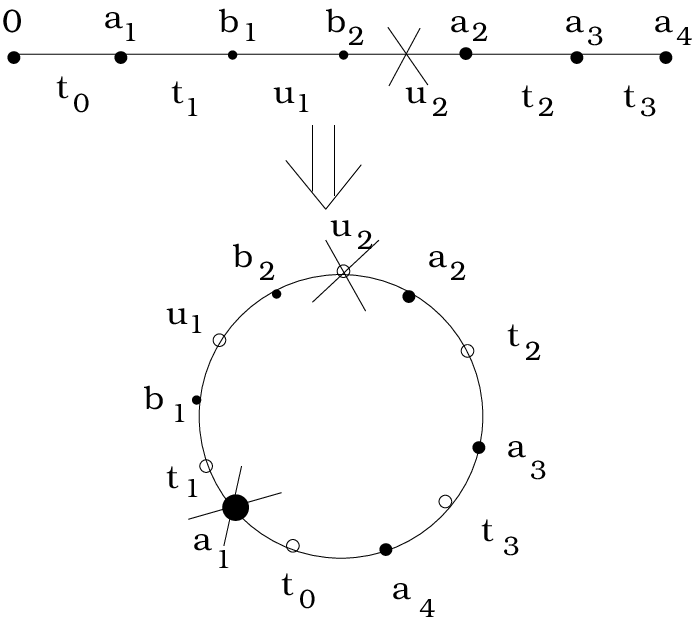}
\end{center}
Observe that the minus  sign in (\ref{1.26.02.l}) comparing to the plus sign in (\ref{1.26.02.r}) 
agree with the fact that the order of the terms on the first picture is 
different from the one prescribed by formula (\ref{ccc3}). So interchanging the 
terms we get a minus sign as well, in according with minus sign in front of (\ref{ccc3}).

So far we did not use the fact that $a_i \in \mu_N$.  Now let us use it to show that 
the contribution of the unique special marked colored segment with no wide gaps is zero. 
Indeed,  $\log^{\cal H} x$ is an $N$-torsion element if $x^N=1$. 
So working over $\Q$ (or at least modulo $N$-torsion) we 
have $\log^{\cal H} x =0$. Thus 
$$
x^t:= {\rm exp}(\log^{\cal H}x  \cdot t) = 1\quad \mbox{provided $x^N=1$}
$$ 
and hence there is no contribution from this term to the restricted 
coproduct $\Delta'$, and hence to $\delta^{\cal H}$. 
We proved that the map $\nu^{\cal H}_{\bullet, \bullet}(\mu_N)$ commutes 
with the cobrackets, and hence 
finished the proof of theorem \ref{4.16.01.pi}.

{\bf Remark}. The very definition of the cobracket in the dihedral Lie algebra 
given in chapter 4 of [G4] emerged from the calculation of  the coproduct 
for the framed multiple polylogarithm Hodge-Tate structures.

{\bf 8. Proof of theorem \ref{2.14.02.1}}. We start from the $w=m$ case

\begin{proposition} \label{1.30.02.5} The dual to the map 
$\nu^{\cal M}_{\bullet, \bullet}(\mu_N)$ from theorem \ref{4.16.01.pi} induces 
isomorphisms
$$
{\rm C}^{\cal M}_{m, m}(\mu_p) = {\rm D}_{m, m}(\mu_p) \quad \mbox{for $m=1, 2,3$}
$$ 
\end{proposition}

{\bf Proof}. For $m=1$ this is given by the classical Bass theorem on cyclotomic units. 
The  $m=2$ case can be deduced from lemma 7.11 in [G4], which claims that 
this is true after tensoring by $\Q_l$. In fact using theorem \ref{4.16.01.pi} it is easy to 
translate the the proof of that lemma into motivic setting.  

In the $m=3$ case we observe that it follows from from section 7.7 [G4], 
see  in particular isomorphism (172) there, that the subspace of 
 ${\rm C}_{-3, -3}(\mu_p)$ generated by the triple commutators 
of elements in ${\rm C}_{-1, -1}(\mu_p)$ is isomorphic to 
${\rm D}_{-3, -3}(\mu_p)$. On the other hand by theorem 
\ref{4.16.01.pi} there is an inclusion 
${\rm C}_{-3, -3}(\mu_p) \subset {\rm D}_{-3, -3}(\mu_p)$. The proposition is proved. 

Theorem \ref{2.14.02.1} in the case $N=1$, $m=1$ just claims that 
$\zeta^{\cal M}(2n)=0$ and $\zeta^{\cal M}(2n+1)\not =0$. In the Hodge realization both 
assertions are well known: the first follows from 
 the distribution relations (for $l=2$), and the second 
from the fact that the real period of  $\zeta^{\cal H}(2n+1)$ is not zero.

In the $N=1$ theorem \ref{2.14.02.1} follows from theorem 2.4 (when $m=2$) 
and theorem 2.10 (when $m=3$) in [G4]. Indeed, those two theorems prove similar 
results in the l-adic setting, and so thanks to  theorem \ref{5.2.02.2} imply 
theorem \ref{2.14.02.1} for $N=1$, $m=2,3$ cases. In fact by repeating the arguments 
used to prove theorems 2.4 and 2.10 in [G4] we get a direct proof of theorem 
\ref{2.14.02.1} in these cases. 

{\bf 9.  Proof of theorem \ref{1.30.02.7}}. 
Combining proposition \ref{1.30.02.5} 
with formulas (169)-(170) in [G4] we get theorem \ref{1.30.02.7}. 

{\bf Remark}. One should correct the apparent misprints in (169)-(170) in [G4] by changing 
$H^i_{(w,2)}$ to $H^i_{(2)}$ and $H^i_{(w,3)}$ to $H^i_{(3)}$. 

Let $d^*_{w,m}(p) = {\rm dim}{\rm C}^{\cal M}_{w,m}(\mu_p)$. 

 \begin{corollary} \label{1.30.02.6} One has for prime $p \geq 5$
$$
d^*_{2,2}(p) = \frac{(p-1)(p-5)}{12}; \quad d^*_{3,3}(p) = \frac{(p-5)(p^2 - 2p - 11)}{48}
$$
\end{corollary}

{\bf Proof}. This is the motivic version of the corollary 2.16 
refined in the depth three case. The proof goes just like the proof of that corollary 
in chapter 7 of [G4], but now in the motivic setting. Namely,  
we calculate the Euler characteristic the $(m,m)$ part of the standard cochain complex 
for the Lie algebra ${\rm C}^{\cal M}_{\bullet, \bullet}(\mu_p)$ using theorem \ref{1.30.02.7}. 

Observe that for $i=2,3$ the group $H^{i}(\Gamma_1(3;p), \Q)$ contain 
the mysterious cuspidal cohomology subgroup $H_{\rm cusp}^i(\Gamma_1(3;p), \Q)$. 
There is no closed 
 formula for its dimension, but thanks to the Poincare duality the cuspidal cohomology for 
$i=2$ and $i=3$  cancel each other 
in the computation of the Euler characteristic. As a result we were 
 able to get a closed formula for $d^*_{3,3}(p)$. The corollary \ref{1.30.02.6} is proved.

{\bf 10. Proof of theorem \ref{4.16.01.qoin}}. By theorem 1.8 in [G7] there is 
canonical surjective homomorphism 
$$
{\cal Z}^{\cal H}_{\bullet}(\mu_N) \lra {\rm Gr}^W\widetilde {\cal Z}(\mu_N)
$$ 
Let us choose a complex primitive $N$-th root of unity. Then thanks to 
theorem \ref{4.16.01.13sa1}a) the Hodge realization functor 
 provides a surjective homomorphism
$$
{\cal Z}^{\cal M}_{\bullet}(\mu_N) \lra {\cal Z}^{\cal H}_{\bullet}(\mu_N) 
$$
By  
theorem \ref{4.16.01.qa} ${\cal Z}^{\cal M}_{\bullet}(\mu_N)$ is a subalgebra of 
${\cal A}^{\cal M}_{\bullet}(S_N)$.  The theorem is proved.

{\bf 11. Understanding ${\cal Z}_{\bullet}(S_N)$}.
 It would be very interesting to calculate 
$$
d_w(N):= {\rm dim}{\cal Z}_{w}(S_N) \quad \mbox{and} \quad 
d_{w,m}(N):= {\rm dim}{\rm Gr}^D_m{\cal Z}_{w}(S_N)
$$ as  functions of 
$N$ for a given  $(w,m)$.  
According to Deligne [D2-3] ${\cal A}_{\bullet}(S_N)$ can  not be smaller then 
${\cal Z}_{\bullet}^{\cal M}(\mu_N)$ for $N = 2,3,4$. This  plus theorem 
\ref{4.16.01.qa} yields
\begin{equation} \label{1.30.02.2}
{\cal A}_{\bullet}(S_N) = 
{\cal Z}_{\bullet}^{\cal M}(\mu_N) \quad \mbox{for $N = 2,3,4$}
\end{equation}
According to [G2] or [G4] this is not true for sufficiently big $N$, e.g. for prime $N\geq 5$. 
The case 
$N=1$ remains open, as well as the case $N=6$.

\begin{conjecture} \label{1.29.02.1} If $p$ is a prime  then 
$d_w(p)$ and $d_{w,m}(p)$ are polynomials in $p$.
\end{conjecture} 

For $p=2,3$ this follows from   (\ref{1.30.02.2}).

Conjecture \ref{1.29.02.1} is equivalent to a similar statement about $d^*_{w,m}(p)$. 
Corollary \ref{1.30.02.6} supports it for $m=w=2$ and $m=w=3$. According to the results of [G10] we
have  the same situation for $m=w=4$.

{\bf 12. Proof of theorem \ref{2.14.02.2}}. We need to check the second isomorphism. 
The crucial fact is the following formula for the coproduct of the classical 
polylogarithm motive. 

\begin{theorem} \label{4.22.02.15} Let $F$ be a number field. Then for any $x \in F$ 
one has 
$$
\Delta {\rm Li}^{\cal M}_n(x)  = {\rm Li}^{\cal M}_n(x) \otimes 1 + 1 
\otimes {\rm Li}^{\cal M}_n(x) + 
\sum_{k=1}^{n-1} {\rm Li}^{\cal M}_{n-k}(x) 
\otimes \frac{(\log^{\cal M}x)^k}{k!} 
$$
\end{theorem}

{\bf Proof}. In the Hodge version this result is well known, see for instance chapter 5 
of [G7]. Using lemma 3.4 from [G7] we transfer it to the motivic setting. 
The theorem is proved. 

It follows from this that for any root of unity $\zeta$ the element 
${\rm Li}^{\cal M}_n(\zeta)$ is primitive 
(observe that it lies in the $\Q$-vector space). Thus it provides an class
$$
{\rm Li}^{\cal M}_n(\zeta) \in {\rm Ext}^1(\Q(0), \Q(n)) = K_{2n-1}(\Q(\zeta)) \otimes \Q 
$$
For $n=1$ theorem \ref{2.14.02.2} is  classical. 
Assume $n>1$. Let $\alpha$ be a residue modulo $N$. Consider the classes 
\begin{equation} \label{4.22.02.20}
{\rm Li}^{\cal M}_n(\zeta_N^{\alpha}); \qquad (\alpha,N)=1, \quad 1 \leq \alpha \leq N/2
\end{equation}
Using the distribution relations we check that  these classes span 
the space genertaed by ${\rm Li}^{\cal M}_n(\zeta_N^{\alpha})$ for all $\alpha$. 
We claim that the elements (\ref{4.22.02.20})  generate $K_{2n-1}(\Q(\zeta)) \otimes \Q$. 
We need only to check that these elements 
are linerly independent. Thanks to lemma 3.4 in [G7] this can be done 
in the Hodge realization, where it is well known 
and established using Dirichlet's formula for 
$\zeta_{\Q(\mu_N)}(n)$ plus the fact that this special value is non zero. 
Theorem \ref{2.14.02.2} is proved. 

{\bf 13. The weak version of Zagier's conjecture on the special values 
of the Dedekind $\zeta$ functions at $s=n$}. It is deduced 
from theorem \ref{4.22.02.15} using the arguments given in  [BD1], or section 12 of [G1]. 
Another proves has been given in [BD2] and [DJ].

\section{Apendix:  ${\bf Li}$-shuffle Hopf algebras}

We define 
for an arbitrary commutative group $G$ a commutative graded 
Hopf algebra $Sh_{\bullet}^{\rm Li}(G)$. Its generators are the symbols 
\begin{equation} \label{3.21.01.1}
\widehat 
{\rm Li}_{n_1, ..., n_m}(x_1, ..., x_m); \qquad x_i \in G, \quad n_i \in \Z_+
\end{equation}
The relations  are given by the  ${\rm Li}$-shuffle product 
 and the inversion formulas. The formula for the coproduct $\Delta$ 
in this Hopf algebra was suggested by the formula for the coproduct in the 
multiple polylogarithm 
Hopf algebra ${\cal Z}_{\bullet}(\C^*)$  obtained in [G7]. 
The inversion formula reflects the one obtained in s. 2.6 in [G7]. 
The key result is that the  Hopf axiom holds. This 
enables us to give another proof of the motivic ${\rm Li}$-shuffle 
relations for generic parameters $x_i$. 

{\bf 1. The ${\bf Li}$-shuffle relations revisited}. 
Recall the generating formal power series for multiple polylogarithms:
$$
 {\rm Li}(x_1, ..., x_m| t_1, ... , t_m):= \sum_{n_i >0}
{\rm Li}_{n_1, ..., n_m}(x_1, ... , x_m)t_1^{n_1-1}... t_m^{n_m-1}
$$ 
The ${\rm Li}$-shuffle relations for these generating series look as follows: 
\begin{equation} \label{11.21.0.1}
{\rm Li}(x_1,...,x_m| t_1,...,t_m)\cdot {\rm Li}(x_{m+1},...,x_{m+n}|t_{m+1},...,t_{m+n}) \quad =
\end{equation}
$$
\sum_{\sigma \in \overline \Sigma_{m,n}}{\rm Li}_{\sigma}(x_{1},...,x_{ m+n }
|t_{1}, ...,t_{m+n}) 
$$
We need to 
 make explicit the term ${\rm Li}_{\sigma}(-|-)$ 
corresponding to  a generalized shuffle $\sigma$. 
For example
$$
{\rm Li}(x_1|t_1) \cdot {\rm Li}(x_2|t_2) =  {\rm Li}(x_1, x_2|t_1,  t_2) + {\rm Li}(x_2, x_1|t_2,  t_1) + 
$$
$$
\frac{1}{t_1-t_2}
\Bigl( {\rm Li}_1(x_1x_2|t_1 ) - {\rm Li}_1(x_1x_2|t_2 )\Bigr) 
$$
Recall that a generalized shuffle of two strings 
\begin{equation} \label{11.22.00.2}
\{x_1,...,x_{m}\} \quad \mbox{and} \quad \{y_{1},...,y_{n}\}
\end{equation}
is given by the following data:

i) A pair of subsets $I,J$ where 
\begin{equation} \label{11.22.00.1}
I= \{1 \leq i_1 < ... < i_p \leq m\}, \qquad J= \{1 \leq j_1 < ... < j_p \leq n\}, \quad p \geq 0 
\end{equation}
 providing  $p$ 
{\rm special pairs} 
$(x_{i_1}, y_{j_1}), ..., (x_{i_p}, y_{j_p})$. Then each of  the complements 
$$
\{x_1,...,x_{m}\} - \{x_{i_1}, ..., x_{i_p}\} \quad \mbox{and} \quad \{y_{1},...,y_{n}\}
- \{y_{j_1}, ..., y_{j_p}\} 
$$
is a union of $p+1$ substrings (some of them might be empty). 

ii)  A choice of $p+1$  permutations $\sigma_0, ..., \sigma_p$,  
where  $\sigma_i$ shuffles the elements of the $i$-th substrings 
of (\ref{11.22.00.2}).

Consider the operator 
$$
D_{ij} f(x_i,x_j|t_i,t_j):= \quad \frac{1}{t_i-t_j} \Bigl(f(x_ix_j|t_i) - f(x_ix_j|t_j)\Bigr)
$$
A generalized shuffle $\sigma$ provides a well defined usual shuffle $\sigma'$ 
obtained by shifting a bit to the right the $y$ variables in each special pair. 
For a usual shuffle $\sigma'$ the element ${\rm Li}_{\sigma'}(-|-)$ has an obvious meaning
Set
$$
{\rm Li}_{\sigma}(x_1, ..., x_{m+n}|t_1, ..., t_{m+n}):= D_{i_1j_1} ... D_{i_pj_p}
{\rm Li}_{\sigma'}(x_1, ..., x_{m+n}|t_1, ..., t_{m+n})
$$

We say that a generalized shuffle has 
depth $-p$ if $|I| = |J| =p$. Denote the subset of all such shuffles by 
$\overline \Sigma^p_{m,n}$, so $\overline  \Sigma_{m,n} = 
\cup_{p\geq 0}\overline  \Sigma^p_{m,n}$.

{\bf 2. The ${\bf Li}$-shuffle Hopf algebra}. 
Let $G$ be a commutative group. We are going to define a 
commutative algebra with unit   $Sh_{\bullet}^{\rm Li}(G)$ 
the over $\Q$. By definition it is  generated by the symbols 
${\widehat {\rm Li}}_{n_1, ..., n_m}(x_1, ... , x_m)$, as in 
(\ref{3.21.01.1}). To present the relations we need to introduce 
some notations. 
Let us package the generators into the 
formal power series in $t_i$:
$$
 {\widehat {\rm Li}}(x_1, ..., x_m| t_1, ... , t_m):= \sum_{n_i >0}
{\widehat {\rm Li}}_{n_1, ..., n_m}(x_1, ... , x_m)t_1^{n_1-1}... t_m^{n_m-1}
$$ 
 Let 
 $x_0 ... x_m =1$. We will also use notations
\begin{equation} \label{12.10.00.2}
{\widehat {\rm Li}}(*, x_1, ..., x_m| t_0: ... :t_m) = {\widehat {\rm Li}}(x_0, x_1, ..., x_m| t_0: ... :t_m) := 
\end{equation}
$$
{\widehat {\rm Li}}(x_1, ..., x_m| t_1-t_0, ... , t_m-t_0)
$$

Let 
$
{\widehat \log} (x) := {\widehat {\rm Li}}_1(x^{-1}) - {\widehat {\rm Li}}_1(x)
$. 
We  employ the notation 
$$
x^t:= \quad e^{t \cdot {\widehat \log} (x)}:= \quad 1+ \sum_{n > 0} 
{\widehat \log}^n (x)\frac{t^n}{n!};
\qquad \widehat B(x|t)  = x^t/t
$$

{\bf Relations}
1. {\it The additivity of the formal logarithm}. 
Then 
$$
{\widehat \log} (xy) = {\widehat \log} (x) + {\widehat \log} (y) 
\qquad \mbox{for any $x,y \in G$}
$$

2. {\it The ${\rm Li}$-shuffle product relations}.
$$
{\widehat {\rm Li}}(*, x_1, ..., x_m| t_0: t_1: ... : t_m) \cdot 
{\widehat {\rm Li}}(*, x_{m+1}, ..., x_{m+n}| t_0: t_{m+1}: ...: 
t_{m+n}) = 
$$
\begin{equation} \label{6.ne111}
\sum_{\sigma \in \overline \Sigma_{m,n}}\widehat{\rm Li}_{\sigma}(*, x_{1},...,x_{ m+n }
|t_0:  t_{1}: ...:t_{m+n}) 
\end{equation}
where $\sigma$ shuffles the variables $t_1, ..., t_{m+n}$. 
Set  
$$
{\widehat B}(x_1, ..., x_m | t_1, ..., t_m):= 
$$
\begin{equation} \label{12.3.00.3}
\sum_{j=1}^m (-1)^{j-1}
{\widehat {\rm Li}}(*, x^{-1}_{j-1}, ..., x^{-1}_1| -t_{j}: -t_{j-1}: ... :
 - t_{1}) \cdot 
\end{equation}
\begin{equation} \label{12.1.00.4}
\widehat B(x_1 ... x_m|t_j) \cdot 
{\widehat {\rm Li}}(*, x_{j+1}, ...,  x_m| t_{j}: t_{j+1}: ... : t_{m}) 
\end{equation}

{\bf Remark}. If $G = \{z \in \C^*;  |z| =1\}$ the inversion formula 
is a formal 
version  the inversion formula derived in  section 2.6 of [G7].

{\it 3. The inversion formula}. 
$$
{\widehat B}(x_1, ..., x_m | t_1, ..., t_m) = 
$$
\begin{equation} \label{12.3.00.1}
\sum_{j=1}^m (-1)^j {\widehat {\rm Li}}(x^{-1}_j, ..., x^{-1}_1 | -t_j, ..., -t_1) 
{\widehat {\rm Li}}(x_{j+1}, ..., x_m | t_{j+1}, ..., t_m) + 
\end{equation}
$$
\sum_{j=1}^m \frac{(-1)^{j}}{t_j}
{\widehat {\rm Li}}(x^{-1}_{j-1}, ..., x^{-1}_1| -t_{j-1}, ...,
- t_{1}) \cdot 
{\widehat {\rm Li}}(x_{j+1}, ..., x_m | t_{j+1}, ..., t_m) 
$$

Set $X_{a\to b}:= \prod_{s=a}^{b-1}x_s$. 
\begin{definition} \label{12.10.00.1}
We define a linear map $\Delta: Sh_{\bullet}^{\rm Li}(G) \lms Sh_{\bullet}^{\rm Li}(G) \otimes 
Sh_{\bullet}^{\rm Li}(G)$ by 
$$
\Delta {\widehat {\rm Li}}(x_0, x_1, ..., x_m | t_0: t_1: ... :t_{m}) = 
$$
$$
\sum {\widehat {\rm Li}}(X_{i_0 \to i_{1}}, X_{i_1 \to i_{2}}, ... , X_{i_k \to m}| 
t_{j_0}: t_{j_1}:  ... :
  t_{j_{k}}) \otimes 
$$
\begin{equation}  \label{ur1001222}
\prod_{p = 0}^k 
\Bigl( X_{i_p \to i_{p+1}}^{t_{j_p}} \cdot (-1)^{j_p-i_p}\cdot {\widehat {\rm Li}}(*, x^{-1}_{{j_p}-1},  ..., x^{-1}_{i_p}|-t_{j_{p}}:  
-t_{j_{p}-1}:  ... : -t_{i_{p}}) \cdot 
\end{equation}
\begin{equation}  \label{11.16.00.2}
{\widehat {\rm Li}}(*, x_{j_{p}+1}, x_{j_{p}+2}, ..., x_{i_{p+1}-1}| 
t_{j_{p}}: t_{j_{p}+1}:  ... : t_{i_{p+1}-1}) \Bigr)  
\end{equation} 
where the sum is over all 
sequences $\{i_{\alpha}\}$ and $\{j_{\alpha}\}$ 
satisfying conditions 
\begin{equation}  \label{ur10*}
i_{\alpha } \leq j_{\alpha } < i_{\alpha +1} \quad \mbox{ for any 
$0\leq \alpha \leq k$}, \qquad j_0 = i_0 = 0
\end{equation} 
\end{definition}

Observe that this map is well defined on the generators. 
This is not granted for free by the formula for $\Delta$ since,  
 according to  definition (\ref{12.10.00.2}),  
the shift  
$t_i \lms t_i +t$ does not change the generating series 
${\widehat {\rm Li}}(x_0, x_1, ..., x_m | t_0: t_1: ... :t_{m})$. However 
since $x_0 \cdot ... \cdot x_m =1$ it also does not change 
 formula (\ref{ur1001222}): 
one has $\prod_{p=1}^k X_{i_p \to i_{p+1}}^{t} = (x_0 \cdot ... \cdot x_m)^t =1$.

\begin{theorem} \label{ur100111}
$\Delta$ is a well defined homomorphism of algebras.
\end{theorem} 

{\bf Proof}. This statement is equivalent to  a collection of statements, one 
 for each weight $w$. We will use the induction on $w$.

Let $\widetilde Sh_{\bullet}^{\rm Li}(G)$ be the graded vector space 
generated by the symbols (\ref{3.21.01.1}) 
subject to the relations 1 and 3. One  shows that 
the formula for $\Delta$ provide a well defined map of graded linear vector spaces 
$\Delta: \widetilde Sh_{\bullet}^{\rm Li}(G) \lms \widetilde Sh_{\bullet}^{\rm Li}(G) \otimes 
\widetilde Sh_{\bullet}^{\rm Li} (G)$. We left the details to the reader.

Now we need to handle the shuffle relations. This is the complicated part of the proof, 
and we provide a detailed account of it. 
We picture a generator ${\widehat {\rm Li}}(x_0, ..., x_m|t_0: ... : t_m)$ 
by a segment subdivided into $m+1$ arcs labelled by  $(x_i|t_i)$. 
 Every  term  of the formula for the coproduct of this generator 
is determined by the following geometric procedure: 
we subdivide this segment into a union of 
$k+1$ little  segments, each being a union of arcs. 
We  mark (by cross) an arc on each of these little segments. 
We assume in addition that this marking is special, 
i.e. in the very left segment we  
always mark the very left arc. 

The product 
\begin{equation} \label{11.14.00.2}
{\widehat {\rm Li}}(*, x_1, ..., x_m| t_0: t_1: ... : t_m) \cdot 
{\widehat {\rm Li}}(*, y_1, ..., y_n| t_0: s_1: ... : s_n) 
\end{equation}
is given by the sum of the terms corresponding to the generalized shuffles of symbols 
$(x_i|t_i)$ and $(y_j|s_j)$. More precisely, a given generalized shuffle of depth $-p$ provides
$2^p$ terms in the formula. 
Thus applying $\Delta$ to this sum we get  
 sum of the terms corresponding to the following data:

{\it a generalized shuffle for (\ref{11.14.00.2}) plus a choice of one of the $2^p$ related 
terms; 

a subdivision of the corresponding segment 
into little marked segments}.

Choose integers $ 1 \leq a \leq  b \leq m$ and $ 1 \leq c \leq d \leq n$. 
Then there is the following subset of $x_i$'s and $y_j$'s:
\begin{equation} \label{11.14.00.1}
\{x_a, ..., x_{b-1}\} \cup \{y_c, y_{c+1}, ..., y_{d-1}\}
\end{equation}
Consider all the terms which correspond to the above data such that  the 
product of $x_i$'s and $y_j$'s on the p-th little segment equals to the product of 
elements in (\ref{11.14.00.1}). 
Then there are four possible types of marking of the $p$-th little segment: 
$$
{\rm i}) \quad (x_i|t_i), \qquad {\rm ii}) \quad (y_j|s_j), \qquad {\rm iii}) 
\quad (x_iy_j|t_i)  \qquad {\rm iv}) \quad (x_iy_j|s_j)
$$ 
where $a \leq i \leq b, c \leq j \leq d$. The corresponding left 
factors of the terms in the coproduct 
look as follows:  
\begin{equation} \label{11.13.00.5}
{\widehat {\rm Li}}( ... , X_{a \to b}\cdot Y_{c \to d}, ... | ... : t_i: ... ) 
\end{equation}
$$
{\widehat {\rm Li}}( ... , X_{a \to b}\cdot Y_{c \to d}, ... | ... : s_j: ... ) 
$$
\begin{equation} \label{11.18.00.1}
{\widehat {\rm Li}}( ... , X_{a \to b}\cdot Y_{c \to d}, ... | ... : t_i: ... ) 
\frac{1}{t_i - s_j}
\end{equation}
$$
{\widehat {\rm Li}}( ... , X_{a \to b}\cdot Y_{c \to d}, ... | ... : s_j: ... ) 
\frac{1}{s_j - t_i}
$$
Let us assume first that $a<b, c<d$. 
We claim that then  the sum of the terms in the formula for 
$\Delta$ over all 
generalized shuffles with the markings of $p$-th little segments of types  i) 
(i.e. in the left factors of these terms  
are as in 
 (\ref{11.14.00.1})) 
equals to 
$$
 (\ref{11.13.00.5}) \otimes (X_{a\to b} \cdot Y_{c\to d})^{t_i}\cdot 
\Bigl(\sum_{k=c}^{d-1}(-1)^{i+k-a - c+1}
\cdot {\widehat {\rm Li}}
(*, x^{-1}_{i-1}, ..., x^{-1}_a|-t_i: -t_{i-1} : ... : -t_a) \cdot
$$
\begin{equation} \label{11.19.00.1}
{\widehat {\rm Li}}(*, y^{-1}_{k},   ..., y^{-1}_c|  -t_i: -s_{k}:     ... : -s_c) \cdot
{\widehat {\rm Li}}(*, y_{k+1},  ..., y_{d-1}| t_i: s_{k+1}:  ...: s_{d-1})  
\end{equation}
$$
{\widehat {\rm Li}}(*, x_{i+1}, ..., x_{b-1}|t_i: t_{i+1}: ...: t_{b-1})\Bigr) \cdot\prod_{q \not = p}(\mbox{contribution of  $q$-th little segment})
$$
To check this formula at any given weight $w$ 
we may use the shuffle product relations of weights $<w$. 
Consider the generalized shuffles related to product (\ref{11.14.00.2}) 
where in addition the arcs on the $p$-th little segment located on the right of 
$(x_i|t_i)$-arc are labelled by the elements of $\{x_{i+1}, ..., x_{b-1}\} 
\cup \{y_{k+1}, ..., y_{d-1}\}$. Then using the shuffle product formula we see that  contribution 
of factor (\ref{11.16.00.2}) 
to $\Delta( \mbox{ sum of such generalized shuffles})$  is 
(factor 3) $\times$ (factor 4) in the sum (\ref{11.19.00.1}). 
The use of the shuffle product formula here is legitimate since we proceed by induction, and 
we may assume that the left factor has a nonzero weight. 
 
The same kind of arguments,   which take into account 
the signs in  (\ref{ur1001222}),
show that  (\ref{ur1001222}) contributes (factor 1) $\times$ (factor 2) 
in the sum (\ref{11.19.00.1}).  
 Formula (\ref{11.19.00.1}) is proved. 

A completely similar considerations for the terms of type iii) lead to a similar 
to (\ref{11.19.00.1}) formula where the (factor 2) $\times$ (factor 3) is replaced by 
\begin{equation} \label{12.10.00.10}
\frac{1}{t_i-s_k}\cdot {\widehat {\rm Li}}(*, y^{-1}_{k-1},   ..., y^{-1}_c|
  -t_i: -s_{k-1}:     ... : -s_c) \cdot 
\end{equation}
$$
{\widehat {\rm Li}}(*, y_{k+1},  ..., y_{d-1}| t_i: s_{k+1}:  ...: s_{d-1}) 
$$
By the inversion formula (\ref{12.3.00.1}) 
$$
\sum_{k=c}^{d-1}\Bigl( (\ref{12.10.00.10}) + (\mbox{factor $2$}) \cdot  
(\mbox{factor $3$) in  
 (\ref{11.19.00.1})}
 \Bigr)  =  \widehat B(y_c, ..., y_{d-1}|s_c-t_i, ..., s_{d-1} - t_i)
$$
 Using definition (\ref{12.3.00.3}) we see that 
$$
Y^{t_i}_{c \to d}\widehat B(y_c, ..., y_{d-1}|s_c-t_i, ..., s_{d-1} - t_i) = 
$$
$$
\sum_{j=c}^{d-1} (-1)^{j-c}\frac{Y^{s_j}_{c \to d}}{s_j - t_i}\cdot
{\widehat {\rm Li}}(y^{-1}_{j-1}, ..., y^{-1}_c| s_{j}-s_{j-1}, ...,
s_{j} - s_{c}) \cdot 
$$
$$
{\widehat {\rm Li}}(y_{j+1}, ...,  y_{d-1}| s_{j+1}- s_{j}, ...,
s_{d-1} - s_{j}) 
$$
Therefore (\ref{11.19.00.1}) plus  a similar expression for the type iii) terms equals to
\begin{equation} \label{11.18.00.2}
\frac{(\ref{11.13.00.5})}{s_j - t_i} \otimes 
\Bigl(\sum_{i=a}^{b-1} (-1)^{i-a}X_{a\to b}^{t_i}\cdot
{\widehat {\rm Li}}(x^{-1}_{i-1}, ..., x^{-1}_a| t_{i}-t_{i-1}, ...,
t_{i} - t_{a}) \cdot
\end{equation}
$$ {\widehat {\rm Li}}(x_{i+1}, ...,  x_{b-1}| t_{i+1}- t_{i}, ...,
t_{b-1} - t_{i})\Bigr)\cdot
$$
$$
\Bigl(\sum_{j=c}^{d-1} (-1)^{j-c}Y^{s_j}_{c \to d}\cdot
{\widehat {\rm Li}}(y^{-1}_{j-1}, ..., y^{-1}_c| s_{j}-s_{j-1}, ...,
s_{j} - s_{c}) \cdot
$$
$$ {\widehat {\rm Li}}(y_{j+1}, ...,  y_{d-1}| s_{j+1}- s_{j}, ...,
s_{d-1} - s_{j}) \Bigr)
$$
The remaining two cases, when the labels are of types ii) or iv), are 
handled the same way.





The computation of the contribution in the case when $a=b$ or $c=d$ is trivial. 

Comparing this with 
$$
\Delta {\widehat {\rm Li}}(*, x_1, ..., x_m| t_0: t_1: ... : t_m) \cdot 
\Delta {\widehat {\rm Li}}(*, y_1, ..., y_n| t_0: s_1: ... : s_n) 
$$
we get the same result. 
Indeed, this expression lives in $Sh_{\bullet}^{I}(G) \otimes Sh_{\bullet}^{I}(G)$, 
and the left factor is  
$$
\sum {\widehat {\rm Li}}(X_{i_0 \to i_1}, ..., X_{i_{k-1} \to i_k}| t_{j_0}: ... : t_{j_{k-1}}) \cdot 
{\widehat {\rm Li}}(Y_{i'_0 \to i'_1}, ..., 
Y_{i'_{l-1} \to i'_l}| s_{j'_0}: ... : s_{j'_{l-1}}) 
$$
where the sum is over all sequences $i_{\alpha}, j_{\alpha}$ and 
$i'_{\beta}, j'_{\beta}$ satisfying the condition (\ref{ur10*}). 
Let us write this as a sum over generalized shuffles, and take a term in the sum we get. 
Then the arcs on the corresponding segment are labelled either by 
\begin{equation} \label{12.10.00.116}
(X_{a \to b}\cdot Y_{c\to d}|t_{i}) \quad \mbox{or} \quad  (X_{a \to b}\cdot Y_{c\to d}|s_{j}) 
 \quad \mbox{where} \quad  
a \leq i < b, \quad c \leq j < d
\end{equation} 
or by $(X_{a \to b}| t_i)$ or $(Y_{c\to d}|s_{j})$.  
In the first case in (\ref{12.10.00.116}) the corresponding term in 
$Sh_{\bullet}^{\rm Li}(G) \otimes Sh_{\bullet}^{\rm Li}(G)$ 
is given by (\ref{11.18.00.2}). The matching pattern 
in the second case in (\ref{12.10.00.116}) 
is completely similar. The other two cases are trivial. 
The theorem is proved. 

\begin{theorem}  $Sh_{\bullet}^{\rm Li}(G)$ is a well defined commutative algebra 
with the coproduct $\Delta$ 
and the product given by the ${\rm Li}$-shuffle formula.  
\end{theorem}

The key ingredient of the proof, the Hopf axiom, 
 is given by theorem \ref{ur100111}. 
To finish the proof it remains  to check 
that $\Delta$ is coassociative, i.e.  $\Delta\circ \Delta =0$. We leave the details 
to the reader.

{\bf 3. Another proof of theorem \ref{4.17.01.77z}}. Theorem \ref{ur100111} 
implies  that 
applying the restricted coproduct $\Delta'$ to an element presenting the left hand side 
of the weight $w$ ${\rm Li}$-shuffle or 
inversion relation we get zero modulo the ${\rm Li}$-shuffle and  
inversion relations of weights $\leq w-1$, plus the obvious additivity relation 1.   
This plus lemma \ref{1.4.02.2} immediately imply that in the Hodge or \'etale setting 
the corresponding element is a constant. 
Let us taking specialization along a curve from the generic point to 
to $x_1 = ... = x_m = 0$. 
We claim that the first column of the period matrix tends to the column $(1, 0, ..., 0)$ 
when $x_i \to 0$. Indeed, it follows from the given in chapter 5 of [G7] 
explicit description 
of the variation of Hodge-Tate structures related to multiple polylogarithms 
that the weight $w$ entry of the first column of the period matrix is given 
by a function of the shape 
$$
{\rm Li}_{k_1, k_2, ...}(x_1 \cdot ... \cdot x_{i_1}, x_{i_1+1} \cdot ...\cdot 
 x_{i_1+i_2}, ...) \qquad w = 
k_1 + ... + k_m
$$ 
Moreover we may assume all of them are given by the convergent power series. Thus 
if $w>0$ each of these functions obviously goes to zero as $x_i \to 0$. 
Thus the equivalence class of the framed object obtained after the specialization 
is zero. Therefore the constant  is zero. 
So the theorem follows. 

There is  a natural homomorphism of 
Hopf $\Q$-algebras 
$
Sh_{\bullet}^{\rm Li}(\C^*)\lra {\cal H}_{\bullet}
$. 

{\bf 4. Reamrks on an integral version of the inversion formula}. Our definition of 
$\widehat B(x|t)$ was motivated by 
Kronecker's formula 
$$
B(\varphi| t) := \sum_{-\infty <k < \infty}
\frac{e^{2\pi i (\varphi k)}}{(k - t)} =  -2\pi i
\frac{e^{2\pi i \{\varphi\}t}}{e^{2\pi i t}-1} =  -2\pi i \cdot \sum_{n \geq 0} B_n(\{\varphi\}) 
\frac{(2 \pi i t)^{n-1}}{n!}
$$
Here $\{\varphi\}$ is the fractional part of $\varphi $.  However this formula 
 and $\widehat B(x|t):= x^t/t$ agree only if we  suppress all the Bernoulli numbers 
except $B_0=1$ from Kronecker's formula. The discussion below shows why this definition 
of $\widehat B(x|t)$ is meaningful if we work modulo torsion, and how one can try to make 
it integral.

The classical formula for the Bernoulli numbers
$$
B_{2n} := -(2\pi i)^{-2n}2 (2n)! \zeta(2n)
$$
suggests the following definitions
$$
B^{\cal M}_{2n} := -2 (2n)! \zeta^{\cal M}(2n); \qquad B^{\cal M}_{1} := 
\log^{\cal M}(-1), \quad 
B^{\cal M}_{0} := 1
$$
$$
  {B}^{\cal M}_n({\log}^{\cal M} (x)) := \sum_{k=0}^n  {n \choose k} B^{\cal M}_k 
({\log}^{\cal M}(x))^{n-k} 
$$
$$
B^{\cal M}(x|t):=  \sum_{n \geq 0} {B}^{\cal M}_n({\log}^{\cal M} (x)) 
\frac{t^{n-1}}{n!} 
$$  
Here the objects we introduce belong to the abelian group of framed mixed 
Tate motives. This group is  supposed to be defined {\it integrally}. 
We are not discussing below how to define this group. Its
 l-adic realization should be given by Galois  $\Z_l$-modules. 
Its Hodge version is given by the equivalence classes of integral Hodge-tate structures. 
So ``the motivic Bernoulli numbers'' ${B}^{\cal M}_n$ are no longer 
numbers but rather weight $n$ framed mixed Tate motives,  
which are torsion elements for $n>0$. 

We claim that one has the inversion formula
$$
{\rm Li}^{\cal M}(x|t) - {\rm Li}^{\cal M}(x^{-1}|-t) = -B^{\cal M}(x|t)
$$

For example 
$$
{\rm Li}^{\cal M}_1(x) - {\rm Li}^{\cal M}_1(x^{-1})  = - \log^{\cal M}_1(1-x) + 
\log^{\cal M}_1(1-x^{-1}) = 
$$
$$
-\log^{\cal M}(x)  +  \log^{\cal M}_1(-1) = -{B}^{\cal M}_1({\log}^{\cal M} (x))
$$
Further
$$
{\rm Li}^{\cal M}_2(x) + {\rm Li}^{\cal M}_2(x^{-1})  \quad =  \quad 
-\frac{1}{2}\Bigl((\log^{\cal M}(x))^2 + 2 \log^{\cal M}(-1)\cdot \log^{\cal M}(x) + {B}^{\cal M}_2\Bigr) = 
$$
$$
-\frac{1}{2}\Bigl((\log^{\cal M}(x))^2 + {B}^{\cal M}_2\Bigr)
$$
since $2 \log^{\cal M}(-1)=0$, and so on. 

The formula we want to demonstrate has been 
checked in the example above  for $n=1$. At $x=1$ this formula 
is trivial for odd $n>1$,   and 
boils down to identity 
$$
{\rm Li}^{\cal M}_{2n}(1) + {\rm Li}^{\cal M}_{2n}(1)  \quad =  \quad 
-\frac{1}{(2n)!} {B}^{\cal M}_{2n}
$$
which is clear by the very definition 
of the Bernoulli framed motives ${B}^{\cal M}_k$. 
Thus it remains to check that it is killed by the restricted coproduct $\Delta'$, 
which we do by induction using the following formula (see [G7]):
$$
\Delta' {\rm Li}^{\cal M}(x|t) = {\rm Li}^{\cal M}(x|t) \otimes (x^t-1)
$$
Namely, 
$$
\Delta' \Bigl({\rm Li}^{\cal M}(x|t) - {\rm Li}^{\cal M}(x^{-1}|-t) \Bigr) 
= \Bigl({\rm Li}^{\cal M}(x|t) - {\rm Li}^{\cal M}(x^{-1}|-t)\Bigr)\otimes (x^t-1) 
$$
so using the induction assumption we can rewrite this as 
$ 
-B^{\cal M}(x|t) \otimes (x^t-1)
$. 
It remains to notice that observing $\Delta'(B^{\cal M}_{n}) =0$ it is easy to check that 
$$
\Delta' (-B^{\cal M}(x|t)) = -B^{\cal M}(x|t) \otimes (x^t-1)
$$
 The theorem is proved. 

The following result provides the  integral motivic version of the 
inversion formula (\ref{12.3.00.1}).

\begin{theorem} \label{}
The  
inversion formula (\ref{12.3.00.1}), where 
$\widehat {\rm Li}$ is changed to  ${\rm Li}^{\cal M}$, and 
$\widehat B(x|t):= B^{\cal M}(x|t)$, holds 
for the depth $m$ motivic 
multiple polylogarithms. 
\end{theorem}

There are three different proofs of this result: 

1) Argue as in the proof of the previous theorem. 

2) Take the proof of proposition 2.7  in section 2.6 in [G7] and make it 
motivic following the same ideas as in chapter 9.

3) Follow the proof of theorem 4.1 in [G4] and  use the motivic 
double shuffle relations.

\vskip 3mm \noindent
{\bf REFERENCES}
\begin{itemize}
\item[{[BBD]}] Beilinson A.A., Bernstein J., Deligne P.: {\it
    Faisceaux pervers}.  Analysis and topology on singular spaces, I (Luminy,
1981), 5--171, Ast{\'e}risque, 100, Soc. Math. France, Paris, 1982.
  \item[{[BD1]}] Beilinson A.A., Deligne P.: {\it Interpr{\'e}tation motivique de la conjecture de Zagier reliant polylogarithmes et r{\'e}gulateurs}. Motives (Seattle, WA, 1991), 97--121, 
Proc. Sympos. Pure Math., 55, Part 2,
AMS, Providence, RI, 1994.
\item[{[BD2]}] Beilinson A.A., Deligne P.: 
{\it Motivic polylogarithms and Zagier's conjecture} 
Manuscript, version of 1992. 
\item[{[BMS]}] Beilinson A.A., MacPherson R.D. Schechtman V.V: {\it
Notes on motivic cohomology}. Duke Math. J., 1987 vol 55 
\item[{[BGSV]}] Beilinson A.A., Goncharov A.A., Schechtman V.V., Varchenko A.N.: {\it Aomoto dilogarithms, mixed Hodge structures and motivic cohomology of a pair of triangles in the plane}, the Grothendieck Feschtrift, Birkhauser, vol 86, 1990, p. 135-171.
\item[{[Be]}] Bernstein J., {\it Lectures on D-modules}. 
\item[{[Br]}] Broadhurst D.J., {\it On the enumeration of irreducible $k$-fold sums and their role in knot theory and field theory} Preprint hep-th/9604128. 
\item[{[D]}] Deligne P.: {\it Le group fondamental de la droite projective moine trois points}, In: Galois groups over $\Q$. 
Publ. MSRI, no. 16 (1989) 79-298.   
\item[{[DG]}] Deligne P., Goncharov A.B.:  In preparation. 
\item[{[D2]}] Deligne P.: {\it A letter to D. Broadhurst}, June 1997.
\item[{[D3]}] Deligne P.: {\it Letter to the author}. July 25, 2000.
\item[{[D4]}] Deligne P.: {\it Th{\'e}orie de Hodge II, III}. Publ. IHES No. 40 (1971), 5--57;
No. 44 (1974), 5--77.
\item[{[DJ]}] De Jeu, Rob: {\it Zagier's conjecture and wedge complexes in algebraic $K$-theory}. Compositio Math. 96 (1995), no. 2,
197--247.
\item[{[Dr]}] Drinfeld V.G.: {\it On quasi-triangular quasi-Hopf algebras and some group related to clousely associated with Gal$(\overline{\Q}/\Q)$}. Leningrad Math. Journal, 1991. (In Russian).
\item[{[E]}] J.Ecalle, {\it ARI/GARI and the arithmetics of 
multizetas: the main aspects}. To appear. 
\item[{[G0]}] Goncharov A.B.: {\it Multiple $\zeta$-numbers,
    hyperlogarithms and mixed Tate motives}, Preprint MSRI 058-93, June 1993.
\item[{[G1]}] Goncharov A.B.: {\it Polylogarithms in arithmetic and geometry}, 
Proc. of the International Congress of Mathematicians, Vol. 1, 2
(Zurich, 1994), 374--387, Birkhauser, Basel, 1995.
\item[{[G2]}] Goncharov A.B.: {\it The double logarithm and Manin's complex for modular curves}. Math. Res. Lett. 4 (1997), no. 5, 617--636.
\item[{[G3]}] Goncharov A.B.: {\it Multiple polylogarithms, cyclotomy and modular complexes
},     Math. Res. Letters, 
 vol. 5. (1998), pp. 497-516. www.math.uiuc.edu/K-theory/ N 297.
\item[{[G4]}] Goncharov A.B.: {\it The dihedral Lie algebras and galois symmetries of 
$\pi_1^{(l)}({\Bbb P}^1 \backslash 0, \mu_N, \infty)$}. Duke Math. J.  vol 100, N3, (2001), pp. 397-487. 
math.AG/0009121. 
\item[{[G5]}] Goncharov A.B.: {\it Multiple $\zeta$-values, Galois groups and geometry of 
modular varieties} Proc. of the Third European Congress of
Mathematicians. Progress in Mathematics, 
vol. 201, p. 361-392. Birkhauser Verlag. (2001)
math.AG/0005069.
\item[{[G6]}] Goncharov A.B.: {\it Geometry of configurations, Polylogarithms and motivic cohomology}.  Adv. in Math.  114 (1995), no. 2, 197--318.
\item[{[G7]}] Goncharov A.B.: {\it Multiple polylogarithms and mixed
    Tate motives}. math.AG/0103059. 
\item[{[G8]}] Goncharov A.B.: {\it Mixed ellitptic motives}. Galois representations in arithmetic algebraic geometry (Durham, 1996), 147--221, London
Math. Soc. Lecture Note Ser., 254, Cambridge Univ. Press, Cambridge, 1998. 
K-theory e-print archive, www.math.uiuc.edu/K-theory/ N 228.
\item[{[G9]}] Goncharov A.B.: {\it Volumes of hyperbolic manifolds 
and mixed Tate motives} J. Amer. Math. Soc. 12 (1999) N2, 569-618. 
math.alg-geom/9601021
\item[{[G10]}] Goncharov A.B.: {\it Multiple L-values, 
geometry of symmetric spaces and Feynman diagrams}. To appear. 
\item[{[GM]}] Goncharov A.B., Manin Yu.I.: {\it 
Multiple $\zeta$--motives and the moduli space ${\cal M}_{0,n}$}. 
 math.AG/0204102. 
\item[{[H1]}] Hain R., 
{\it The geometry of mixed Hodge structure on $\pi_1$}.  
 Algebraic geometry, Bowdoin, 1985
(Brunswick, Maine, 1985), 247--282, Proc. Sympos. Pure Math., 46, Part 2, Amer. Math. Soc., Providence, RI, 1987. 
\item[{[HM]}] Hain R., Matsumoto M..: {\it Weighted Completion of
    Galois Groups and Some 
Conjectures of Deligne}. math.AG/0006158. 
\item[{[Ih]}] Ihara Y.: {\it Braids, Galois groups, and some arithmetic functions}. 
Proc. of the Int. Congress of Mathematicians, Vol. I, II
(Kyoto, 1990), 99--120, Math. Soc. Japan, Tokyo, 1991.
\item[{[IK]}] Ihara K.,  Kaneko M.: {\it Derivation relations and regularized double shuffle relations for multiple zeta values}. Preprint 2001.
\item[{[L]}] Levine, M: {\it Mixed motives}. Mathematical Surveys and
  Monographs, 57. American 
Mathematical Society, Providence, RI, 1998
\item[{[L1]}] Levine, M: {\it Tate motives and the vanishing conjectures for algebraic 
$K$-theory}. Algebraic $K$-theory and algebraic topology (Lake
Louise, AB, 1991), 167--188, NATO Adv. Sci. Inst. Ser. C 
Math. Phys. Sci., 407, 
Kluwer Acad. Publ., Dordrecht, 1993.
\item[{[M]}] Milne J.S. {\it \'Etale cohomology} Princeton University Press,
Princeton, NJ, 1980.
\item[{[R]}] Racinet G.: {\it S\'eries g\'en\'eratrices non commutatives de polyzetas et associateurs de Drinfel'd}. Th\'ese, (2000). 
\item[{[R1]}] Racinet G.: {\it Doubles melanges des polylogarithmes multiples aux racines de l'unit'e}. math.QA/0202142.
\item[{[T]}] Terasoma T.: {\it Multiple zeta values and mixed tate
    motives}. math.AG/0104231. 
\item[{[Z]}] Zagier D.: {\it Values of zeta functions and their
    applications}. First European 
Congress of Mathematics, Vol. II (Paris, 1992), 497--512, Progr.
Math., 120, Birkhauser, Basel, 1994. 
\item[{[Ve]}] Verdier, J.-L. 
Sp{\'e}cialisation de faisceaux et monodromie mod{\'e}r{\'e}e. 
Analysis and topology on singular spaces, II, III (Luminy, 1981), 332-364, 
Ast{\'e}risque, 101-102, 
Soc. Math. France, Paris, 1983. 
\item[{[V]}] Voevodsky V.:  {\it Triangulated 
category of motives over a field} in 
Cycles, transfers, and motivic homology 
theories. Annals of Mathematics Studies,
143. Princeton University Press, Princeton, NJ, 2000. 
\item[{[Vol]}] Vologodsky V.:  
{\it Hodge structure on the fundamental group and its application to p-adic integration}. math.AG/0108109. 
\item[{[W]}] Woitkowiyak Z.: {\it Cosimplicial objects in algebraic geometry}. 
Algebraic $K$-theory and algebraic topology (Lake Louise, AB, 1991),
287--327, NATO Adv. Sci. Inst. Ser. C Math. Phys. Sci., 407, 
Kluwer Acad. Publ., 
Dordrecht, 1993.

\end{itemize}

\end{document}